%% file: 2degreegeneral-revision.tex
\def\version{07 June 2021}	        	%

\documentclass[reqno,11pt]{amsart}

\usepackage[T1]{fontenc}
\usepackage[utf8]{inputenc}
\usepackage{amsmath} 
\usepackage[mathscr]{eucal}
\usepackage{amssymb,bbm}
\usepackage{srcltx} 
\usepackage{dsfont}
\usepackage{hyperref}
\usepackage{color}
\usepackage{enumerate}
\usepackage{tikz}
\usepackage{boldline}
\usepackage{comment}
\usepackage{lscape}
\usepackage{graphicx}
\usepackage{subcaption}

\numberwithin{equation}{section}
 
\def\la{\lambda}

\def\emptyset{\varnothing} 

\def\I{\tau}

\def\d{{\rm d}} 
\def\e{\varepsilon} 
 
\def\L{\Lambda}

\newfam\Bbbfam 
\font\tenBbb=msbm10 
\font\sevenBbb=msbm7 
\font\fiveBbb=msbm5 
\textfont\Bbbfam=\tenBbb 
\scriptfont\Bbbfam=\sevenBbb 
\scriptscriptfont\Bbbfam=\fiveBbb

\def\one{\mathds{1}}
\def\kNN{$k$NN}
\def\U{$\mathbf U$}
\def\B{$\mathbf B$}
\def\NN{\mathrm{v}}
\def\NNN{\mathbf{v}}

\def\bs{\boldsymbol}

\def\X{{\mathbf{X}}}
\def\ZZ{{\mathbf{z}}}

\newcommand{\R}     {\mathbb{R}} 
\newcommand{\Z}     {\mathbb{Z}} 
\newcommand{\N}     {\mathbb{N}} 
 
\renewcommand{\P}   {\mathbb{P}} 
 
\newcommand{\E}     {\mathbb{E}} 
\newcommand{\Q}     {\mathbb{Q}}

\def\1{{\mathchoice {1\mskip-4mu\mathrm l}      
{1\mskip-4mu\mathrm l} 
{1\mskip-4.5mu\mathrm l} {1\mskip-5mu\mathrm l}}} 
 
\def\comment#1{} 
\newtheoremstyle{thm}{2ex}{2ex}{\itshape\rmfamily}{} 
{\bfseries\rmfamily}{}{1.7ex}{} 
 
\newtheoremstyle{rem}{1.3ex}{1.3ex}{\rmfamily}{} 
{\itshape\rmfamily}{}{1.5ex}{}

\renewcommand{\theequation}{\thesection.\arabic{equation}} 
 
\newtheorem{theorem}{Theorem}[section] 
\newtheorem{lemma}[theorem]{Lemma} 
\newtheorem{prop}[theorem] {Proposition}

\theoremstyle{definition}
\newtheorem{defn}[theorem] {Definition} 
\newtheorem{example}[theorem] {Example}

\renewcommand{\d}{{\rm d}} 
 
\newcommand{\eps}{\varepsilon}

\def\bs{\boldsymbol}

\newcommand{\Bcal}  {{\mathcal B}}

\newcommand{\Mcal}   {{\mathcal M }}

\newcommand\numberthis{\addtocounter{equation}{1}\tag{\theequation}}
\renewcommand{\e}   {{\operatorname e }}

\definecolor{Red}{rgb}{1,0,0}

 \title[Absence of percolation in stationary graphs with degrees bounded by two]{Absence of percolation in graphs based on stationary point processes with degrees bounded by two}
 \author[Benedikt Jahnel and Andr\'as T\'obi\'as]{}
 
\begin{document}
 \maketitle
 \centerline{{\sc Benedikt Jahnel\footnote{Weierstrass Institute Berlin, Mohrenstra\ss e 39, 10117 Berlin, Germany, \texttt{Jahnel@wias-berlin.de}} and Andr\'as T\'obi\'as\footnote{TU Berlin, Stra\ss e des 17.~Juni 136, 10623 Berlin, Germany, \texttt{Tobias@math.tu-berlin.de}}}}

\renewcommand{\thefootnote}{}

\bigskip

\centerline{\small(\version)} 
\vspace{.5cm} 
 
\begin{quote} 
{\small {\bf Abstract:}} We consider undirected graphs that arise as deterministic functions of stationary point processes such that each point has degree bounded by two. For a large class of point processes and edge-drawing rules, we show that the arising graph has no infinite connected component, almost surely. In particular, this extends our previous result for signal-to-interference ratio graphs based on stabilizing Cox point processes and verifies the conjecture of Balister and Bollob\'as that the bidirectional $k$-nearest neighbor graph of a two-dimensional homogeneous Poisson point process does not percolate for $k=2$.
\end{quote}




\setcounter{tocdepth}{3}

\section{Introduction}\label{sec-Intro}
\input{2degree-Intro.tex}

\section{Model definition and main result}\label{sec-model}
\input{2degree-Model.tex}

\section{Proof of Theorem~\ref{theorem-nopercolation}}\label{sec-proof}
\input{2degree-Proof.tex}


\section{Examples, discussion and extensions}\label{sec-examples}
\input{2degree-Examples.tex}

\section*{Acknowledgements}
The authors thank A.~Hinsen and C.~Hirsch for interesting discussions and comments and D.~Dereudre for pointing out the reference~\cite{S18}. The authors thank two anonymous reviewers for useful comments on an earlier version of the manuscript, in particular for pointing out the reference~\cite{HS13}. Further, the second author thanks a number of participants of the spring school \emph{Complex Networks} at TU Darmstadt for their insightful questions that enlightened that our methods for SINR graphs can be lifted to the \B-2NN graph. This work was funded by the Deutsche Forschungsgemeinschaft (DFG, German Research Foundation) under Germany's Excellence Strategy MATH+: The Berlin Mathematics Research Center, EXC-2046/1 project ID: 390685689.

\begin{comment}{
\newpage
%
Checklist:
\begin{itemize}
    \item Applications of Gibbs p.p.\ for telecommunications? We could speak about my Gibbsian model with Wolfgang and say something about the thermodynamical limit...
    \item Need to decide if we want to mention the fact that the reduced Palm measure P0! is absolutely continuous with respect to P.
    \item Delete this list in the end
\end{itemize}}\end{comment}
\end{document}

%% file: 2degree-Intro.tex
Continuum percolation was introduced by Gilbert~\cite{G61} in order to model connectivity in large telecommunication networks. In his graph model, the vertices form a homogeneous Poisson point process (PPP) in $\R^2$, and two points are connected whenever their distance is less than a fixed connection radius $r>0$. He showed that this model undergoes a phase transition: if the spatial intensity $\lambda>0$ of the PPP is sufficiently small, then the graph consists of finite components only, almost surely, whereas for large enough $\lambda$, the graph \emph{percolates}, i.e., it has an unbounded connected component, also almost surely.

\begin{figure}[!htpb]
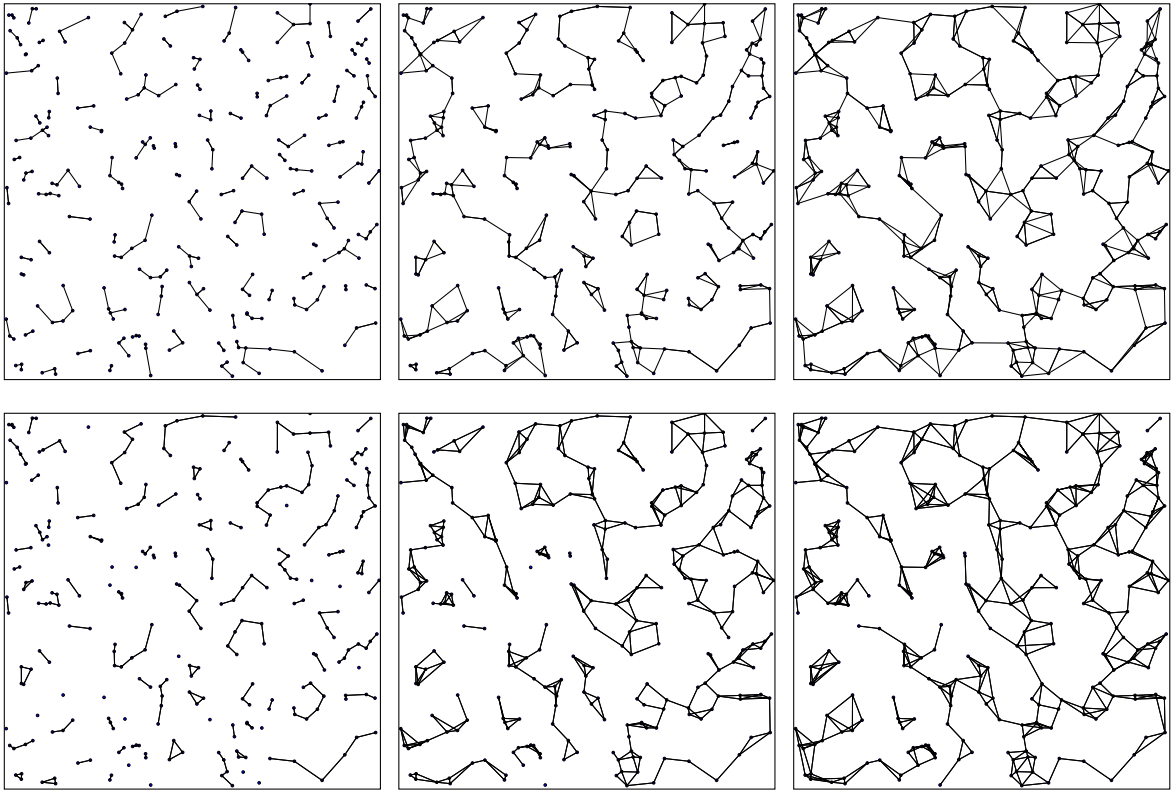

\begin{subfigure}{0.32\textwidth}
\input{Fig-U1NN.tex}
\end{subfigure}
\begin{subfigure}{0.32\textwidth}
\input{Fig-U2NN.tex}
\end{subfigure}
\begin{subfigure}{0.32\textwidth}
\input{Fig-U3NN.tex}
\end{subfigure}
\\\vspace{0.4cm}
\begin{subfigure}{0.32\textwidth}
\input{Fig-B2NN.tex}
\end{subfigure}
\begin{subfigure}{0.32\textwidth}
\input{Fig-B4NN.tex}
\end{subfigure}
\begin{subfigure}{0.32\textwidth}
\input{Fig-B5NN.tex}
\end{subfigure}
\caption{Top: Realizations of \U-\kNN~graphs based on a PPP with $k=1$ (left), $k=2$ (center) and $k=3$ (right). Bottom: Realizations of \B-\kNN~graphs based on a PPP with $k=2$ (left), $k=4$ (center) and $k=5$ (right).}
\label{Fig-kNN}
\end{figure}

This model has been widely extended, for instance to the case of random connection radii and for various point processes, see~\cite{MR96,BY13,CD14,GKP16,HJC17,J16,S13,JTC20}. A drawback of Gilbert's model is that it allows for an arbitrarily large degree of the vertices, whereas for many applications, it is a reasonable assumption that the vertices should have bounded degree. Incorporating this property, Häggström and Meester \cite{HM96} studied percolation in the so-called \emph{undirected $k$-nearest neighbor} (\U-\kNN) \emph{graph}, based on a stationary PPP in $\R^d$, $d \geq 1$, see top line of Figure~\ref{Fig-kNN}. Here, all points of the point process are connected to their $k$-nearest neighbors, for some fixed $k \in \N$. This results in a graph that is the undirected variant of a directed graph with out-degrees bounded by $k$, which itself also has degrees larger than $k$. Let us write $k_{\mathbf U,d}$ for the minimum of all $k\in \N$ such that the \U-\kNN-graph of the stationary PPP in $\R^d$ percolates with positive probability. It was shown in \cite{HM96} that $k_{\mathbf U,d} >1$ for all $d \in\N$, however, $k_{\mathbf U,d}=2$ for all sufficiently large $d$. This was complemented in \cite{TY07} by the assertion that $k_{\mathbf U,d}<\infty$ for all $d \geq 2$.

Balister and Bollobás \cite{BB08} studied the case $d=2$. They also introduced another undirected graph, which is contained in the \U-\kNN~graph, called the \emph{bidirectional $k$-nearest neighbor} (\B-\kNN) \emph{graph}, see bottom line of Figure~\ref{Fig-kNN}. Here, one connects two points of the point process if and only if they are mutually among the $k$-nearest neighbors of each other. This graph has in fact degrees bounded by $k$, which immediately implies that there is no percolation for $k=1$, whatever the vertex set is (note that in the PPP case this also follows from the results of \cite{HM96}). Define the critical out-degree $k_{\mathbf B,d}$ analogously to $k_{\mathbf U,d}$ but with \U~replaced by \B.  It was shown in \cite{BB08} that $k_{\mathbf U,2} \leq 11$ and $k_{\mathbf B,2} \leq 15$. Further, `high-confidence results' of \cite{BB08} indicate $k_{\mathbf U,2}=3$ and $k_{\mathbf B,2}=5$. By `high-confidence results', the authors of that paper meant that these assertions follow once one shows that a certain deterministic integral exceeds a certain deterministic value, however, the integrals were only evaluated via Monte--Carlo methods so far. Hence, from a mathematical point of view, these are still open conjectures, but there is very strong numerical evidence that they are true.

Another line of research on percolation of bounded-degree spatial graphs with unbounded-range dependences, which is also close to applications in wireless networks, is \emph{signal-to-interference plus noise ratio (SINR)} percolation, introduced in \cite{DBT05,DF06}. Here, a transmission in the network is considered successful if and only if, measured at the receiver, the incoming signal power of the transmitter is larger than a given threshold times the interference (sum of signal powers) coming from all other users plus some external noise. Then, the SINR graph is constructed by drawing an edge between two vertices whenever the transmission between them is successful in both directions, see Section~\ref{sec-examplesf} for more details. This graph has bounded degrees  (see \cite[Theorem 1]{DBT05}), where the smallest degree bound $k$ depends on the model parameters. If the transmitted signal powers are all equal, then the SINR graph is contained in the corresponding \B-\kNN~graph (see \cite[Lemma 4.1.13]{T19}) and hence also of the \U-\kNN~graph.

In general, if in an undirected graph all degrees are bounded by $k=2$, all infinite connected components must be path graphs (no cycles, no branchings), infinite in one or two directions, which makes the graph similar to a one-dimensional continuum percolation model, indicating that under rather general conditions, there should be no infinite connected component. Certainly, there are deterministic point processes where percolation is possible, but a little bit of randomness can be expected to suffice for non-percolation. In our recent paper \cite{JT19}, we showed that in SINR graphs based on general stationary Cox point processes (CPPs) in any dimension, under rather general choices of the parameters resulting in degrees bounded by 2, there is no percolation.

In the present paper, moving away from the particular setting of SINR graphs, we present analogous results in a general framework, extending the methods of the proof of \cite[Theorem 2.2]{JT19}. We consider a generalization of the \B-\kNN~graph, called the $f$-\kNN~graph. Here, points of the underlying marked point process are connected by an edge whenever they are mutually among the $k$ nearest neighbors of each other with respect to an ordering that may also depend on some marks of the points, apart from the (not necessarily Euclidean) distance of the points. The ordering is expressed in terms of a function $f$, hence the name $f$-\kNN~graph. We show that under suitable conditions on the underlying stationary marked point process, the $f$-\kNN~graph does not percolate for $k=2$. This in particular implies non-percolation of subgraphs of the $f$-\kNN~graph depending on additional randomness. In fact, our results extend to all stationary point processes that are \emph{deletion-tolerant} in the sense of \cite{HS13}. This includes all CPPs satisfying a basic nondegeneracy condition, and 
a large class of Gibbs point processes. 

As a special case, our results imply that for general stationary CPPs, the \B-2NN~graph does not percolate. This in particular implies that $k_{\mathbf B,2} \geq 3$, which provides a partial verification of the high-confidence results of \cite{BB08}. Note that this result does not follow from \cite[Theorem 2.2]{JT19} because in general, if the SINR graph is contained in a \B-$2$NN~graph, it is a proper subgraph of it with substantially less edges. After stating and proving our main results, we also present examples of graphs with degrees bounded by 2 that are not contained in an $f$-$2$NN~graph but where our proof techniques are also applicable, and also ones where they are not applicable.

Our setting is also related to the line of research on \emph{outdegree-one} graphs, which were introduced in \cite{CDS20}. However, our results do not follow from the results of that paper, and also not the other way around. We will comment on the similarities and differences of the two models in Section~\ref{sec-outdegree1}. 

The rest of the paper is organized as follows. In Section~\ref{sec-model} we present our setting and main result. In Section~\ref{sec-proof} we provide the proofs for our main result. Section~\ref{sec-examples} is devoted to examples, extensions of our methods, and discussions.

%% file: 2degree-Model.tex
In this section we present our model definition and main results. Our setting is as follows. Let $d \in \N$, and let $\Vert \cdot \Vert$ be an arbitrary norm on $\R^d$. Further, let $\mathcal B(\R^d)$ denote the Borel-$\sigma$-algebra of $\R^d$ (clearly, $\Vert \cdot \Vert$ generates the standard topology of $\R^d$, which generates $\Bcal(\R^d)$). Moreover, consider the measurable space $(M,\Mcal)$, which serves as a mark space.

Next, let $\mathbf X=\{(X_i,P_i)\}_{i \in I}$ be a marked point process in $\R^d \times M$, so that $X=\{X_i\}_{i \in I}$ is a stationary point process in $\R^d$ with finite intensity $\lambda=\E[X([0,1]^d)]$, that is \emph{nonequidistant}. This means that for all $i,j,k,l\in I$, $\Vert X_i-X_j\Vert=\Vert X_k-X_l\Vert>0$ implies $\{ i,j\} = \{ k,l \}$ and $\Vert X_i\Vert=\Vert X_j\Vert$ implies $i= j$, almost surely.  Clearly, this property implies that the point process $X$ is \emph{simple}, i.e., $\P(X_i\neq X_j,\forall i,j \in I\text{ with }i \neq j)=1$. 
For illustration, note that the randomly shifted lattice $\Z^d+U$, where $U$ is a uniform random variable in $[0,1]^d$, is a simple stationary but not nonequidistant point process on $\R^d$.  

Next, we introduce a total ordering of the points. For this, let $f \colon [0,\infty)\times M \to [0,\infty)$ be a measurable function such that $v \mapsto f(v,q)$ is monotone decreasing for all $q \in M$. We call such a function an {\em ordering function}.  
\begin{defn}\label{defn-closer}
Let $f$ be an ordering function and $(v_1,q_1),(v_2,q_2),(v_3,q_3)\in\R^d \times M$. We say that \emph{$v_2$ is $f$-closer to $v_1$ than to $v_3$} if one of the following conditions is satisfied:
\begin{enumerate}
    \item $f(\Vert v_1-v_2 \Vert,q_1) < f(\Vert v_3-v_2 \Vert,q_3)$, or
    \item\label{tiebreaking} $f(\Vert v_1-v_2 \Vert,q_1) = f(\Vert v_3-v_2 \Vert,q_3)$ and $\Vert v_1-v_2\Vert<\Vert v_3-v_2\Vert$.
\end{enumerate}
\end{defn}
Then, it is elementary to verify the following lemma.
\begin{lemma}
Let $f$ be an ordering function. For $\mathbf X$ defined as above and $X$ nonequidistant, almost surely, the following holds. For all $i \in I$, the relation ``$X_i$ is $f$-closer to $X_j$ than to $X_l$'' is a total ordering (i.e., irreflexive, antisymmetric and transitive, with any two elements being comparable) on the set $\{ (j,l) \in I\times I \colon j \neq i \text{ and } l \neq i \}$, which we call the $f$-ordering.
\end{lemma}
Thus, if $\bs x=\{(x_i,p_i)\}_{i \in I}$ is a deterministic, locally finite, infinite, and nonequidistant set of points in $\R^d \times M$ (for some countable index set $I$) and $v_o\in x:=\{x_i\}_{i \in I}$, we can represent $x$ as $x= \{\NN^f_n(v_o,\bs x)\}_{n\in \N_0}$, where $\NN^f_0(v_o,\bs x)=v_o$, and $\NN^f_n(v_o,\bs x)$ is the $n$-th nearest neighbor of $v_o$ in $x$ with respect to the $f$-ordering for any $n\in\N_0$. 
Next, we build a graph based on the $f$-ordering. 
\begin{defn}\label{defn-fkNN} Let $f$ be an ordering function, $k \in \N$ and $\mathbf X$ defined as above with $X$ nonequidistant, almost surely. The \emph{$f$-$k$-nearest neighbor} ($f$-\kNN) \emph{graph} $g_{k,f}(\X)$ is the random graph having vertex set $X$ and for all $i \in I$ and $n \in \{ 1,\ldots,k\}$ an edge between $X_i$ and $\NN^f_n(X_i, \X)$ whenever $X_i \in \{ \NN^f_1(\NN^f_n(X_i,\X),\X), \ldots, \NN^f_k(\NN^f_n(X_i,\X),\X) \}$.
\end{defn}
As the next example shows, the \B-\kNN~graph is an $f$-\kNN~graph for a point process with trivial marks. Let us write $\{ \star \}$ for the one-point measurable space (with $\Mcal=\{ \emptyset, \{ \star \} \}$).
\begin{example}
Consider a nonequidistant  point process $X$ in $\R^d$, $d \geq 1$, and equip $X$ with trivial marks in $M=\{ \star \}$. Then, $f(v,q)=f(v)=v$ yields the \B-\kNN~graph based on $X$.
\end{example}
We will explain the relations between $f$-\kNN~graphs and SINR graphs in Section~\ref{sec-examplesf}.
%

Apart from the basic requirement of being nonequidistant, the property of \emph{deletion-tolerance} introduced in \cite{HS13} is the most important requirement on the marked point process. An \emph{$\X$-point} is an $(\R^d\times M)$-valued random variable $\ZZ$, defined on the same probability space as $\X$, such that $\ZZ\in \X$ a.s., and one says that $\X$ is \emph{deletion-tolerant} if for any $\X$-point $\ZZ$, the distribution of $\X \setminus \{ \ZZ \}$ is absolutely continuous with respect to the one of $\X$. See \cite[Theorem 1.1]{HS13} for the equivalent formulations of this property for non-marked point processes. We will present examples of marked point processes that are deletion-tolerant below.  Equipped with the above definitions, we are now able to state our main result.

\begin{theorem}\label{theorem-nopercolation}
Let $f$ be an ordering function and let the deletion-tolerant marked point process $\mathbf X$ be such that the underlying point process $X$ is stationary, nonequidistant and has a finite intensity. Then, we have
\[ \P(g_{2,f}(\mathbf X)\text{ percolates}) = 0. \]
\end{theorem}
The proof of this theorem is carried out in Section~\ref{sec-proof}.  In Section~\ref{sec-examplesf} we discuss the relation between the proof of Theorem~\ref{theorem-nopercolation} and the one of \cite[Theorem 2.2]{JT19}. In Section~\ref{sec-examplesnotf} we will explain how it extends to other graphs that are defined similarly, have degrees bounded by 2, but are not subgraphs of $f$-\kNN~graphs. A key ingredient in the proof of Theorem~\ref{theorem-nopercolation} and its aforementioned extensions is the so-called \emph{edge-preserving property} of the underlying graph, which we will introduce in Definition~\ref{defn-edgepreserving} below.

Note that deletion-tolerance is satisfied by many point processes, as is shown by the following proposition, see also  \cite{HS13} for more examples without marks. 
\begin{prop}\label{prop-thinning}
Any i.i.d.-marked Cox point process $\mathbf X$ on $\R^d\times M$ is deletion-tolerant. Further, any infinite-volume Gibbs point process $\X$ on $\R^d\times M$ based on an Hamiltonian $H$ is deletion-tolerant, whenever $M$ is a Polish space equipped with the associated Borel $\sigma$-algebra, and for all bounded measurable $\L\subset\R^d$ and almost-all boundary conditions ${\boldsymbol y}$
$$Z_\L({\boldsymbol y})=\int \mathcal  P_\L(\d {\boldsymbol x}_\L)\exp(-H_\L({\boldsymbol x}_\L\cup {\boldsymbol y}_{\L^{\rm c}}))<\infty,$$
where $\mathcal P_\L$ is an i.i.d.-marked Poisson point process in $\L$ and $\bs y_{\L^{\rm c}}=\{(v,q)\in \bs y\colon v\in \R^d\setminus \L)$. 
\end{prop}
The proof of this proposition is presented in Section~\ref{sec-thinningexamples}. 
Note that the class of stationarity and nonequidistant CPPs includes the homogeneous PPPs. The class of Gibbs point processes satisfying the above condition is very rich and in particular includes the classical examples of superstable Hamiltonians, see~\cite{D19} and references therein. 
Moreover, there are also well-known point processes that are not deletion-tolerant. A very straightforward class of such processes is the following. As introduced in \cite{GP12,KL20}, we say that the point process $X$ is \emph{number rigid} if for any $\mathcal K \subset \R^d$ compact, there exists a deterministic measurable function $h_\mathcal K$ such that,
\[ \# (X \cap \mathcal K) = h_\mathcal K(X \setminus \mathcal K), \]
almost surely, i.e., $X$ outside $\mathcal K$ determines the number of points of $X$ in $\mathcal K$. The following proposition states that number rigid point processes fail to be deletion-tolerant.
\begin{prop}\label{prop-rigidity}
If the non-marked version $X$ of the marked point process $\mathbf X$ is stationary and number rigid with positive intensity, then $\mathbf X$ is not deletion-tolerant.
\end{prop}
This proposition follows immediately from results of \cite{HS13}, see Section~\ref{sec-thinningexamples} for a proof.

According to \cite{GP12}, the Ginibre ensemble and the Gaussian zero process are number rigid point processes in $\R^2$, which are also stationary, nonequidistant, and of positive intensity. Hence, the proof of Theorem~\ref{theorem-nopercolation} is not applicable for them. We nevertheless believe that they satisfy the assertion of the theorem, but the proof would require additional arguments.

%% file: 2degree-Proof.tex
The proof of the Theorem~\ref{theorem-nopercolation} proceeds along the following line of arguments. We first show that with probability 1, clusters, i.e., maximally connected components, are either finite or they consist only of points of degree 2, see Lemma~\ref{lemma-noheads} below. Next, we assume for a contradiction that there exists an infinite cluster with positive probability. Then, we introduce a procedure that removes a finite subprocess from the infinite cluster that is closest to the origin in a certain sense associated with the $f$-ordering. (Below we will provide a formal definition of a finite subprocess of $\mathbf X$). In the resulting configuration, 
the infinite cluster still remains infinite, but it contains a vertex of degree 1. Hence, the probability that the process takes values in the set of the resulting configurations is zero. Then it remains to show that also the probability that the process takes place in the set of original configurations is zero, which leads to the desired contradiction. This last step relies on the deletion-tolerance of $\X$ and uses a certain result from \cite{HS13} that we will recall below.

\medskip
Note that for the proof of Theorem~\ref{theorem-nopercolation}, we can assume that the intensity $\lambda$ of the underlying stationary point process is positive, since otherwise Theorem~\ref{theorem-nopercolation} is trivially true. We start with the following, previously proven lemma, which excludes existence of infinite clusters that have a degree-one point in the case of general random graphs based on stationary marked point processes.
\begin{lemma}{\cite[Lemma 5.4]{JT19}} \label{lemma-noheads}
Let $g(\X)$ be a random graph based on a marked point process $\X$ such that the degree of all $X_i\in X$, $\deg(X_i)$, is bounded by $2$, almost surely. Let $X$ be stationary and nonequidistant with a finite intensity, and consider the point process of degree-one points in infinite clusters
\[ \mathcal X_0=\sum_{i\in I}\delta_{X_i}\one\{\deg(X_i)=1,\, X_i\text{ is part of an infinite cluster in }g(\X)\}.\]
Then, $\P(\mathcal  X_0(\R^d)=0)=1$. 
\end{lemma}
The proof is based on a certain variant of the \emph{mass-transport principle} (see \cite[Section 4.2]{CDS20} for instance). Informally speaking, the proof goes as follows: If there was an infinite cluster having a point of degree one, then by stationarity, the point process of degree-1 points of infinite clusters $\mathcal X_0$ would have to have a positive density. This however leads to a contradiction because any infinite cluster can only contain at most one degree-1 point and must contain infinitely many degree-2 points, which implies that the aforementioned density must be equal to zero. We refer the reader to \cite[Section 5.2]{JT19} for further details.


%
In will be convenient to assume that $\mathbf X$ takes values in the set $\mathbf N$ of marked point configurations $\mathbf x $ in $\R^d  \times M$ such that $ x=\{ x_i \colon (x_i,p_i) \in  \mathbf x\}$ is an infinite, locally-finite, nonequidistant point configuration on $\R^d$. We will denote by $\rm N$ the set of such point configurations $x$ and equip $\mathbf N$ and ${\rm N}$ with the corresponding evaluation $\sigma$-fields.  

One essential property that we will use in the proof of Theorem~\ref{theorem-nopercolation} is the so-called \emph{edge-preserving property} of the underlying $f$-\kNN~graph. Informally speaking, this means that if we remove points from a configuration, then edges between remaining points are preserved. Now we provide the definition of this property corresponding to our configuration spaces.
\begin{defn}\label{defn-edgepreserving}
Let $g \colon \mathbf N \to \mathrm N \times (\mathrm N \times \mathrm N)$, $\mathbf x \mapsto g(\mathbf x)=(x,E_g(\mathbf x))$ be a function that maps a marked point configuration $\mathbf x$ to a graph with vertex set $x$. We say that $g$ is \emph{edge-preserving} if for all $\mathbf y, \mathbf x \in \mathbf N$ with $\mathbf y \subseteq \mathbf x$, for all $(v_1,q_1),(v_2,q_2) \in \mathbf y$ such that $(v_1,v_2) \in E_g(\mathbf x)$, one has $(v_1,v_2) \in E_g(\mathbf y)$.
\end{defn}
Let us, for the remainder of this section, fix an ordering function $f$. It is then easy to see that the $f$-\kNN~graph $g_{k,f} \colon \mathbf x \mapsto g_{k,f}(\mathbf x)$ is edge-preserving for all $k\in \N$. See Sections~\ref{sec-examplesf} and~\ref{sec-examplesnotf} for further examples of edge-preserving graphs. 

Now, for $\mathbf x \in \mathbf N$ and $v_o \in x$, we consider the sequence $(\NNN_n(v_o,\mathbf x))_{n \in \N_0}$ of the marked points of $\mathbf x$ ordered in increasing $f$-order of $\mathbf x$, measured from $v_o$. 
Then,  non-boldface $\NN_i(v_o,\mathbf x)$, as defined in Section~\ref{sec-model}, is the first component of $\NNN_i(v_o,\mathbf x)$, which we call the \emph{$i$-th nearest $f$-neighbor of $v_o$}. In particular, recall that $\NN_0(v_o,\mathbf x)=v_o$. 

Next, if $v_o$ has degree two in $g_{2,f}(\mathbf x)$, then $v_o$ is connected by an edge to both $\NN_1(v_o,\mathbf x)$ and $\NN_2(v_o,\mathbf x)$.
Moreover, both $\NN_1(v_o,\mathbf x)$ and $\NN_2(v_o,\mathbf x)$ also have $v_o$ as one of their first two nearest $f$-neighbors, i.e., 
\[ v_o \in\big \{ \NN_1(\NN_i(v_o,\mathbf x),\mathbf x), \NN_2(\NN_i(v_o,\mathbf x),\mathbf x) \big\}, \]
for all $i\in\{1,2\}$. 
These $k$-nearest $f$-neighbor relations hold almost surely and in particular for every nonequidistant configuration $\mathbf x$. Thanks to the choice of our configuration space $\mathbf N$, we can entirely exclude configurations that offend the degree bound or the $k$-nearest $f$-neighbor relations or are not nonequidistant. In particular, the $f$-2NN graph $g_{2,f}(\mathbf x)$ is well-defined for all $\mathbf x \in \mathbf N$. 

Given this, for $\mathbf x \in \mathbf N$, we define $\ell(\mathbf x) \in [0,\infty]$ as the number of infinite clusters in $g_{2,f}(\mathbf x)$, and in case $\ell(\mathbf x)>0$, we let $\ZZ(\mathbf x)=(z(\mathbf x),r(\mathbf x))$ denote the closest point of $x$ to the origin such that $z(\mathbf x)$ has degree two and is contained in an infinite cluster in $g_{2,f}(\mathbf x)$, and we write $\mathcal C_1(\mathbf x)$ for the cluster containing $z(\mathbf x)$. 

Now, in order to show that the probability that $g_{2,f}(\mathbf X)$ percolates is zero, it suffices to verify that
\[ \P(\ell(\mathbf X)\ge 1)=0. \numberthis\label{onotinaninfinitecluster} \]
Next, for $\mathbf x \in \mathbf N$ with $\ell(\mathbf x) \geq 1$, we define the third-nearest $f$-neighbor of $z(\mathbf x)$ within the infinite cluster $\mathcal C_1(\mathbf x)$ by
\[ \I(\mathbf x) = \inf \{ i \geq 3 \colon \NN_i(z(\mathbf x),\mathbf x) \in \mathcal C_1(\mathbf x) \}. \]
Indeed, among the points $\NN_j(z(\mathbf x),\mathbf x)$, $j < \I(\mathbf x)$, precisely two (namely, $\NN_1(z(\mathbf x),\mathbf x)$ and $\NN_2(z(\mathbf x),\mathbf x)$) are contained in $\mathcal C_1(\mathbf x)$. In case $\ell(\mathbf x)=0$, we put $\I(\mathbf x)=\infty$. Then, we have the following assertion.
\begin{prop}\label{lemma-inprop2deg}
Under the assumptions of Theorem~\ref{theorem-nopercolation}, for any natural number $i \geq 3$, we have 
\[ \P(\{ \ell(\mathbf X)\ge 1 \}\cap \{ \I(\mathbf X)=i \})=0. \numberthis\label{there'snobody} \]
\end{prop}
Before proving Proposition~\ref{lemma-inprop2deg}, let us show why it implies Theorem~\ref{theorem-nopercolation}.
\begin{proof}[Proof of Theorem~\ref{theorem-nopercolation}]
Noting that $\{ \ell(\mathbf X)\ge 1  \}\subset \{\I(\mathbf X)<\infty\}$ and using a union bound, Proposition~\ref{lemma-inprop2deg} implies $\P(\ell(\mathbf X)\ge 1)=0$, which is \eqref{onotinaninfinitecluster}, and thus finishes the proof of Theorem~\ref{theorem-nopercolation}.
\end{proof}
%
In order to verify Proposition~\ref{lemma-inprop2deg}, let us now present a preliminary result. First, we say that a point process $\mathbf Y$ is a \emph{finite subprocess} of $\mathbf X$ if $\mathbf Y$ is defined on the same probability space as $\X$, satisfies $\# \mathbf Y<\infty$ and $\mathbf Y \subset \mathbf X$ almost surely. Then we have the following claim: If $\mathbf X$ is deletion-tolerant, then for any subprocess $\mathbf Y$ of $\mathbf X$, the law of $\mathbf X \setminus \mathbf Y$ is absolutely continuous with respect to the one of $\mathbf X$.
Indeed, for $M=\R^{l}$, $l\geq 0$, the claim follows from \cite[Theorem 1.1]{HS13}. Now we verify it in the case of general $M$ analogously to the proof presented for the case $M=\R^{l}$ in \cite[Section 3]{HS13}. Let us assume that for some measurable subset $\mathcal A$ of $\mathbf N$ we have
$\P(\mathbf X \setminus \mathbf Y \in \mathcal A)>0$. Then, there must exist some $n\in\N$ such that $\P(\{ \mathbf X \setminus \mathbf Y \in \mathcal A \} \cap \{ \# \mathbf Y = n \})>0$. Set $\mathbf Y_n=\mathbf Y$ if $\# \mathbf Y = n$ and set $\mathbf Y_n = \emptyset$ otherwise. Then, it follows that $\P(\mathbf X \setminus \mathbf Y_n \in \mathcal A)>0$. Since $\mathbf X$ is deletion-tolerant, it follows that $\P(\mathbf X \in \mathcal A)>0$, as claimed.

\begin{comment}{
\begin{lemma}[{\cite[Theorem 1.1]{HS13}}]\label{theorem-HS13}
Let $X$ satisfy the assumptions of Theorem~\ref{theorem-nopercolation} apart from the one of deletion-tolerance. Then the following assertions are equivalent.
\begin{enumerate}
    \item\label{deletiontolerance} $X$ is deletion-tolerant.
    \item\label{finiteprocess} For any finite subprocess $Y$ of $X$, the law of $ X \setminus Y$ is absolutely continuous with respect to the one of $X$.
    \item\label{voidprobnonzero} For any bounded $S \in \mathcal B(\R^d)$, we have that almost surely
    \[ \P\big(X \cap S=\emptyset|X \cap S^c\big) >0.\numberthis\label{void} \]
\end{enumerate}
\end{lemma}
}
\end{comment}

Next, we observe the following. For $\mathbf x\in \mathbf N$ with $ \ell(\mathbf x)\ge 1$, by definition, we have that $z(\mathbf x)$ is connected by an edge both to $\NN_1(z(\mathbf x), \mathbf x)$ and $\NN_2(z(\mathbf x), \mathbf x)$ in $g_{2,f}(\mathbf x)$. Further, since the degrees are bounded by two, for such $\mathbf x$, $\NN_1(z(\mathbf x), \mathbf x)$ and $\NN_2(z(\mathbf x), \mathbf x)$ have no further joint neighbor in $g_{2,f}(\mathbf x)$ since otherwise $\mathcal C_1(\mathbf x)$ has a loop and cannot be infinite.
Hence, for any $i \geq 3$, there exists $l \in \{1,2\}$ such that $\NN_i(z(\mathbf x), \mathbf x)$ and $\NN_l(z(\mathbf x), \mathbf x)$ are not connected by an edge in $g_{2,f}(\mathbf x)$. Let us denote the corresponding $\NN_l(z(\mathbf x), \mathbf x)$ by $m_i(\mathbf x)$, and define $m_i(\mathbf x)=\NN_1(z(\mathbf x), \mathbf x)$ if neither $\NN_1(z(\mathbf x), \mathbf x)$ nor $\NN_2(z(\mathbf x), \mathbf x)$ is connected to $\NN_i(z(\mathbf x), \mathbf x)$ by an edge. The element of $\{ \NN_1(z(\mathbf x), \mathbf x),\NN_2( z(\mathbf x), \mathbf x) \}$ not being equal to $m_i(\mathbf x)$ is denoted by $n_i(\mathbf x)$.
We will write $q_i(\mathbf x)$ for the mark of $m_i(\mathbf x)$. 

Let us define the set of configurations in which the infinite cluster closest to the origin is one-armed
\[ B = \{\mathbf x \in \mathbf N \colon  \ell(\mathbf x)\ge 1 \text{ and }\mathcal C_1(\mathbf x) \text{ contains a point of degree one}\}, \]
and note that $B \subset \{\mathbf x \in \mathbf N \colon \ell(\mathbf x) \geq 1\}$. With this, for $\mathbf x \in \mathbf N$, let us define the finite point configuration
\[ \mathbf y(\mathbf x) = 
\begin{cases}
 \{\NNN_3(z(\mathbf x),\mathbf x), \dots, \NNN_{\I(\mathbf x)-1}(z(\mathbf x),\mathbf x),(m_{\I(\mathbf x)}(\mathbf x),q_{\I(\mathbf x)}(\mathbf x))\}& \text{ if } \ell(\mathbf x) \geq 1\text{ and } \mathbf x \notin B, \\
\emptyset, & \text{ if } \ell(\mathbf x) =0\text{ or } \mathbf x \in B.
\end{cases} \]
Let us fix $i \geq 3$. Let $\mathbf x \in \mathbf N$ be such that $\ell(\mathbf x) \geq 1$ and $\I(\mathbf x)=i$. Then, we define a thinned configuration
\begin{align*} 
\mathbf x^i =\mathbf x\setminus \mathbf y(\mathbf x)
\end{align*}
and put $x^i=\{ v \colon (v,q) \in \mathbf x^i  \}$ for the correspoding non-marked configuration. 

We claim that $\{\ell(\mathbf x)\ge 1\}\cap \{\I(\mathbf x)=i\}$ is a subset of $\{ \ell(\mathbf x^i) \geq 1 \}$. For this, first note that the removal of finitely-many points of the point process can change the edge structure of the remaining points. Nevertheless, it cannot remove edges that were already present before the removal of points since $g_{2,f}$ is edge-preserving. 
Hence, all edges between two points of $x^i$ in $g_{2,f}(\mathbf x)$ also exist in $g_{2,f}(\mathbf x^i)$. In particular, if $g_{2,f}(\mathbf x)$ contains an infinite cluster, then $g_{2,f}(\mathbf x^i)$ contains all but at most finitely many edges of this cluster. This implies that $\ell(\mathbf x^i)\ge 1$, hence the claim.

\begin{figure}[!htpb]
\centering
\begin{tikzpicture}[scale=3,
        shorten >=1pt, auto, thick,
        node distance=1.8cm,
    main node/.style={circle,draw,font=\sffamily\large\bfseries}
                            ]

			\node[main node] (1) {$z$};
			\node[main node,color=red] (2) [right of=1] {$\NN_1=m_3$};
			\node[main node] (3) [left of=1] {$\NN_2=n_i$};
			\node[main node] (5) [right of=2] {};
			\node[main node,color=white] (8) [right of=5] {};
			\node[main node] (6) [left of=3] {};
			\node[main node] (4) [left of=6] {$\NN_i$};
			\node[main node,color=white] (7) [left of=4] {};
			\node[main node,color=red] (9) [below of=1] {$\NN_3$};
			\node[main node,color=white] (10) [below of=3] {\color{red}$\dots$\color{white}};
			\node[main node,color=red] (11) [below of=6] {$\NN_{i-1}$};

			\path[every node/.style={font=\sffamily\small}]
			  (1) edge [color=red] (2)
			      edge (3)   
			  (6) edge [dashdotted] (4)
			  (4) edge (7)
			  (2) edge [color=red] (5)
			  (5) edge [dotted] (8)
			  (3) edge [dashdotted] (6)
			  (3) edge [dashed,bend left] (4)
			  (4) edge [dotted, bend left] (8);
\end{tikzpicture}
\caption{An illustration of the removal of the finite subconfiguration $\mathbf y(\mathbf x)=\{ (m_i,q_i),\NNN_3(z(\mathbf x)),\ldots,\NNN_{n-1}(z(\mathbf x))\}$ from the configuration $\mathbf x$ satisfying $\ell(\mathbf x)\geq 1$, $\I(\mathbf x)=i \geq 3$ and $\mathbf x \notin B$. The point $\NN_i=\NN_i(z(\mathbf x),\mathbf x)$ is contained in the infinite cluster $\mathcal C_1=\mathcal C_1(\mathbf x)$ including $z=z(\mathbf x)$, and it is not a neighbor of $m_i=m_i(\mathbf x)$, which equals $\NN_1=\NN_1(z(\mathbf x),\mathbf x)$ here, whereas $\NN_2=\NN_2(z(\mathbf x),\mathbf x)=n_i=n_i(\mathbf x)$. (In the figure we use the short-hand notations $\NN_j=\NN_j(z(\mathbf x),\mathbf x)$ for all values of $j$.) Hence, if $\NN_i$ has degree two in $\mathcal C_1$,  then there are various possibilities respecting the degree bound of two to connect $\NN_i$ to $\mathcal C_1$ so that it is not connected to $m_i$ by an edge. $\NN_i$ can either be a direct neighbor of $\NN_2$ (see dashed line) or a further point of the path from $z$ to infinity starting with the edge from $z$ to $\NN_2$ (dash-dotted lines) or a non-direct neighbor of $\NN_1$ on the path from $z$ to infinity starting with the edge from $z$ to $\NN_1$ (dotted lines). \\
Now, removing the finite subconfiguration $\mathbf y(\mathbf x)$ from the realization (the non-marked points corresponding to $\mathbf y$ are colored red in the figure), in the resulting $f$-2NN graph of $\mathbf x^i=\mathbf x \setminus \mathbf y(\mathbf x)$, both edges adjacent to $\NN_i$ are preserved. Also all edges from $z$ to infinity starting with the edge from $z$ to $\NN_2$ are preserved, hence $z$ is still contained in an infinite cluster, but the edge from $z$ to $\NN_1$ is removed. In the resulting new configuration, the second-nearest $f$-neighbor towards $z$ is $\NN_i$, and hence this is the only point of the configuration that could be connected to $z$ by an edge. But $\NN_i$ still cannot have degree 3 or more, hence it cannot be connected to $z$. Thus, in the new configuration, $z$ is in an infinite cluster and has degree 1.}
\label{fig-NN3}
\end{figure}
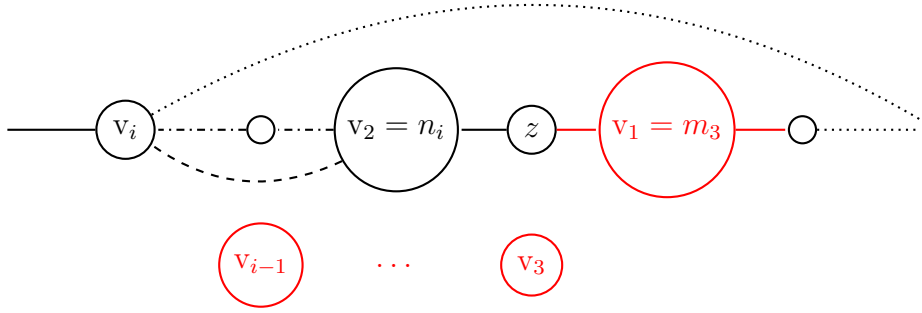

Then, the next claim is that for $\mathbf x \in \mathbf N$ with $\ell(\mathbf x)\geq 1$ and $\I(\mathbf x)=i$, we have that $\mathbf x^i$ is contained in $B$. The proof of this claim is clear if $\mathbf x \in B$ since then $\mathbf x^i=\mathbf x$. Thus, we now consider the case $\mathbf x\notin B$. This case is illustrated in Figure~\ref{fig-NN3}.
Recall that for $\mathbf x \in \mathbf N$, $z(\mathbf x)$ cannot have degree higher than two in $g_{2,f}(\mathbf x^i)$, whereas it has degree at least one and its cluster $\mathcal C_1(\mathbf x^i)$ is infinite in $g_{2,f}(\mathbf x^i)$. Note also that the edge between $z(\mathbf x)$ and $n_i(\mathbf x)$ still exists in $g_{2,f}(\mathbf x^i)$. Further, if $z(\mathbf x)$ has degree two in $g_{2,f}(\mathbf x^i)$, then it is connected to the second-nearest $f$-neighbor towards $z(\mathbf x)$ in $\mathbf x^i$, which is $\NN_2(z(\mathbf x),\mathbf x^i)=\NN_i(z(\mathbf x),\mathbf x)$, whereas $\NN_1(z(\mathbf x),\mathbf x^i)=n_i(\mathbf x)$. Now, since $\mathbf x \notin B$, $\ell(\mathbf x) \ge 1$ and $\NN_i(z(\mathbf x),\mathbf x) \in \mathcal C_1(\mathbf x)$, it follows that $\NN_i(z(\mathbf x),\mathbf x)$ has degree equal to two in $g_{2,f}(\mathbf x)$. Further, it is neither connected to $m_i(\mathbf x)$ by an edge nor to $z(\mathbf x)$ in this graph. Hence, both edges adjacent to $\NN_i(z(\mathbf x),\mathbf x)$ also exist in $g_{2,f}(\mathbf x^i)$. But since $\NN_i(z(\mathbf x),\mathbf x)$ has degree at most two in $g_{2,f}(\mathbf x^i)$, it follows that $z(\mathbf x)$ and $\NN_i(z(\mathbf x),\mathbf x)$ are not connected by an edge in this graph. Hence, $\mathbf x^i \in B$, which implies the claim.

Note that by Lemma~\ref{lemma-noheads}, the set $B$ is a $\P$-null set, i.e.,
\[ \P(\mathbf X \in \{ \mathbf x^i \colon \mathbf x\in \mathbf N, \ell(\mathbf x )\ge1 \text{ and }\I(\mathbf x)=i \})=0.\numberthis\label{isunlikely} \] 
This implies \eqref{there'snobody}. Next, we have the following lemma. %
\begin{lemma}\label{lemma-omegaomegai} 
Under the assumptions of Theorem~\ref{theorem-nopercolation}, for any $i \geq 3$, $\P(\{\ell(\mathbf X)\ge1\}\cap \{ \I(\mathbf X)=i \})>0$ implies $\P(\mathbf X \in \{ \mathbf x^i \colon \mathbf x \in \mathbf N, \ell(\mathbf x)\text{ and } \I(\mathbf x)=i \}) >0$.
\end{lemma}
Before proving the lemma, let us now explain why it implies Proposition~\ref{lemma-inprop2deg}.
\begin{proof}[Proof of Proposition~\ref{lemma-inprop2deg}]
By Lemma~\ref{lemma-omegaomegai}, where we show that if the collection of thinned configurations is contained in a $\P$-null set, also the non-thinned configurations form a $\P$-null set, we see that \eqref{isunlikely} implies \eqref{there'snobody}, which concludes the proof of Proposition~\ref{lemma-inprop2deg}.
\end{proof}
Finally, it remains to verify Lemma~\ref{lemma-omegaomegai}. Its proof strongly relies on deletion-tolerance. 
\begin{proof}[Proof of Lemma~\ref{lemma-omegaomegai}]

Let us fix $i \geq 3$ and assume that $\P(\{\ell(\mathbf X)\ge1\}\cap \{\I(\mathbf X)=i \})>0$.
For the finite subprocess $\mathbf Y=\mathbf y(\mathbf X)$ of $\mathbf X$ this means that
\[ \begin{aligned}
0 &< \P(\{\ell(\mathbf X)\ge1\}\cap \{\I(\mathbf X)=i \})=\P\big( \mathbf X \in \{ \mathbf x \in \mathbf N \colon \ell(\mathbf x)\ge 1,\I(\mathbf x)=i\} \big) 
\\ & \leq \P\big( \mathbf X \setminus \mathbf Y \in \big\{ \mathbf x^i \colon \ell(\mathbf x)\ge 1,\I(\mathbf x)=i\}\big).
\end{aligned}\]
Next, since $\mathbf X$ is assumed to be deletion-tolerant and $\mathbf Y$ is a finite subprocess of $\mathbf X$, by the equivalent characterization of deletion-tolerance as presented below Proposition~\ref{lemma-inprop2deg}, it follows that the law of $\mathbf X \setminus \mathbf Y$ is absolutely continuous with respect to the one of $\mathbf X$. This yields
\[ 0 < \P\big( \mathbf X  \in \big\{ \mathbf x^i \colon \ell(\mathbf x)\ge 1,\I(\mathbf x)=i \big\}\big). \]
This implies the lemma.
\end{proof}

%% file: 2degree-Examples.tex
\subsection{Examples of deletion-tolerance}\label{sec-thinningexamples}
In this section, we verify Propositions~\ref{prop-thinning} and \ref{prop-rigidity}. We first carry our the proof of Proposition~\ref{prop-thinning}.
\begin{comment}{In general, Proposition~\ref{prop-thinning} claims that under $\mathbb P'$, the distribution of $\mathbf X^{K,p}$ is absolutely continuous with respect to the one of $\mathbf X$. Let $F$ be an element of the evaluation $\sigma$-algebra of $\mathbf N$ such that $\P'(\mathbf X^{K,p} \in F)>0$. We have to show that then also $\P(\mathbf X \in F)>0$. Under the assumption that $\P'(\mathbf X^{K,p} \in F)>0$, by continuity of measures, we can find $K,l \in \N$ such that
\[ \eps:=\P'(\mathbf X^{K,p} \in F, \# (\mathbf X^{K,p} \cap (B_K(o) \times M)) = l) >0. \numberthis\label{eps1prob}\]
In other words, we have $0<\eps = \P'(\mathbf X^{K,p} \in G)$ where $G = \{ \bs \omega \in F \colon \# (\bs \omega \cap ( B_K(o) \times M))=l \}$.
Thus,
\[
\begin{aligned}
 \P(\mathbf X \in F)  &\geq \P'(\mathbf X \in G, \mathbf X^{K,p} = \mathbf X) 
 \geq \P'(\mathbf X^{K,p} \in G)\P'(\mathbf X^{K,p} = \mathbf X | \mathbf X^{K,p} \in G)\\
 & = \eps \P'(\mathbf X^{K,p}= \mathbf X | \mathbf X^{K,p} \in G),
\end{aligned}
\] 
 and further,
\[ 
\begin{aligned}
 \P'(\mathbf X^{K,p} = \mathbf X | \mathbf X^{K,p} \in G) & 
 \geq \P'(\mathbf X^{K,p} = \mathbf X, \mathbf X^{K,p} \in G)
 = \P'(\mathbf X' = \emptyset, \mathbf X^{K,p} \in G).
\end{aligned} \numberthis\label{probofcond} \]
%
Now we perform a case distinction depending of the type of the point process. We start with the case of CPPs.}\end{comment}

\begin{proof}[Proof of Proposition~\ref{prop-thinning}]
First, let $\X$ be an i.i.d.~marked CPP, where $X$ is stationary and non-equidistant. Then, due to the i.i.d. ~markings, $\X$ can be seen as a PPP with random directing measure $\Gamma(\d v)\otimes \mathcal Q(\d q)$, where $\mathcal Q$ is the distribution of the marks and $\Gamma$ is the (random) intensity of the unmarked CPP, see e.g.~\cite[Section 13]{LP17}. Consider an event $\mathcal A$ such that $\P(\X\in \mathcal A)=0$, then, by the definition of deletion tolerance, it suffices to show that for any bounded measurable set $B\subset \R^d$, we have 
$$\E[\sum_{X_i\in X\cap B}\one\{\X\setminus \X_i\in \mathcal A\}]=0.$$
But, using the Mecke formula for PPPs in general measurable spaces, see~\cite[Theorem 13.8]{LP17}, we have that 
\[
\begin{aligned} 
\E[\sum_{X_i\in X\cap B}\one\{\X\setminus \X_i\in \mathcal A\}]&=\E[\Gamma(B)\P( \X\in \mathcal A|\Gamma)]=0,
\end{aligned}
\]
as desired. 
%

\medskip
Let now $\X$ be a stationary Gibbs point processes. 
As in the proof of Proposition~\ref{prop-rigidity} below, we want to employ an equivalent characterization of deletion tolerance via almost-sure positivity of conditional void probabilities, see~\cite[Theorem 1.1]{HS13} for the case of unmarked point processes. In order to lift this to the level of marked point processes
we extend the arguments as exhibited in~\cite{HS13}. However, first note that by our assumptions, for all bounded measurable $\L\subset\R^d$, and almost all realizations $\X_{\L^{\rm c}}$ of $\X$ in $\L^{\rm c}$, we have
$$\P(\X_\L=\emptyset| \X_{\L^{\rm c}})\ge \e^{-\la|\L|}Z^{-1}_\L(\X_{\L^{\rm c}})>0,$$
where $\la>0$ denotes the intensity of the underlying unmarked Poisson point process. 

Using this, on a general level, if for some bounded $\L\subset \R^d$ and measurable $\mathcal A\subset \R^d\times M$, we have $\P(\X_{\L^{\rm c}}\in \mathcal A)>0$, then $\P(\X\in \mathcal A)\ge \P(\X_{\L^{\rm c}}\in \mathcal A,\X_\L=\emptyset)=\E[\one\{\X_{\L^{\rm c}}\in \mathcal A\}\P(\X_\L=\emptyset| \X_{\L^{\rm c}})]>0$ and hence, the law of $\X_{\L^{\rm c}}$ is absolutely continuous with respect to the law of $\X$. In order to deduce deletion tolerance from this absolute continuity, we consider open balls $B(\bs x,r)$ in $\R^d\times M$, where $\bs x\in \Q^d\times D$, with $D$ a countable-dense subset of the Polish space $M$, and $r\ge 0$ is rational. We let $C$ denote the set of unions of finitely many such balls. Then, by the local finiteness of $\X$, for all $\X$-points $\bs Z$, there exists a $C$-valued random variable $\bs \Gamma$ such that $\P(\X_{\bs\Gamma}=\delta_{\bs Z})=1$. But then, if for some measurable $\mathcal A\subset \R^d\times M$ we have $\P(\X-\delta_{\bs Z}\in \mathcal A)>0$, it follows that 
\[
\begin{aligned} 
0&<\P(\X-\delta_{\bs Z}\in \mathcal A, \X_{\bs\Gamma}=\delta_{\bs Z})=\P(\X_{\bs\Gamma^{\rm  c}}\in \mathcal A)=\sum_{\Gamma\in C}\P(\X_{\bs\Gamma^{\rm  c}}\in \mathcal A, \bs\Gamma=\Gamma),
\end{aligned}
\]
and hence $\P(\X_{\Gamma^{\rm  c}}\in \mathcal A)>0$ for some $\Gamma$. But this implies $\P(\X\in \mathcal A)>0$ by the absolute continuity asserted above. 
\end{proof}

Finally, we prove Proposition~\ref{prop-rigidity}. The proof is based on the following lemma.
\begin{lemma}\label{lemma-markednonmarked}
Let the marked point process $\mathbf X$ be deletion-tolerant. Then, the underlying point process $X$ is also deletion-tolerant.
\end{lemma}
In order to keep the proof short, we again assume that $\mathbf X$ takes values in the configuration space $\mathbf N$ and $X$ in the configuration space $\mathrm N$, which were defined in Section~\ref{sec-proof}. 
\begin{proof}
Let $A$ be any measurable subset of $\mathbf N$. Further, let $Y$ be any $X$-point. Let $P$ be the mark of $Y$ (defined realizationwise). Then $(Y,P)$ is an $\mathbf X$-point and we have
\[ \P(X \in A)=\P\big(X \in A, (P_i)_{i\in\N} \in M^{\N} \big) = \P(\mathbf X \in A \times M^{\N}). \]
Since $A\times M^{\N}$ is a measurable subset of $\mathbf N$ and $\mathbf X$ is deletion-tolerant, we have
\[ 0<\P(\mathbf X \setminus \{ (Y,P) \} \in A \times M^{\N})
= \P(X \setminus \{ Y \} \in A). \]
Hence, the distribution of $X \setminus \{Y\}$ is absolutely continuous with respect to the one of $X$, as wanted.
\end{proof}
\begin{proof}[Proof of Proposition~\ref{prop-rigidity}]
Thanks to Lemma~\ref{lemma-markednonmarked}, it suffices to show that $X$ is not deletion-tolerant. Further, as already mentioned in the proof of Proposition~\ref{prop-thinning}, the following assertion is known from \cite[Theorem 1.1]{HS13}: If $X$ is deletion-tolerant, then for any bounded $S \in \mathcal B(\R^d)$, we have that almost surely
\[ \P\big(X \cap S=\emptyset|X \cap S^c\big) >0.\numberthis\label{void} \]

Assume now that the conditions of the proposition are satisfied. Let $\mathcal K \in \mathcal B(\R^d)$ be compact. Since $X$ is stationary with positive intensity, there exists some $k \in \N$ such that
\[ \P\big(\#( X \cap \mathcal K) =k)>0. \]
Let us now consider the measurable set $F=h_{\mathcal K}^{-1}(k)$. We have \[ \P(\# (X \cap \mathcal K^c)\in F) >0. \]
Hence, the conditional probability
\[ \P( X \cap \mathcal K=\emptyset | X \cap \mathcal K^c \in F)  \numberthis\label{zeroprob} \]
is well-defined. However, since $\{ X \cap \mathcal K^c \in F \} \subseteq \{ \#(X \cap \mathcal K)=k \}$, \eqref{zeroprob} equals zero. This shows that for our choice of $X$ and for $S=\mathcal K$, the assertion~\eqref{void} does not hold, and hence $X$ is not deletion-tolerant.
\end{proof}

\begin{comment}
{In the Poisson case it's easy to show that the r.h.s.~is nonzero because we have
\[
\begin{aligned}
\P'&(\mathbf X' = \emptyset, \mathbf X^{K,p} \in F, \# (\mathbf X^{K,p} \cap (B_K(o) \times M)) = l) \\ & = \P'(\mathbf X' = \emptyset) \P'(\mathbf X^{K,p} \in F, \# (\mathbf X^{K,p} \cap (B_K(o) \times M)) = l) \\ 
& = \eps_1 \e^{-\lambda(1-p)|B_K(o)|} >0.
\end{aligned} \]
Similarly, in the Cox case we obtain
\[
\begin{aligned}
\P'&(\mathbf X' = \emptyset, \mathbf X^{K,p} \in F, \# (\mathbf X^{K,p} \cap (B_K(o) \times M)) = l) \\ & = \E'\big( \P'(\mathbf X^{K,p} \in F, \# (\mathbf X^{K,p} \cap (B_K(o) \times M)) = l |\Lambda ) \e^{-\lambda(1-p)\Lambda(B_K(o))}  \big).
\end{aligned} \]
\color{black}}\end{comment}
\color{black}
\color{black}

\subsection{SINR graphs as subgraphs of $f$-\kNN~graphs}\label{sec-examplesf}
Let us briefly summarize the relation between Theorem~\ref{theorem-nopercolation} and \cite[Theorem 2.2]{JT19}. Indeed, the two proofs are similar, in particular, two steps of the proof, Lemma~\ref{lemma-noheads} and the part of Proposition~\ref{prop-thinning} regarding the Cox case already appeared in \cite{JT19}. However, in \cite{JT19} we focused on the particular case of SINR graphs based on stationary and nonequidistant CPPs, having concrete applications in telecommunications in mind, and we did not aim at checking whether our proof works also for a wider class of point processes or graphs. Thus, the main novelty in this paper is not that we exploit new proof techniques (although the proofs of Proposition~\ref{prop-rigidity} and the part of Proposition~\ref{prop-thinning} regarding Gibbs point processes have no analogues in \cite{JT19}). Instead, we highlight that apart from the general combinatorial condition of working with undirected and stationary random graphs with degrees bounded by two, two properties are crucial for a straightforward generalization of the proof in \cite{JT19}: (1) the deletion tolerance of the underlying point process and (2) the edge-preserving property of the graph. The latter observation allows for the extensions of Theorem~\ref{theorem-nopercolation} presented in Section~\ref{sec-examplesnotf}.

This puts the result into a general framework and allows for generalizations of the result both with respect to the type of graph and with respect to the kind of point process. Here, let us note that the SINR graph is not a special case of an $f$-\kNN~graph, but a proper subgraph of an $f$-2NN~graph under particular choices of the parameters, which is itself edge-preserving.
Non-percolation in $f$-2NN~graphs was not even known before the present paper in the simplest case represented by the \B-$2$NN~graph. 
\begin{comment}{
First, the assertion of the present paper is stronger since SINR graphs with degrees bounded by 2 are proper subgraphs of $f$-2NN~graphs under particular choices of the parameters (see Section~\ref{sec-examplesf}). While many arguments of the two proofs are similar,
the general condition of deletion-tolerance (instead of considering Cox point processes only) extends the applicability of our results to a larger class of point processes. Further, while the edge-preserving property (which holds also for SINR graphs, see Section~\ref{sec-examplesf}) was also used in \cite{JT19}, it was not exploited there that it holds for many other bounded-degree graphs apart from the SINR graph itself. The latter observation allows for the extensions of Theorem~\ref{theorem-nopercolation} presented in Section~\ref{sec-examplesnotf}}\end{comment}

In order to make the relation between $f$-\kNN~graphs and SINR graphs explicit we recall the definition and interpretation of the latter graphs. 
Let $M=N=[0,\infty)$, $\Vert \cdot \Vert=\Vert \cdot \Vert_2$, $P_o$ be a nonnegative random variable,  and $\mathbf X=\{(X_i,P_i)\}_{i\in I}$ an i.i.d.~marked point process in $\R^d \times [0,\infty)$ such that all $P_i$ are distributed as $P_o$.  Let $\ell \colon (0,\infty) \to [0,\infty)$, the so-called \emph{path-loss function}, be a monotone decreasing function. Typical examples of path-loss functions correspond to \emph{Hertzian propagation} (see~\cite{DBT05,DF06}), e.g., for $\alpha>d$, the unbounded function $\ell(r)=r^{-\alpha}$, its truncated variant $\ell(r)=\min \{ 1,r^{-\alpha} \}$, and its ``shifted'' variant $\ell(r)=(1+r)^{-\alpha}$. Now, define $f(x,p)=1/(p\ell(\Vert x\Vert))$. In a telecommunication context, for $(X_i,P_i) \in \mathbf X$ and $x \in \R^d$, $P_i$ expresses the signal power transmitted by a device at spatial position $X_i$, and $\ell$ describes propagation of signal strength over distance.
Note that $\ell$ need not be strictly decreasing, which gives relevance to  Part~\eqref{tiebreaking} of Definition~\ref{defn-closer} in order to make the $f$-ordering well-defined. We observe that in case $P_o$ is almost surely equal to a fixed positive constant, then the arising $f$-\kNN~graph is the \B-\kNN~graph.  

In this setting, the \emph{SINR graph} \cite{DBT05} is usually introduced in the following way. Let $N_o$ be another nonnegative random variable independent of $\mathbf X$.  Choose two further parameters $\gamma, \tau >0$, the so-called \emph{interference-cancellation factor} and the \emph{SINR threshold}, respectively, and for $i,j \in I$, $i \neq j$, connect $X_i$ and $X_j$ by an edge whenever the \emph{SINR constraint} is satisfied in both directions, i.e., 
\[ P_i \ell(|X_i-X_j|) >  \tau\Big (N_o +\gamma \sum_{k \in I \setminus \{ i,j\}} P_k \ell(|X_k-X_j|)\Big), \numberthis\label{SINR>tau} \]
and the same holds with the roles of $i$ and $j$ interchanged.
Then,  it is known from \cite[Theorem 1]{DBT05} that if $\mathbf X$ is a simple point process (even if not stationarity or not nonequidistant), all degrees in the SINR graph are less than $k=1+1/(\tau\gamma)$. Using the elementary arguments of \cite[Lemma 4.1.13]{T19}, one can easily verify that if $X=\{X_i\}_{i\in I}$ is also nonequidistant, then the SINR graph is a subgraph of the $f$-\kNN~graph of the present example. Further, if $N_o>0$ is deterministic, then the SINR graph has bounded edge lengths and hence is a subgraph of the Gilbert graph introduced in \cite{G61}. The same assertion holds also for $N_o=0$ in case $\ell$ has bounded support. 
For positive assertions about percolation in SINR graphs based on various kinds of point processes, we refer the reader e.g.\ to \cite{DBT05,DF06,BY13,T19,T18,L19,JT19}. 
 
 Hence, for point processes satisfying the conditions of  Theorem~\ref{theorem-nopercolation}, there is no percolation in the SINR graph if its degrees are bounded by 2, which is always the case if $\gamma \geq 1/(2\tau)$. Thanks to Proposition~\ref{prop-thinning}, in particular, Gibbs point processes are covered by this result. To the best of our knowledge, there have been no results about SINR percolation for Gibbs point processes before (apart from the degree bounds themselves). Regarding non-percolation in case degrees are bounded by two, the case of CPPs was handled in \cite[Theorem 2.2]{JT19}. Here, based on the observation \cite[Section 5.2, Proof of Proposition 5.8]{JT19} that SINR graphs are edge-preserving on their own right in the sense of Definition~\ref{defn-edgepreserving}, we carried out a certain variant of the proof of Theorem~\ref{theorem-nopercolation} for the SINR graph directly, with no direct reference to $f$-\kNN~(or even \B-\kNN) graphs.

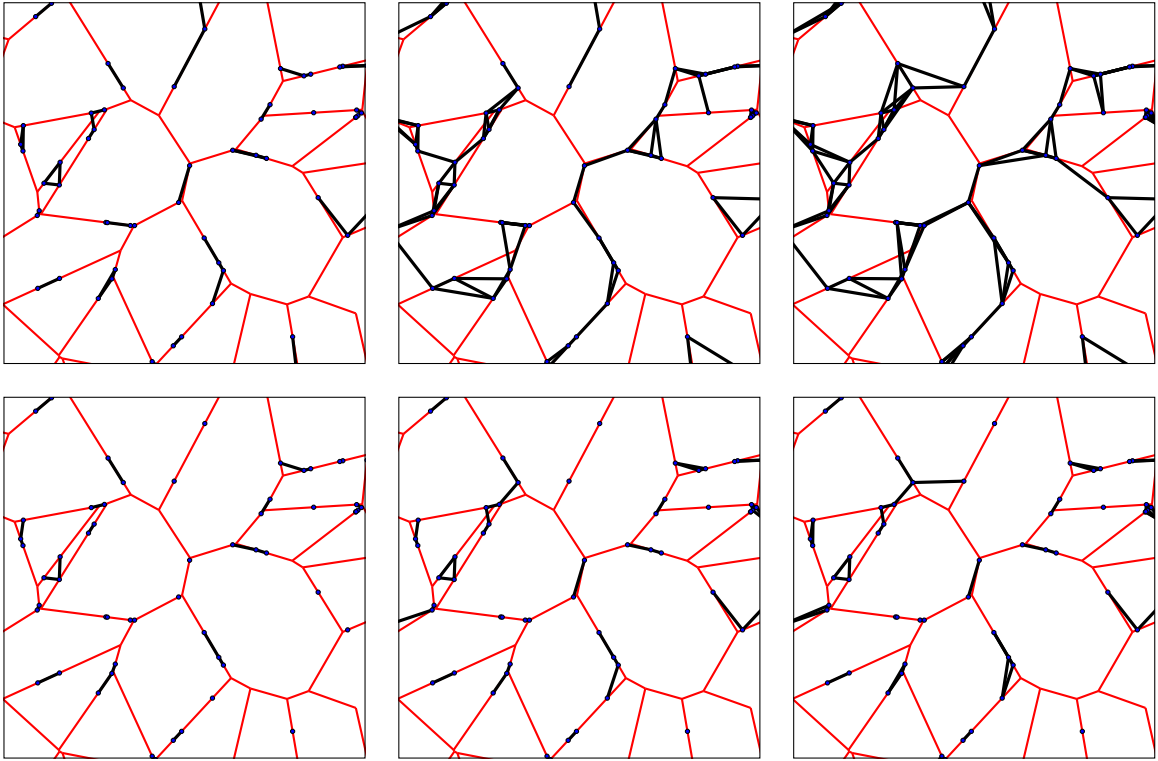
\begin{figure}[!htpb]
\begin{subfigure}{0.32\textwidth}
\input{VoronoiB2.tex}
\end{subfigure}
\begin{subfigure}{0.32\textwidth}
\input{VoronoiB4.tex}
\end{subfigure}
\begin{subfigure}{0.32\textwidth}
\input{VoronoiB5.tex}
\end{subfigure}
\\\vspace{0.4cm}
\begin{subfigure}{0.32\textwidth}
\input{VoronoiSINR2.tex}
\end{subfigure}
\begin{subfigure}{0.32\textwidth}
\input{VoronoiSINR4.tex}
\end{subfigure}
\begin{subfigure}{0.32\textwidth}
\input{VoronoiSINR5.tex}
\end{subfigure}
\caption{\B-\kNN~graphs (in the first line) and SINR graphs with degree bound $k$ (in corresponding order in the second line) for $k=2,4,5$, for $X$ being a stationary CPP. The random intensity measure $\Lambda$ is given as the edge-length measure (i.e., the one-dimensional Hausdorff measure) of a two-dimensional Poisson--Voronoi tessellation. The simulation leads to the conjecture that the smallest $k$ such that the \B-\kNN~graph percolates is $k=5$, which would mean that it equals the one of the two-dimensional PPP (which is 5 according to the high-confidence results of \cite{BB08}). Further, it is known from \cite{T18} that in this case, for large enough $\lambda$ and accordingly chosen small $\gamma>0$, there is also percolation in the SINR graph. However, it does not seem to be the case that this already happens when the degree bound equals 5, as the simulation suggests. Here, $\gamma$ is just slightly bigger than $1/(5\tau)$, i.e., a small further increase of $\gamma$ would increase the degree bound to 6, but the SINR graph is still much less connected than the corresponding \B-5${\rm NN}$ graph.}
\label{Fig-B}
\end{figure}

The aforementioned positive results on SINR percolation guarantee that for various kinds of point processes, the SINR graph percolates with positive probability for some positive $\gamma$ given that $\lambda$ is sufficiently large, while all the other parameters (depending on the type of point process) are fixed. In other words, we know that $k^*<\infty$ holds for the infimum $k^*$ of all degree bounds $k$ such that there exists an SINR graph with largest degree equal to $k$ that percolates. There are multiple interesting open questions related to this. First, what is the smallest value of such $k^*$, and how does it depend on the type of the point process?  The main results of the present paper imply that for stationary and nonequidistant point processes that are deletion-tolerant, we have $k^* \geq 3$. Further, according to the high-confidence results of \cite{BB08}, $k^* \geq 5$ for the two-dimensional PPP. Second, is the smallest such degree bound the same for SINR graphs as for the underlying \B-\kNN~graph? While the relationship between Gilbert graphs and SINR graphs is clear (namely, Gilbert graphs are increasing limits of SINR graphs as $\gamma \downarrow 0$), we are not aware of results stating that the \B-\kNN~graph is an increasing limit of certain SINR graphs with degree bound $k$, and such a result may not be true in general. Namely, it may be the case that an SINR constraint of the form~\eqref{SINR>tau} with degree bound $k \in \N$ poses stronger restrictions on the edges of the graph than a \B-\kNN~constraint for the same $k$. We defer the investigation of such questions to future work, noting that numerical evidence indicates that the two critical degree bounds are not the same in general, see e.g.~Figure~\ref{Fig-B}. 

\subsection{Extensions and limitations of the proof of Theorem~\ref{theorem-nopercolation}}\label{sec-examplesnotf}
We now present examples of graphs with degrees bounded by two that are not contained in an $f$-$2$NN graph but have rather similar properties to it, to the extent that the proof techniques of Theorem~\ref{theorem-nopercolation} are applicable to it.
\begin{example}[Locally furthest neighbors]
The edge-preserving property of $f$-\kNN~graphs (see Definition~\ref{defn-edgepreserving}) also holds if we replace the ``$k$-nearest neighbors with respect to the $f$-ordering'' in their definition by ``$k$-furthest neighbors in a bounded (possibly random) set shifted to the point, w.r.t.~$f$-ordering''. For the sake of simplicity of notation, let us explain how this works in the case of \B-\kNN~graph. The case of general $f$-\kNN~graphs can be handled similarly, taking into account also the marks and using the $f$-ordering instead of the ordering of Euclidean norms. We assume throughout this discussion that the point process $X$ is stationary, nonequidistant, and deletion-tolerant.

Let us fix a deterministic measurable set $A \subseteq \R^d$ of finite Lebesgue measure and define a random graph with vertex set $\mathbf X$ via connecting two different points $X_i,X_j$ of the point process $X$ by an edge whenever $X_j$ is one of the $k \in \N$ furthest neighbors in $(A+X_i) \setminus \{ X_i \}$ and the same holds with the roles of $i$ and $j$ interchanged. It is easy to see that this graph is well-defined and edge-preserving. Clearly, for $k=2$ it has degrees bounded by two. Hence, non-percolation of the graph can be verified along the lines of the proof of Theorem~\ref{theorem-nopercolation} also in the case of this graph. If $A$ is bounded, then the graph has bounded edge lengths (unlike the \B-\kNN~graph). 

This approach can be extended to the case when the deterministic set $A+X_i$ is replaced by a random set $A_{X_i}$ in such a way that $\{A_{X_i}\}_{i \in I}$ are stationary, since the edge-preserving property and the degree bound of two are still preserved. E.g., the proof techniques of Theorem~\ref{theorem-nopercolation} are still applicable if $\{A_{X_i}\}_{i \in I}$ is a Boolean model with random radii based on $X=\{X_i\}_{i\in I}$ \cite{MR96}. Instead of connecting $X_i$ to the two furthest neighbors in $X \cap A_{X_i}$ by an edge, one can also connect it to the two nearest neighbors in $X \cap A_{X_i}$ and obtain the same result. We refrain from presenting further details here.
\end{example}
The next example shows that there are graphs defined very similarly to the $f$-2NN~graph for which our methods are not applicable.
\begin{example}[$k_1$-th and $k_2$-th nearest neighbors]
Let $k_1,k_2 \in \N$ such that $ k_1 < k_2$. Similarly to Definition~\ref{defn-fkNN}, the $f$-\emph{$k_1$-th or $k_2$-th}-\emph{nearest neighbor} (\emph{$f$-$(k_1,k_2)$NN}) \emph{graph} $g_{(k_1,k_2),f}(\X)$ is defined as the random graph having vertex set $X$ and for all $i \in I$ and $n \in \{ 1,2\}$ an edge between $X_i$ and $\NN^f_n(X_i,\X)$ whenever $X_i \in \{ \NN^f_{k_1}(\NN^f_n(X_i,\X),\X), \ldots, \NN^f_{k_2}(\NN^f_n(X_i,\X),\X) \}$. In the case $k_1=1$ and $k_2=2$, $g_{(1,2),f}(\X)$ is equal to the $f$-2NN graph $g_{2,f}(\X)$. However, it is easy to see that if $(k_1,k_2) \neq (1,2)$, then the $f$-$(k_1,k_2)$NN graph is in general not edge-preserving in the sense of Definition~\ref{defn-fkNN}. This is a major obstacle for generalizing the proof of Theorem~\ref{theorem-nopercolation} to the case $(k_1,k_2) \neq (1,2)$, despite the fact that many of the proof ingredients of the theorem are still available in this case.
\end{example}
\subsection{Relation of our model to outdegree-one graphs}\label{sec-outdegree1}
In the setting of outdegree-one graphs \cite{CDS20}, one considers directed percolation in a directed graph arising as a deterministic and stationary function of a PPP in $\R^d$, where each vertex has precisely one out-degree. It was shown in \cite{CDS20} that under certain stabilization and looping conditions of the edge-drawing mechanism, this model does not percolate, in the sense that the out-component (or the in-component) of any vertex is almost-surely finite, see also~\cite{H16}. In~\cite{S18}, it was shown that this result is applicable for the example of the \emph{$k$-th nearest neighbor graph}, where the outgoing edge of a vertex points to the $k$-th nearest neighbor of the vertex in the point process. This setting looks rather similar to the one that we are considering but is still different from it, for at least two reasons. First, although it is tempting to think that the \B-\kNN~graph can be obtained as a deterministic transformation of a (stationary) outdegree-one graph satisfying the conditions of \cite{CDS20}, we were not able to find such an outdegree-one graph. Second, the $k$-th nearest neighbor graph is only contained in the \U-\kNN~graph, not in the \B-\kNN~one; in particular, the results of \cite{S18} cannot be derived from our ones.

%% file: VoronoiB2.tex
\begin{tikzpicture}[scale=8] 
 \begin{scope} 
\clip(0.2,0.2) rectangle (0.8,0.8);
\draw[red, thick] (0.8340383503796273,0.7758938836779715)--(0.822841928213028,0.8123889598178787);
\draw[red, thick] (0.16402452764015735,0.3738974712746203)--(0.17418575103719508,0.3127497753447539);
\draw[red, thick] (0.148397993788087,0.28793013543623097)--(0.17418575103719508,0.3127497753447539);
\draw[red, thick] (0.16402452764015735,0.3738974712746203)--(0.18588103912022907,0.40276003273976374);
\draw[red, thick] (0.18588103912022907,0.40276003273976374)--(0.13638721173493218,0.6043853598410402);
\draw[red, thick] (0.0072051335010195755,0.8663122452666777)--(0.07225031953352958,0.7988984371260657);
\draw[red, thick] (0.06964041109724696,0.7722567622110954)--(0.07225031953352958,0.7988984371260657);
\draw[red, thick] (0.822841928213028,0.8123889598178787)--(0.767274141410929,0.8903491846659611);
\draw[red, thick] (0.767274141410929,0.8903491846659611)--(0.6347894819842471,0.8155039577188894);
\draw[red, thick] (0.6347894819842471,0.8155039577188894)--(0.5646470338381138,0.8144491937677845);
\draw[red, thick] (0.148397993788087,0.28793013543623097)--(0.13656917844354743,0.2001593013620731);
\draw[red, thick] (0.19799452036188725,0.2983935547149186)--(0.17418575103719508,0.3127497753447539);
\draw[red, thick] (0.19799452036188725,0.2983935547149186)--(0.39386035741641695,0.3879950389869909);
\draw[red, thick] (0.2614932165060834,0.4483621656048845)--(0.25909724262654277,0.44854516125233673);
\draw[red, thick] (0.2614932165060834,0.4483621656048845)--(0.4157501376828479,0.4287256297513061);
\draw[red, thick] (0.18588103912022907,0.40276003273976374)--(0.25909724262654277,0.44854516125233673);
\draw[red, thick] (0.4157501376828479,0.4287256297513061)--(0.39386035741641695,0.3879950389869909);
\draw[red, thick] (0.2125410814435029,0.932206125956973)--(0.1503709212537891,0.9000494683595954);
\draw[red, thick] (0.2125410814435029,0.932206125956973)--(0.049259228651907655,0.9730210070877323);
\draw[red, thick] (0.1503709212537891,0.9000494683595954)--(0.07950484344558738,0.9060196908690696);
\draw[red, thick] (0.07950484344558738,0.9060196908690696)--(0.03418643383988759,0.9411415638172064);
\draw[red, thick] (0.049259228651907655,0.9730210070877323)--(0.03418643383988759,0.9411415638172064);
\draw[red, thick] (0.3078280725794059,0.9664502126130512)--(0.2125410814435029,0.932206125956973);
\draw[red, thick] (0.0072051335010195755,0.8663122452666777)--(0.07950484344558738,0.9060196908690696);
\draw[red, thick] (0.07225031953352958,0.7988984371260657)--(0.1503709212537891,0.9000494683595954);
\draw[red, thick] (0.20785657217564413,0.7391137029871157)--(0.299736887125702,0.8166411601178638);
\draw[red, thick] (0.20785657217564413,0.7391137029871157)--(0.1626026765847196,0.6149341693854876);
\draw[red, thick] (0.299736887125702,0.8166411601178638)--(0.410824171102973,0.6379558481429204);
\draw[red, thick] (0.1626026765847196,0.6149341693854876)--(0.21771196590924116,0.5929681031496008);
\draw[red, thick] (0.21771196590924116,0.5929681031496008)--(0.3581748076390561,0.6189381158954018);
\draw[red, thick] (0.3581748076390561,0.6189381158954018)--(0.37128451309713617,0.6232642739680077);
\draw[red, thick] (0.410824171102973,0.6379558481429204)--(0.37128451309713617,0.6232642739680077);
\draw[red, thick] (0.3078280725794059,0.9664502126130512)--(0.299736887125702,0.8166411601178638);
\draw[red, thick] (0.06964041109724696,0.7722567622110954)--(0.20785657217564413,0.7391137029871157);
\draw[red, thick] (0.13638721173493218,0.6043853598410402)--(0.1626026765847196,0.6149341693854876);
\draw[red, thick] (0.5646470338381138,0.8144491937677845)--(0.4573568415412147,0.6126174690235537);
\draw[red, thick] (0.4573568415412147,0.6126174690235537)--(0.410824171102973,0.6379558481429204);
\draw[red, thick] (0.25909724262654277,0.44854516125233673)--(0.2553892027923351,0.48539597554160624);
\draw[red, thick] (0.2553892027923351,0.48539597554160624)--(0.21771196590924116,0.5929681031496008);
\draw[red, thick] (0.2553892027923351,0.48539597554160624)--(0.3581748076390561,0.6189381158954018);
\draw[red, thick] (0.2614932165060834,0.4483621656048845)--(0.37128451309713617,0.6232642739680077);
\draw[red, thick] (0.5092587354725757,0.5327225311464738)--(0.4573568415412147,0.6126174690235537);
\draw[red, thick] (0.5092587354725757,0.5327225311464738)--(0.4959932207799768,0.47088604498338843);
\draw[red, thick] (0.4959932207799768,0.47088604498338843)--(0.4157501376828479,0.4287256297513061);
\draw[red, thick] (0.45995291048752523,0.028965958412618953)--(0.33570477350580863,0.0944446314103765);
\draw[red, thick] (0.13656917844354743,0.2001593013620731)--(0.22603564326488124,0.11592720692047373);
\draw[red, thick] (0.679280231333423,0.07685494809943949)--(0.7039360788752047,0.13691025596481718);
\draw[red, thick] (0.679280231333423,0.07685494809943949)--(0.5401162699392226,0.10592861733374764);
\draw[red, thick] (0.7039360788752047,0.13691025596481718)--(0.692984094105217,0.16657177838006781);
\draw[red, thick] (0.692984094105217,0.16657177838006781)--(0.5705030851078778,0.14956730017396283);
\draw[red, thick] (0.5705030851078778,0.14956730017396283)--(0.5401162699392226,0.10592861733374764);
\draw[red, thick] (0.6704261078411188,0.2983594041110991)--(0.692984094105217,0.16657177838006781);
\draw[red, thick] (0.6704261078411188,0.2983594041110991)--(0.609626416494816,0.3156521591180685);
\draw[red, thick] (0.5705030851078778,0.14956730017396283)--(0.609626416494816,0.3156521591180685);
\draw[red, thick] (0.45995291048752523,0.028965958412618953)--(0.4765623254864447,0.08494160121449679);
\draw[red, thick] (0.4765623254864447,0.08494160121449679)--(0.5401162699392226,0.10592861733374764);
\draw[red, thick] (0.38150587398438834,0.343981855497234)--(0.4497840091409947,0.1953790520155802);
\draw[red, thick] (0.38150587398438834,0.343981855497234)--(0.29113000848283804,0.21354918994579553);
\draw[red, thick] (0.297180677272192,0.20845266990432376)--(0.2922904993448189,0.21081962561616713);
\draw[red, thick] (0.297180677272192,0.20845266990432376)--(0.4415822666940058,0.18039376697821286);
\draw[red, thick] (0.4415822666940058,0.18039376697821286)--(0.4497840091409947,0.1953790520155802);
\draw[red, thick] (0.29113000848283804,0.21354918994579553)--(0.2922904993448189,0.21081962561616713);
\draw[red, thick] (0.4959932207799768,0.47088604498338843)--(0.5773694923567108,0.33302816738559193);
\draw[red, thick] (0.39386035741641695,0.3879950389869909)--(0.38150587398438834,0.343981855497234);
\draw[red, thick] (0.5773694923567108,0.33302816738559193)--(0.4497840091409947,0.1953790520155802);
\draw[red, thick] (0.19799452036188725,0.2983935547149186)--(0.29113000848283804,0.21354918994579553);
\draw[red, thick] (0.33570477350580863,0.0944446314103765)--(0.297180677272192,0.20845266990432376);
\draw[red, thick] (0.22603564326488124,0.11592720692047373)--(0.2922904993448189,0.21081962561616713);
\draw[red, thick] (0.4765623254864447,0.08494160121449679)--(0.4415822666940058,0.18039376697821286);
\draw[red, thick] (0.609626416494816,0.3156521591180685)--(0.5773694923567108,0.33302816738559193);
\draw[red, thick] (0.6966078171943033,0.517084675521012)--(0.8394170635303543,0.5386328082928883);
\draw[red, thick] (0.6966078171943033,0.517084675521012)--(0.7629987978835259,0.4097388057906942);
\draw[red, thick] (0.8394170635303543,0.5386328082928883)--(0.8199448304662863,0.43374447248530434);
\draw[red, thick] (0.7629987978835259,0.4097388057906942)--(0.8199448304662863,0.43374447248530434);
\draw[red, thick] (0.5092587354725757,0.5327225311464738)--(0.5854609328882098,0.5563630707966651);
\draw[red, thick] (0.6704261078411188,0.2983594041110991)--(0.706255398851686,0.31142037155747987);
\draw[red, thick] (0.706255398851686,0.31142037155747987)--(0.7629987978835259,0.4097388057906942);
\draw[red, thick] (0.6966078171943033,0.517084675521012)--(0.6791059840045244,0.528146008831672);
\draw[red, thick] (0.5854609328882098,0.5563630707966651)--(0.6791059840045244,0.528146008831672);
\draw[red, thick] (0.8394170635303543,0.5386328082928883)--(0.8421316549547286,0.5424862865901321);
\draw[red, thick] (0.6791059840045244,0.528146008831672)--(0.7961709115289796,0.6183095407912882);
\draw[red, thick] (0.8421316549547286,0.5424862865901321)--(0.7961709115289796,0.6183095407912882);
\draw[red, thick] (0.8180743974620243,0.1360181565046783)--(0.7847271846057482,0.2833695475015336);
\draw[red, thick] (0.7039360788752047,0.13691025596481718)--(0.8180743974620243,0.1360181565046783);
\draw[red, thick] (0.706255398851686,0.31142037155747987)--(0.7847271846057482,0.2833695475015336);
\draw[red, thick] (0.8340383503796273,0.7758938836779715)--(0.8378339831936373,0.7542743252554434);
\draw[red, thick] (0.663654233906528,0.6697057286556454)--(0.6320097831139766,0.6120438657752406);
\draw[red, thick] (0.663654233906528,0.6697057286556454)--(0.8042869452731249,0.7046892265869095);
\draw[red, thick] (0.6320097831139766,0.6120438657752406)--(0.7954628450760505,0.6220772182985786);
\draw[red, thick] (0.8042869452731249,0.7046892265869095)--(0.7954628450760505,0.6220772182985786);
\draw[red, thick] (0.6347894819842471,0.8155039577188894)--(0.663654233906528,0.6697057286556454);
\draw[red, thick] (0.5854609328882098,0.5563630707966651)--(0.6320097831139766,0.6120438657752406);
\draw[red, thick] (0.8378339831936373,0.7542743252554434)--(0.8042869452731249,0.7046892265869095);
\draw[red, thick] (0.7961709115289796,0.6183095407912882)--(0.7954628450760505,0.6220772182985786);
\draw[black, very thick] (0.05789114625497006,0.7591447714783306)--(0.025416319820836576,0.7229033871050434);
\draw[black, very thick] (0.05789114625497006,0.7591447714783306)--(0.07059719372583992,0.7820235001551941);
\draw[black, very thick] (0.025416319820836576,0.7229033871050434)--(0.07059719372583992,0.7820235001551941);
\draw[black, very thick] (0.816871669547315,0.8207650844535032)--(0.797244196123064,0.8483019422161482);
\draw[black, very thick] (0.816871669547315,0.8207650844535032)--(0.84050884537445,0.7523984527019172);
\draw[black, very thick] (0.797244196123064,0.8483019422161482)--(0.7024556644193323,0.853730955492912);
\draw[black, very thick] (0.5150163133378375,0.8659430761214151)--(0.47164836244225017,0.9109390813944592);
\draw[black, very thick] (0.5150163133378375,0.8659430761214151)--(0.5337666871494476,0.7563578295951302);
\draw[black, very thick] (0.29224699803261267,0.34151062925235404)--(0.25616296056720894,0.325003496798576);
\draw[black, very thick] (0.3693375759616444,0.43463383746958717)--(0.37181151999371226,0.4343189103317002);
\draw[black, very thick] (0.37181151999371226,0.4343189103317002)--(0.4100600707523283,0.4294499616155625);
\draw[black, very thick] (0.4100600707523283,0.4294499616155625)--(0.4169904242473743,0.42937728712224627);
\draw[black, very thick] (0.25532816389877266,0.4461882000106892)--(0.2585381570265892,0.4541014013513852);
\draw[black, very thick] (0.15937905866236854,0.9047088031490487)--(0.1590067661594326,0.9455878746166037);
\draw[black, very thick] (0.15937905866236854,0.9047088031490487)--(0.17361789056520271,0.9419355930212361);
\draw[black, very thick] (0.1590067661594326,0.9455878746166037)--(0.17361789056520271,0.9419355930212361);
\draw[black, very thick] (0.2921345429365123,0.9608102969999692)--(0.30700725570527704,0.9512527117166875);
\draw[black, very thick] (0.2921345429365123,0.9608102969999692)--(0.304643617111923,0.9074897229510144);
\draw[black, very thick] (0.2526856156454257,0.7769398874936125)--(0.27892638396565794,0.7990815186159101);
\draw[black, very thick] (0.19387030084463816,0.7007344954508157)--(0.18563030063821997,0.7444433737589387);
\draw[black, very thick] (0.37293603650876134,0.6988993954844914)--(0.3982187058250188,0.6582318996866904);
\draw[black, very thick] (0.18134071416677522,0.6074653564194161)--(0.1726754312359303,0.6109192601286578);
\draw[black, very thick] (0.34509732830012835,0.6165202358288133)--(0.36683993239497326,0.6217975775006084);
\draw[black, very thick] (0.34509732830012835,0.6165202358288133)--(0.3498485323164778,0.5891158637590141);
\draw[black, very thick] (0.232149332632162,0.5956374112256204)--(0.23156648132676858,0.5534121318594152);
\draw[black, very thick] (0.232149332632162,0.5956374112256204)--(0.2277040583383429,0.5644397202114);
\draw[black, very thick] (0.3010446550260852,0.8408546046323935)--(0.29981130217595947,0.8180189617063004);
\draw[black, very thick] (0.3010446550260852,0.8408546046323935)--(0.30055791219730454,0.8318425158098488);
\draw[black, very thick] (0.30700725570527704,0.9512527117166875)--(0.304643617111923,0.9074897229510144);
\draw[black, very thick] (0.29981130217595947,0.8180189617063004)--(0.30055791219730454,0.8318425158098488);
\draw[black, very thick] (0.5337666871494476,0.7563578295951302)--(0.4828113417610102,0.6605018615470278);
\draw[black, very thick] (0.23156648132676858,0.5534121318594152)--(0.2277040583383429,0.5644397202114);
\draw[black, very thick] (0.26648526291374264,0.4998123096218199)--(0.293302842223799,0.5346545129197736);
\draw[black, very thick] (0.26648526291374264,0.4998123096218199)--(0.2919655586376268,0.49690587581737233);
\draw[black, very thick] (0.293302842223799,0.5346545129197736)--(0.2919655586376268,0.49690587581737233);
\draw[black, very thick] (0.3498485323164778,0.5891158637590141)--(0.34045815137782165,0.574156595807007);
\draw[black, very thick] (0.5084688879484911,0.5290406992216823)--(0.49024128124118205,0.4678639258176544);
\draw[black, very thick] (0.38355244884048856,0.06922894316707184)--(0.37508174221498614,0.07369299906895496);
\draw[black, very thick] (0.44145691925859737,0.03871333161972445)--(0.4744519247018341,0.0778293073440714);
\draw[black, very thick] (0.37508174221498614,0.07369299906895496)--(0.35676393151346963,0.08334647127462609);
\draw[black, very thick] (0.35676393151346963,0.08334647127462609)--(0.33402090497369247,0.09314723041436684);
\draw[black, very thick] (0.27052606386950184,0.04422519848677735)--(0.25415869288763704,0.031614330708849434);
\draw[black, very thick] (0.19012226347068706,-0.017724988674975257)--(0.17858300600531926,-0.004973681536212152);
\draw[black, very thick] (0.15362496192050873,0.1841013922061019)--(0.153557661738525,0.18416475488949516);
\draw[black, very thick] (0.15362496192050873,0.1841013922061019)--(0.14846174817708205,0.18896252434680957);
\draw[black, very thick] (0.153557661738525,0.18416475488949516)--(0.14846174817708205,0.18896252434680957);
\draw[black, very thick] (0.1872786858015116,0.01718136555514127)--(0.18417778144474048,0.00928081221598432);
\draw[black, very thick] (0.18417778144474048,0.00928081221598432)--(0.17858300600531926,-0.004973681536212152);
\draw[black, very thick] (0.694675762855622,0.11435450611503778)--(0.6796560855289289,0.2444362909933815);
\draw[black, very thick] (0.694675762855622,0.11435450611503778)--(0.8118874384593696,0.1633564836689559);
\draw[black, very thick] (0.4744519247018341,0.0778293073440714)--(0.5263999923459467,0.10139917795079811);
\draw[black, very thick] (0.44625527759791844,0.20305910066999244)--(0.45162040934298203,0.19736030299255442);
\draw[black, very thick] (0.35675114120119855,0.3082552180370413)--(0.37932348476575223,0.34083217785400804);
\draw[black, very thick] (0.35675114120119855,0.3082552180370413)--(0.3849643126762702,0.3563026374685184);
\draw[black, very thick] (0.37932348476575223,0.34083217785400804)--(0.3849643126762702,0.3563026374685184);
\draw[black, very thick] (0.43025304052789864,0.18259516696622943)--(0.45162040934298203,0.19736030299255442);
\draw[black, very thick] (0.43025304052789864,0.18259516696622943)--(0.44500451138548675,0.17105528445839943);
\draw[black, very thick] (0.5646392290876473,0.35459424606590373)--(0.5570056583155196,0.3675261222104846);
\draw[black, very thick] (0.5646392290876473,0.35459424606590373)--(0.5464741713477421,0.2996958987076048);
\draw[black, very thick] (0.5570056583155196,0.3675261222104846)--(0.5326754741324017,0.40874338970209434);
\draw[black, very thick] (0.48139088714400735,0.22947900357015216)--(0.4952009553583087,0.24437837632781984);
\draw[black, very thick] (0.7221119059162853,0.4758477701896109)--(0.7711711320631728,0.41318386283041797);
\draw[black, very thick] (0.8272393746070328,0.4730369662060027)--(0.7711711320631728,0.41318386283041797);
\draw[black, very thick] (0.5800834392674282,0.554694787525681)--(0.6188652259059673,0.5462977123528099);
\draw[black, very thick] (0.5800834392674282,0.554694787525681)--(0.6361928686502474,0.5410765590383083);
\draw[black, very thick] (0.6188652259059673,0.5462977123528099)--(0.6361928686502474,0.5410765590383083);
\draw[black, very thick] (0.7880511358690263,0.6120556809230696)--(0.7843733421016718,0.6092230402432683);
\draw[black, very thick] (0.7880511358690263,0.6120556809230696)--(0.7851149531728918,0.6097942298808102);
\draw[black, very thick] (0.7843733421016718,0.6092230402432683)--(0.7851149531728918,0.6097942298808102);
\draw[black, very thick] (0.7940377748297177,0.6166665966317921)--(0.7860546895085898,0.6214997109685254);
\draw[black, very thick] (0.8851632825205141,0.5493531197369946)--(0.9421973439307684,0.5584544111718442);
\draw[black, very thick] (1.0241664694608577,0.57153475091379)--(1.0378006206639054,0.5737104399756852);
\draw[black, very thick] (1.0605948498489768,0.5980528201790608)--(1.0418135837704467,0.6112240638453353);
\draw[black, very thick] (1.0605948498489768,0.5980528201790608)--(1.0588598410834567,0.5992695764466222);
\draw[black, very thick] (0.9625363197915181,0.6668209632607429)--(0.9168260577530807,0.6988774282956411);
\draw[black, very thick] (0.9625363197915181,0.6668209632607429)--(0.9338936187195976,0.6869080007107792);
\draw[black, very thick] (1.0418135837704467,0.6112240638453353)--(1.0588598410834567,0.5992695764466222);
\draw[black, very thick] (0.9168260577530807,0.6988774282956411)--(0.9338936187195976,0.6869080007107792);
\draw[black, very thick] (0.642114287806097,0.6304560849618224)--(0.6273026388087081,0.6064132747783626);
\draw[black, very thick] (0.7630519194473339,0.6944316881617858)--(0.7580659591276837,0.6931913911438812);
\draw[black, very thick] (0.7630519194473339,0.6944316881617858)--(0.8032960357300718,0.6954122436021932);
\draw[black, very thick] (0.7580659591276837,0.6931913911438812)--(0.8032960357300718,0.6954122436021932);
\draw[black, very thick] (0.7097088513509334,0.6811621785700015)--(0.6984784752351365,0.6783685337927906);
\draw[black, very thick] (0.6984784752351365,0.6783685337927906)--(0.6595268144771232,0.6905536626697233);
\fill[blue] (0.1813061501361128,0.42139697266953957) circle (0.1pt);\draw(0.1813061501361128,0.42139697266953957) circle (0.1pt);
\fill[blue] (0.14872563829292249,0.5541217325950919) circle (0.1pt);\draw(0.14872563829292249,0.5541217325950919) circle (0.1pt);
\fill[blue] (0.05789114625497006,0.7591447714783306) circle (0.1pt);\draw(0.05789114625497006,0.7591447714783306) circle (0.1pt);
\fill[blue] (0.025416319820836576,0.7229033871050434) circle (0.1pt);\draw(0.025416319820836576,0.7229033871050434) circle (0.1pt);
\fill[blue] (0.07059719372583992,0.7820235001551941) circle (0.1pt);\draw(0.07059719372583992,0.7820235001551941) circle (0.1pt);
\fill[blue] (0.816871669547315,0.8207650844535032) circle (0.1pt);\draw(0.816871669547315,0.8207650844535032) circle (0.1pt);
\fill[blue] (0.797244196123064,0.8483019422161482) circle (0.1pt);\draw(0.797244196123064,0.8483019422161482) circle (0.1pt);
\fill[blue] (0.5150163133378375,0.8659430761214151) circle (0.1pt);\draw(0.5150163133378375,0.8659430761214151) circle (0.1pt);
\fill[blue] (0.47164836244225017,0.9109390813944592) circle (0.1pt);\draw(0.47164836244225017,0.9109390813944592) circle (0.1pt);
\fill[blue] (0.38653207447931925,0.9992506770062238) circle (0.1pt);\draw(0.38653207447931925,0.9992506770062238) circle (0.1pt);
\fill[blue] (0.7024556644193323,0.853730955492912) circle (0.1pt);\draw(0.7024556644193323,0.853730955492912) circle (0.1pt);
\fill[blue] (0.29224699803261267,0.34151062925235404) circle (0.1pt);\draw(0.29224699803261267,0.34151062925235404) circle (0.1pt);
\fill[blue] (0.25616296056720894,0.325003496798576) circle (0.1pt);\draw(0.25616296056720894,0.325003496798576) circle (0.1pt);
\fill[blue] (0.3693375759616444,0.43463383746958717) circle (0.1pt);\draw(0.3693375759616444,0.43463383746958717) circle (0.1pt);
\fill[blue] (0.37181151999371226,0.4343189103317002) circle (0.1pt);\draw(0.37181151999371226,0.4343189103317002) circle (0.1pt);
\fill[blue] (0.4100600707523283,0.4294499616155625) circle (0.1pt);\draw(0.4100600707523283,0.4294499616155625) circle (0.1pt);
\fill[blue] (0.25532816389877266,0.4461882000106892) circle (0.1pt);\draw(0.25532816389877266,0.4461882000106892) circle (0.1pt);
\fill[blue] (0.15937905866236854,0.9047088031490487) circle (0.1pt);\draw(0.15937905866236854,0.9047088031490487) circle (0.1pt);
\fill[blue] (0.1590067661594326,0.9455878746166037) circle (0.1pt);\draw(0.1590067661594326,0.9455878746166037) circle (0.1pt);
\fill[blue] (0.17361789056520271,0.9419355930212361) circle (0.1pt);\draw(0.17361789056520271,0.9419355930212361) circle (0.1pt);
\fill[blue] (0.06284152846679965,0.9189337987707326) circle (0.1pt);\draw(0.06284152846679965,0.9189337987707326) circle (0.1pt);
\fill[blue] (0.2921345429365123,0.9608102969999692) circle (0.1pt);\draw(0.2921345429365123,0.9608102969999692) circle (0.1pt);
\fill[blue] (0.2526856156454257,0.7769398874936125) circle (0.1pt);\draw(0.2526856156454257,0.7769398874936125) circle (0.1pt);
\fill[blue] (0.27892638396565794,0.7990815186159101) circle (0.1pt);\draw(0.27892638396565794,0.7990815186159101) circle (0.1pt);
\fill[blue] (0.19387030084463816,0.7007344954508157) circle (0.1pt);\draw(0.19387030084463816,0.7007344954508157) circle (0.1pt);
\fill[blue] (0.37293603650876134,0.6988993954844914) circle (0.1pt);\draw(0.37293603650876134,0.6988993954844914) circle (0.1pt);
\fill[blue] (0.3982187058250188,0.6582318996866904) circle (0.1pt);\draw(0.3982187058250188,0.6582318996866904) circle (0.1pt);
\fill[blue] (0.18134071416677522,0.6074653564194161) circle (0.1pt);\draw(0.18134071416677522,0.6074653564194161) circle (0.1pt);
\fill[blue] (0.1726754312359303,0.6109192601286578) circle (0.1pt);\draw(0.1726754312359303,0.6109192601286578) circle (0.1pt);
\fill[blue] (0.34509732830012835,0.6165202358288133) circle (0.1pt);\draw(0.34509732830012835,0.6165202358288133) circle (0.1pt);
\fill[blue] (0.232149332632162,0.5956374112256204) circle (0.1pt);\draw(0.232149332632162,0.5956374112256204) circle (0.1pt);
\fill[blue] (0.36683993239497326,0.6217975775006084) circle (0.1pt);\draw(0.36683993239497326,0.6217975775006084) circle (0.1pt);
\fill[blue] (0.3010446550260852,0.8408546046323935) circle (0.1pt);\draw(0.3010446550260852,0.8408546046323935) circle (0.1pt);
\fill[blue] (0.30700725570527704,0.9512527117166875) circle (0.1pt);\draw(0.30700725570527704,0.9512527117166875) circle (0.1pt);
\fill[blue] (0.304643617111923,0.9074897229510144) circle (0.1pt);\draw(0.304643617111923,0.9074897229510144) circle (0.1pt);
\fill[blue] (0.29981130217595947,0.8180189617063004) circle (0.1pt);\draw(0.29981130217595947,0.8180189617063004) circle (0.1pt);
\fill[blue] (0.30055791219730454,0.8318425158098488) circle (0.1pt);\draw(0.30055791219730454,0.8318425158098488) circle (0.1pt);
\fill[blue] (0.18563030063821997,0.7444433737589387) circle (0.1pt);\draw(0.18563030063821997,0.7444433737589387) circle (0.1pt);
\fill[blue] (0.5337666871494476,0.7563578295951302) circle (0.1pt);\draw(0.5337666871494476,0.7563578295951302) circle (0.1pt);
\fill[blue] (0.4828113417610102,0.6605018615470278) circle (0.1pt);\draw(0.4828113417610102,0.6605018615470278) circle (0.1pt);
\fill[blue] (0.2585381570265892,0.4541014013513852) circle (0.1pt);\draw(0.2585381570265892,0.4541014013513852) circle (0.1pt);
\fill[blue] (0.23156648132676858,0.5534121318594152) circle (0.1pt);\draw(0.23156648132676858,0.5534121318594152) circle (0.1pt);
\fill[blue] (0.2277040583383429,0.5644397202114) circle (0.1pt);\draw(0.2277040583383429,0.5644397202114) circle (0.1pt);
\fill[blue] (0.26648526291374264,0.4998123096218199) circle (0.1pt);\draw(0.26648526291374264,0.4998123096218199) circle (0.1pt);
\fill[blue] (0.293302842223799,0.5346545129197736) circle (0.1pt);\draw(0.293302842223799,0.5346545129197736) circle (0.1pt);
\fill[blue] (0.3498485323164778,0.5891158637590141) circle (0.1pt);\draw(0.3498485323164778,0.5891158637590141) circle (0.1pt);
\fill[blue] (0.2919655586376268,0.49690587581737233) circle (0.1pt);\draw(0.2919655586376268,0.49690587581737233) circle (0.1pt);
\fill[blue] (0.34045815137782165,0.574156595807007) circle (0.1pt);\draw(0.34045815137782165,0.574156595807007) circle (0.1pt);
\fill[blue] (0.5084688879484911,0.5290406992216823) circle (0.1pt);\draw(0.5084688879484911,0.5290406992216823) circle (0.1pt);
\fill[blue] (0.4169904242473743,0.42937728712224627) circle (0.1pt);\draw(0.4169904242473743,0.42937728712224627) circle (0.1pt);
\fill[blue] (0.49024128124118205,0.4678639258176544) circle (0.1pt);\draw(0.49024128124118205,0.4678639258176544) circle (0.1pt);
\fill[blue] (0.38355244884048856,0.06922894316707184) circle (0.1pt);\draw(0.38355244884048856,0.06922894316707184) circle (0.1pt);
\fill[blue] (0.44145691925859737,0.03871333161972445) circle (0.1pt);\draw(0.44145691925859737,0.03871333161972445) circle (0.1pt);
\fill[blue] (0.37508174221498614,0.07369299906895496) circle (0.1pt);\draw(0.37508174221498614,0.07369299906895496) circle (0.1pt);
\fill[blue] (0.35676393151346963,0.08334647127462609) circle (0.1pt);\draw(0.35676393151346963,0.08334647127462609) circle (0.1pt);
\fill[blue] (0.27052606386950184,0.04422519848677735) circle (0.1pt);\draw(0.27052606386950184,0.04422519848677735) circle (0.1pt);
\fill[blue] (0.33402090497369247,0.09314723041436684) circle (0.1pt);\draw(0.33402090497369247,0.09314723041436684) circle (0.1pt);
\fill[blue] (0.19012226347068706,-0.017724988674975257) circle (0.1pt);\draw(0.19012226347068706,-0.017724988674975257) circle (0.1pt);
\fill[blue] (0.25415869288763704,0.031614330708849434) circle (0.1pt);\draw(0.25415869288763704,0.031614330708849434) circle (0.1pt);
\fill[blue] (0.19976901885816142,0.14065706287384386) circle (0.1pt);\draw(0.19976901885816142,0.14065706287384386) circle (0.1pt);
\fill[blue] (0.15362496192050873,0.1841013922061019) circle (0.1pt);\draw(0.15362496192050873,0.1841013922061019) circle (0.1pt);
\fill[blue] (0.153557661738525,0.18416475488949516) circle (0.1pt);\draw(0.153557661738525,0.18416475488949516) circle (0.1pt);
\fill[blue] (0.14846174817708205,0.18896252434680957) circle (0.1pt);\draw(0.14846174817708205,0.18896252434680957) circle (0.1pt);
\fill[blue] (0.1872786858015116,0.01718136555514127) circle (0.1pt);\draw(0.1872786858015116,0.01718136555514127) circle (0.1pt);
\fill[blue] (0.18417778144474048,0.00928081221598432) circle (0.1pt);\draw(0.18417778144474048,0.00928081221598432) circle (0.1pt);
\fill[blue] (0.17858300600531926,-0.004973681536212152) circle (0.1pt);\draw(0.17858300600531926,-0.004973681536212152) circle (0.1pt);
\fill[blue] (0.694675762855622,0.11435450611503778) circle (0.1pt);\draw(0.694675762855622,0.11435450611503778) circle (0.1pt);
\fill[blue] (0.6796560855289289,0.2444362909933815) circle (0.1pt);\draw(0.6796560855289289,0.2444362909933815) circle (0.1pt);
\fill[blue] (0.4744519247018341,0.0778293073440714) circle (0.1pt);\draw(0.4744519247018341,0.0778293073440714) circle (0.1pt);
\fill[blue] (0.5263999923459467,0.10139917795079811) circle (0.1pt);\draw(0.5263999923459467,0.10139917795079811) circle (0.1pt);
\fill[blue] (0.44625527759791844,0.20305910066999244) circle (0.1pt);\draw(0.44625527759791844,0.20305910066999244) circle (0.1pt);
\fill[blue] (0.35675114120119855,0.3082552180370413) circle (0.1pt);\draw(0.35675114120119855,0.3082552180370413) circle (0.1pt);
\fill[blue] (0.37932348476575223,0.34083217785400804) circle (0.1pt);\draw(0.37932348476575223,0.34083217785400804) circle (0.1pt);
\fill[blue] (0.3781940958333602,0.1927108234591293) circle (0.1pt);\draw(0.3781940958333602,0.1927108234591293) circle (0.1pt);
\fill[blue] (0.43025304052789864,0.18259516696622943) circle (0.1pt);\draw(0.43025304052789864,0.18259516696622943) circle (0.1pt);
\fill[blue] (0.5646392290876473,0.35459424606590373) circle (0.1pt);\draw(0.5646392290876473,0.35459424606590373) circle (0.1pt);
\fill[blue] (0.5570056583155196,0.3675261222104846) circle (0.1pt);\draw(0.5570056583155196,0.3675261222104846) circle (0.1pt);
\fill[blue] (0.5326754741324017,0.40874338970209434) circle (0.1pt);\draw(0.5326754741324017,0.40874338970209434) circle (0.1pt);
\fill[blue] (0.3849643126762702,0.3563026374685184) circle (0.1pt);\draw(0.3849643126762702,0.3563026374685184) circle (0.1pt);
\fill[blue] (0.48139088714400735,0.22947900357015216) circle (0.1pt);\draw(0.48139088714400735,0.22947900357015216) circle (0.1pt);
\fill[blue] (0.4952009553583087,0.24437837632781984) circle (0.1pt);\draw(0.4952009553583087,0.24437837632781984) circle (0.1pt);
\fill[blue] (0.5464741713477421,0.2996958987076048) circle (0.1pt);\draw(0.5464741713477421,0.2996958987076048) circle (0.1pt);
\fill[blue] (0.45162040934298203,0.19736030299255442) circle (0.1pt);\draw(0.45162040934298203,0.19736030299255442) circle (0.1pt);
\fill[blue] (0.44500451138548675,0.17105528445839943) circle (0.1pt);\draw(0.44500451138548675,0.17105528445839943) circle (0.1pt);
\fill[blue] (0.45877901890888895,0.13346796697448662) circle (0.1pt);\draw(0.45877901890888895,0.13346796697448662) circle (0.1pt);
\fill[blue] (0.7221119059162853,0.4758477701896109) circle (0.1pt);\draw(0.7221119059162853,0.4758477701896109) circle (0.1pt);
\fill[blue] (0.8272393746070328,0.4730369662060027) circle (0.1pt);\draw(0.8272393746070328,0.4730369662060027) circle (0.1pt);
\fill[blue] (0.7711711320631728,0.41318386283041797) circle (0.1pt);\draw(0.7711711320631728,0.41318386283041797) circle (0.1pt);
\fill[blue] (0.5800834392674282,0.554694787525681) circle (0.1pt);\draw(0.5800834392674282,0.554694787525681) circle (0.1pt);
\fill[blue] (0.6188652259059673,0.5462977123528099) circle (0.1pt);\draw(0.6188652259059673,0.5462977123528099) circle (0.1pt);
\fill[blue] (0.6361928686502474,0.5410765590383083) circle (0.1pt);\draw(0.6361928686502474,0.5410765590383083) circle (0.1pt);
\fill[blue] (0.7880511358690263,0.6120556809230696) circle (0.1pt);\draw(0.7880511358690263,0.6120556809230696) circle (0.1pt);
\fill[blue] (0.7843733421016718,0.6092230402432683) circle (0.1pt);\draw(0.7843733421016718,0.6092230402432683) circle (0.1pt);
\fill[blue] (0.7940377748297177,0.6166665966317921) circle (0.1pt);\draw(0.7940377748297177,0.6166665966317921) circle (0.1pt);
\fill[blue] (0.7851149531728918,0.6097942298808102) circle (0.1pt);\draw(0.7851149531728918,0.6097942298808102) circle (0.1pt);
\fill[blue] (0.8108312686319776,0.5941237704915503) circle (0.1pt);\draw(0.8108312686319776,0.5941237704915503) circle (0.1pt);
\fill[blue] (0.8118874384593696,0.1633564836689559) circle (0.1pt);\draw(0.8118874384593696,0.1633564836689559) circle (0.1pt);
\fill[blue] (0.8851632825205141,0.5493531197369946) circle (0.1pt);\draw(0.8851632825205141,0.5493531197369946) circle (0.1pt);
\fill[blue] (1.0241664694608577,0.57153475091379) circle (0.1pt);\draw(1.0241664694608577,0.57153475091379) circle (0.1pt);
\fill[blue] (1.0378006206639054,0.5737104399756852) circle (0.1pt);\draw(1.0378006206639054,0.5737104399756852) circle (0.1pt);
\fill[blue] (1.0823433206558173,0.580818404653238) circle (0.1pt);\draw(1.0823433206558173,0.580818404653238) circle (0.1pt);
\fill[blue] (0.9421973439307684,0.5584544111718442) circle (0.1pt);\draw(0.9421973439307684,0.5584544111718442) circle (0.1pt);
\fill[blue] (1.0605948498489768,0.5980528201790608) circle (0.1pt);\draw(1.0605948498489768,0.5980528201790608) circle (0.1pt);
\fill[blue] (0.9625363197915181,0.6668209632607429) circle (0.1pt);\draw(0.9625363197915181,0.6668209632607429) circle (0.1pt);
\fill[blue] (1.0418135837704467,0.6112240638453353) circle (0.1pt);\draw(1.0418135837704467,0.6112240638453353) circle (0.1pt);
\fill[blue] (0.9168260577530807,0.6988774282956411) circle (0.1pt);\draw(0.9168260577530807,0.6988774282956411) circle (0.1pt);
\fill[blue] (1.0588598410834567,0.5992695764466222) circle (0.1pt);\draw(1.0588598410834567,0.5992695764466222) circle (0.1pt);
\fill[blue] (0.9338936187195976,0.6869080007107792) circle (0.1pt);\draw(0.9338936187195976,0.6869080007107792) circle (0.1pt);
\fill[blue] (0.84050884537445,0.7523984527019172) circle (0.1pt);\draw(0.84050884537445,0.7523984527019172) circle (0.1pt);
\fill[blue] (0.642114287806097,0.6304560849618224) circle (0.1pt);\draw(0.642114287806097,0.6304560849618224) circle (0.1pt);
\fill[blue] (0.7630519194473339,0.6944316881617858) circle (0.1pt);\draw(0.7630519194473339,0.6944316881617858) circle (0.1pt);
\fill[blue] (0.7580659591276837,0.6931913911438812) circle (0.1pt);\draw(0.7580659591276837,0.6931913911438812) circle (0.1pt);
\fill[blue] (0.7097088513509334,0.6811621785700015) circle (0.1pt);\draw(0.7097088513509334,0.6811621785700015) circle (0.1pt);
\fill[blue] (0.6984784752351365,0.6783685337927906) circle (0.1pt);\draw(0.6984784752351365,0.6783685337927906) circle (0.1pt);
\fill[blue] (0.714648332376192,0.6171165250868014) circle (0.1pt);\draw(0.714648332376192,0.6171165250868014) circle (0.1pt);
\fill[blue] (0.7860546895085898,0.6214997109685254) circle (0.1pt);\draw(0.7860546895085898,0.6214997109685254) circle (0.1pt);
\fill[blue] (0.8032960357300718,0.6954122436021932) circle (0.1pt);\draw(0.8032960357300718,0.6954122436021932) circle (0.1pt);
\fill[blue] (0.6595268144771232,0.6905536626697233) circle (0.1pt);\draw(0.6595268144771232,0.6905536626697233) circle (0.1pt);
\fill[blue] (0.6273026388087081,0.6064132747783626) circle (0.1pt);\draw(0.6273026388087081,0.6064132747783626) circle (0.1pt);
 \end{scope} 
 \draw (0.2,0.2) rectangle (0.8,0.8); 
 \end{tikzpicture} 

%% file: VoronoiB4.tex
\begin{tikzpicture}[scale=8] 
 \begin{scope} 
\clip(0.2,0.2) rectangle (0.8,0.8);
\draw[red, thick] (0.8340383503796273,0.7758938836779715)--(0.822841928213028,0.8123889598178787);
\draw[red, thick] (0.16402452764015735,0.3738974712746203)--(0.17418575103719508,0.3127497753447539);
\draw[red, thick] (0.148397993788087,0.28793013543623097)--(0.17418575103719508,0.3127497753447539);
\draw[red, thick] (0.16402452764015735,0.3738974712746203)--(0.18588103912022907,0.40276003273976374);
\draw[red, thick] (0.18588103912022907,0.40276003273976374)--(0.13638721173493218,0.6043853598410402);
\draw[red, thick] (0.0072051335010195755,0.8663122452666777)--(0.07225031953352958,0.7988984371260657);
\draw[red, thick] (0.06964041109724696,0.7722567622110954)--(0.07225031953352958,0.7988984371260657);
\draw[red, thick] (0.822841928213028,0.8123889598178787)--(0.767274141410929,0.8903491846659611);
\draw[red, thick] (0.767274141410929,0.8903491846659611)--(0.6347894819842471,0.8155039577188894);
\draw[red, thick] (0.6347894819842471,0.8155039577188894)--(0.5646470338381138,0.8144491937677845);
\draw[red, thick] (0.148397993788087,0.28793013543623097)--(0.13656917844354743,0.2001593013620731);
\draw[red, thick] (0.19799452036188725,0.2983935547149186)--(0.17418575103719508,0.3127497753447539);
\draw[red, thick] (0.19799452036188725,0.2983935547149186)--(0.39386035741641695,0.3879950389869909);
\draw[red, thick] (0.2614932165060834,0.4483621656048845)--(0.25909724262654277,0.44854516125233673);
\draw[red, thick] (0.2614932165060834,0.4483621656048845)--(0.4157501376828479,0.4287256297513061);
\draw[red, thick] (0.18588103912022907,0.40276003273976374)--(0.25909724262654277,0.44854516125233673);
\draw[red, thick] (0.4157501376828479,0.4287256297513061)--(0.39386035741641695,0.3879950389869909);
\draw[red, thick] (0.2125410814435029,0.932206125956973)--(0.1503709212537891,0.9000494683595954);
\draw[red, thick] (0.2125410814435029,0.932206125956973)--(0.049259228651907655,0.9730210070877323);
\draw[red, thick] (0.1503709212537891,0.9000494683595954)--(0.07950484344558738,0.9060196908690696);
\draw[red, thick] (0.07950484344558738,0.9060196908690696)--(0.03418643383988759,0.9411415638172064);
\draw[red, thick] (0.049259228651907655,0.9730210070877323)--(0.03418643383988759,0.9411415638172064);
\draw[red, thick] (0.3078280725794059,0.9664502126130512)--(0.2125410814435029,0.932206125956973);
\draw[red, thick] (0.0072051335010195755,0.8663122452666777)--(0.07950484344558738,0.9060196908690696);
\draw[red, thick] (0.07225031953352958,0.7988984371260657)--(0.1503709212537891,0.9000494683595954);
\draw[red, thick] (0.20785657217564413,0.7391137029871157)--(0.299736887125702,0.8166411601178638);
\draw[red, thick] (0.20785657217564413,0.7391137029871157)--(0.1626026765847196,0.6149341693854876);
\draw[red, thick] (0.299736887125702,0.8166411601178638)--(0.410824171102973,0.6379558481429204);
\draw[red, thick] (0.1626026765847196,0.6149341693854876)--(0.21771196590924116,0.5929681031496008);
\draw[red, thick] (0.21771196590924116,0.5929681031496008)--(0.3581748076390561,0.6189381158954018);
\draw[red, thick] (0.3581748076390561,0.6189381158954018)--(0.37128451309713617,0.6232642739680077);
\draw[red, thick] (0.410824171102973,0.6379558481429204)--(0.37128451309713617,0.6232642739680077);
\draw[red, thick] (0.3078280725794059,0.9664502126130512)--(0.299736887125702,0.8166411601178638);
\draw[red, thick] (0.06964041109724696,0.7722567622110954)--(0.20785657217564413,0.7391137029871157);
\draw[red, thick] (0.13638721173493218,0.6043853598410402)--(0.1626026765847196,0.6149341693854876);
\draw[red, thick] (0.5646470338381138,0.8144491937677845)--(0.4573568415412147,0.6126174690235537);
\draw[red, thick] (0.4573568415412147,0.6126174690235537)--(0.410824171102973,0.6379558481429204);
\draw[red, thick] (0.25909724262654277,0.44854516125233673)--(0.2553892027923351,0.48539597554160624);
\draw[red, thick] (0.2553892027923351,0.48539597554160624)--(0.21771196590924116,0.5929681031496008);
\draw[red, thick] (0.2553892027923351,0.48539597554160624)--(0.3581748076390561,0.6189381158954018);
\draw[red, thick] (0.2614932165060834,0.4483621656048845)--(0.37128451309713617,0.6232642739680077);
\draw[red, thick] (0.5092587354725757,0.5327225311464738)--(0.4573568415412147,0.6126174690235537);
\draw[red, thick] (0.5092587354725757,0.5327225311464738)--(0.4959932207799768,0.47088604498338843);
\draw[red, thick] (0.4959932207799768,0.47088604498338843)--(0.4157501376828479,0.4287256297513061);
\draw[red, thick] (0.45995291048752523,0.028965958412618953)--(0.33570477350580863,0.0944446314103765);
\draw[red, thick] (0.13656917844354743,0.2001593013620731)--(0.22603564326488124,0.11592720692047373);
\draw[red, thick] (0.679280231333423,0.07685494809943949)--(0.7039360788752047,0.13691025596481718);
\draw[red, thick] (0.679280231333423,0.07685494809943949)--(0.5401162699392226,0.10592861733374764);
\draw[red, thick] (0.7039360788752047,0.13691025596481718)--(0.692984094105217,0.16657177838006781);
\draw[red, thick] (0.692984094105217,0.16657177838006781)--(0.5705030851078778,0.14956730017396283);
\draw[red, thick] (0.5705030851078778,0.14956730017396283)--(0.5401162699392226,0.10592861733374764);
\draw[red, thick] (0.6704261078411188,0.2983594041110991)--(0.692984094105217,0.16657177838006781);
\draw[red, thick] (0.6704261078411188,0.2983594041110991)--(0.609626416494816,0.3156521591180685);
\draw[red, thick] (0.5705030851078778,0.14956730017396283)--(0.609626416494816,0.3156521591180685);
\draw[red, thick] (0.45995291048752523,0.028965958412618953)--(0.4765623254864447,0.08494160121449679);
\draw[red, thick] (0.4765623254864447,0.08494160121449679)--(0.5401162699392226,0.10592861733374764);
\draw[red, thick] (0.38150587398438834,0.343981855497234)--(0.4497840091409947,0.1953790520155802);
\draw[red, thick] (0.38150587398438834,0.343981855497234)--(0.29113000848283804,0.21354918994579553);
\draw[red, thick] (0.297180677272192,0.20845266990432376)--(0.2922904993448189,0.21081962561616713);
\draw[red, thick] (0.297180677272192,0.20845266990432376)--(0.4415822666940058,0.18039376697821286);
\draw[red, thick] (0.4415822666940058,0.18039376697821286)--(0.4497840091409947,0.1953790520155802);
\draw[red, thick] (0.29113000848283804,0.21354918994579553)--(0.2922904993448189,0.21081962561616713);
\draw[red, thick] (0.4959932207799768,0.47088604498338843)--(0.5773694923567108,0.33302816738559193);
\draw[red, thick] (0.39386035741641695,0.3879950389869909)--(0.38150587398438834,0.343981855497234);
\draw[red, thick] (0.5773694923567108,0.33302816738559193)--(0.4497840091409947,0.1953790520155802);
\draw[red, thick] (0.19799452036188725,0.2983935547149186)--(0.29113000848283804,0.21354918994579553);
\draw[red, thick] (0.33570477350580863,0.0944446314103765)--(0.297180677272192,0.20845266990432376);
\draw[red, thick] (0.22603564326488124,0.11592720692047373)--(0.2922904993448189,0.21081962561616713);
\draw[red, thick] (0.4765623254864447,0.08494160121449679)--(0.4415822666940058,0.18039376697821286);
\draw[red, thick] (0.609626416494816,0.3156521591180685)--(0.5773694923567108,0.33302816738559193);
\draw[red, thick] (0.6966078171943033,0.517084675521012)--(0.8394170635303543,0.5386328082928883);
\draw[red, thick] (0.6966078171943033,0.517084675521012)--(0.7629987978835259,0.4097388057906942);
\draw[red, thick] (0.8394170635303543,0.5386328082928883)--(0.8199448304662863,0.43374447248530434);
\draw[red, thick] (0.7629987978835259,0.4097388057906942)--(0.8199448304662863,0.43374447248530434);
\draw[red, thick] (0.5092587354725757,0.5327225311464738)--(0.5854609328882098,0.5563630707966651);
\draw[red, thick] (0.6704261078411188,0.2983594041110991)--(0.706255398851686,0.31142037155747987);
\draw[red, thick] (0.706255398851686,0.31142037155747987)--(0.7629987978835259,0.4097388057906942);
\draw[red, thick] (0.6966078171943033,0.517084675521012)--(0.6791059840045244,0.528146008831672);
\draw[red, thick] (0.5854609328882098,0.5563630707966651)--(0.6791059840045244,0.528146008831672);
\draw[red, thick] (0.8394170635303543,0.5386328082928883)--(0.8421316549547286,0.5424862865901321);
\draw[red, thick] (0.6791059840045244,0.528146008831672)--(0.7961709115289796,0.6183095407912882);
\draw[red, thick] (0.8421316549547286,0.5424862865901321)--(0.7961709115289796,0.6183095407912882);
\draw[red, thick] (0.8180743974620243,0.1360181565046783)--(0.7847271846057482,0.2833695475015336);
\draw[red, thick] (0.7039360788752047,0.13691025596481718)--(0.8180743974620243,0.1360181565046783);
\draw[red, thick] (0.706255398851686,0.31142037155747987)--(0.7847271846057482,0.2833695475015336);
\draw[red, thick] (0.8340383503796273,0.7758938836779715)--(0.8378339831936373,0.7542743252554434);
\draw[red, thick] (0.663654233906528,0.6697057286556454)--(0.6320097831139766,0.6120438657752406);
\draw[red, thick] (0.663654233906528,0.6697057286556454)--(0.8042869452731249,0.7046892265869095);
\draw[red, thick] (0.6320097831139766,0.6120438657752406)--(0.7954628450760505,0.6220772182985786);
\draw[red, thick] (0.8042869452731249,0.7046892265869095)--(0.7954628450760505,0.6220772182985786);
\draw[red, thick] (0.6347894819842471,0.8155039577188894)--(0.663654233906528,0.6697057286556454);
\draw[red, thick] (0.5854609328882098,0.5563630707966651)--(0.6320097831139766,0.6120438657752406);
\draw[red, thick] (0.8378339831936373,0.7542743252554434)--(0.8042869452731249,0.7046892265869095);
\draw[red, thick] (0.7961709115289796,0.6183095407912882)--(0.7954628450760505,0.6220772182985786);
\draw[black, very thick] (0.1813061501361128,0.42139697266953957)--(0.25616296056720894,0.325003496798576);
\draw[black, very thick] (0.1813061501361128,0.42139697266953957)--(0.25532816389877266,0.4461882000106892);
\draw[black, very thick] (0.1813061501361128,0.42139697266953957)--(0.2585381570265892,0.4541014013513852);
\draw[black, very thick] (0.14872563829292249,0.5541217325950919)--(0.18134071416677522,0.6074653564194161);
\draw[black, very thick] (0.14872563829292249,0.5541217325950919)--(0.1726754312359303,0.6109192601286578);
\draw[black, very thick] (0.05789114625497006,0.7591447714783306)--(0.025416319820836576,0.7229033871050434);
\draw[black, very thick] (0.05789114625497006,0.7591447714783306)--(0.07059719372583992,0.7820235001551941);
\draw[black, very thick] (0.025416319820836576,0.7229033871050434)--(0.07059719372583992,0.7820235001551941);
\draw[black, very thick] (0.07059719372583992,0.7820235001551941)--(0.06284152846679965,0.9189337987707326);
\draw[black, very thick] (0.07059719372583992,0.7820235001551941)--(0.18563030063821997,0.7444433737589387);
\draw[black, very thick] (0.816871669547315,0.8207650844535032)--(0.797244196123064,0.8483019422161482);
\draw[black, very thick] (0.816871669547315,0.8207650844535032)--(0.7024556644193323,0.853730955492912);
\draw[black, very thick] (0.816871669547315,0.8207650844535032)--(0.84050884537445,0.7523984527019172);
\draw[black, very thick] (0.797244196123064,0.8483019422161482)--(0.7024556644193323,0.853730955492912);
\draw[black, very thick] (0.5150163133378375,0.8659430761214151)--(0.47164836244225017,0.9109390813944592);
\draw[black, very thick] (0.5150163133378375,0.8659430761214151)--(0.5337666871494476,0.7563578295951302);
\draw[black, very thick] (0.47164836244225017,0.9109390813944592)--(0.38653207447931925,0.9992506770062238);
\draw[black, very thick] (0.38653207447931925,0.9992506770062238)--(0.2921345429365123,0.9608102969999692);
\draw[black, very thick] (0.38653207447931925,0.9992506770062238)--(0.30700725570527704,0.9512527117166875);
\draw[black, very thick] (0.29224699803261267,0.34151062925235404)--(0.25616296056720894,0.325003496798576);
\draw[black, very thick] (0.29224699803261267,0.34151062925235404)--(0.35675114120119855,0.3082552180370413);
\draw[black, very thick] (0.29224699803261267,0.34151062925235404)--(0.37932348476575223,0.34083217785400804);
\draw[black, very thick] (0.25616296056720894,0.325003496798576)--(0.35675114120119855,0.3082552180370413);
\draw[black, very thick] (0.3693375759616444,0.43463383746958717)--(0.37181151999371226,0.4343189103317002);
\draw[black, very thick] (0.3693375759616444,0.43463383746958717)--(0.4100600707523283,0.4294499616155625);
\draw[black, very thick] (0.3693375759616444,0.43463383746958717)--(0.4169904242473743,0.42937728712224627);
\draw[black, very thick] (0.37181151999371226,0.4343189103317002)--(0.4100600707523283,0.4294499616155625);
\draw[black, very thick] (0.37181151999371226,0.4343189103317002)--(0.4169904242473743,0.42937728712224627);
\draw[black, very thick] (0.37181151999371226,0.4343189103317002)--(0.3849643126762702,0.3563026374685184);
\draw[black, very thick] (0.4100600707523283,0.4294499616155625)--(0.4169904242473743,0.42937728712224627);
\draw[black, very thick] (0.4100600707523283,0.4294499616155625)--(0.3849643126762702,0.3563026374685184);
\draw[black, very thick] (0.25532816389877266,0.4461882000106892)--(0.2585381570265892,0.4541014013513852);
\draw[black, very thick] (0.25532816389877266,0.4461882000106892)--(0.26648526291374264,0.4998123096218199);
\draw[black, very thick] (0.25532816389877266,0.4461882000106892)--(0.2919655586376268,0.49690587581737233);
\draw[black, very thick] (0.15937905866236854,0.9047088031490487)--(0.1590067661594326,0.9455878746166037);
\draw[black, very thick] (0.15937905866236854,0.9047088031490487)--(0.17361789056520271,0.9419355930212361);
\draw[black, very thick] (0.15937905866236854,0.9047088031490487)--(0.06284152846679965,0.9189337987707326);
\draw[black, very thick] (0.1590067661594326,0.9455878746166037)--(0.17361789056520271,0.9419355930212361);
\draw[black, very thick] (0.1590067661594326,0.9455878746166037)--(0.06284152846679965,0.9189337987707326);
\draw[black, very thick] (0.17361789056520271,0.9419355930212361)--(0.06284152846679965,0.9189337987707326);
\draw[black, very thick] (0.17361789056520271,0.9419355930212361)--(0.2921345429365123,0.9608102969999692);
\draw[black, very thick] (0.2921345429365123,0.9608102969999692)--(0.30700725570527704,0.9512527117166875);
\draw[black, very thick] (0.2921345429365123,0.9608102969999692)--(0.304643617111923,0.9074897229510144);
\draw[black, very thick] (0.2526856156454257,0.7769398874936125)--(0.27892638396565794,0.7990815186159101);
\draw[black, very thick] (0.2526856156454257,0.7769398874936125)--(0.29981130217595947,0.8180189617063004);
\draw[black, very thick] (0.2526856156454257,0.7769398874936125)--(0.30055791219730454,0.8318425158098488);
\draw[black, very thick] (0.2526856156454257,0.7769398874936125)--(0.18563030063821997,0.7444433737589387);
\draw[black, very thick] (0.27892638396565794,0.7990815186159101)--(0.3010446550260852,0.8408546046323935);
\draw[black, very thick] (0.27892638396565794,0.7990815186159101)--(0.29981130217595947,0.8180189617063004);
\draw[black, very thick] (0.27892638396565794,0.7990815186159101)--(0.30055791219730454,0.8318425158098488);
\draw[black, very thick] (0.19387030084463816,0.7007344954508157)--(0.18563030063821997,0.7444433737589387);
\draw[black, very thick] (0.37293603650876134,0.6988993954844914)--(0.3982187058250188,0.6582318996866904);
\draw[black, very thick] (0.3982187058250188,0.6582318996866904)--(0.34509732830012835,0.6165202358288133);
\draw[black, very thick] (0.3982187058250188,0.6582318996866904)--(0.36683993239497326,0.6217975775006084);
\draw[black, very thick] (0.18134071416677522,0.6074653564194161)--(0.1726754312359303,0.6109192601286578);
\draw[black, very thick] (0.18134071416677522,0.6074653564194161)--(0.232149332632162,0.5956374112256204);
\draw[black, very thick] (0.18134071416677522,0.6074653564194161)--(0.2277040583383429,0.5644397202114);
\draw[black, very thick] (0.1726754312359303,0.6109192601286578)--(0.232149332632162,0.5956374112256204);
\draw[black, very thick] (0.1726754312359303,0.6109192601286578)--(0.2277040583383429,0.5644397202114);
\draw[black, very thick] (0.34509732830012835,0.6165202358288133)--(0.36683993239497326,0.6217975775006084);
\draw[black, very thick] (0.34509732830012835,0.6165202358288133)--(0.3498485323164778,0.5891158637590141);
\draw[black, very thick] (0.34509732830012835,0.6165202358288133)--(0.34045815137782165,0.574156595807007);
\draw[black, very thick] (0.232149332632162,0.5956374112256204)--(0.23156648132676858,0.5534121318594152);
\draw[black, very thick] (0.232149332632162,0.5956374112256204)--(0.2277040583383429,0.5644397202114);
\draw[black, very thick] (0.36683993239497326,0.6217975775006084)--(0.3498485323164778,0.5891158637590141);
\draw[black, very thick] (0.36683993239497326,0.6217975775006084)--(0.34045815137782165,0.574156595807007);
\draw[black, very thick] (0.3010446550260852,0.8408546046323935)--(0.304643617111923,0.9074897229510144);
\draw[black, very thick] (0.3010446550260852,0.8408546046323935)--(0.29981130217595947,0.8180189617063004);
\draw[black, very thick] (0.3010446550260852,0.8408546046323935)--(0.30055791219730454,0.8318425158098488);
\draw[black, very thick] (0.30700725570527704,0.9512527117166875)--(0.304643617111923,0.9074897229510144);
\draw[black, very thick] (0.29981130217595947,0.8180189617063004)--(0.30055791219730454,0.8318425158098488);
\draw[black, very thick] (0.5337666871494476,0.7563578295951302)--(0.4828113417610102,0.6605018615470278);
\draw[black, very thick] (0.2585381570265892,0.4541014013513852)--(0.26648526291374264,0.4998123096218199);
\draw[black, very thick] (0.2585381570265892,0.4541014013513852)--(0.2919655586376268,0.49690587581737233);
\draw[black, very thick] (0.23156648132676858,0.5534121318594152)--(0.2277040583383429,0.5644397202114);
\draw[black, very thick] (0.23156648132676858,0.5534121318594152)--(0.293302842223799,0.5346545129197736);
\draw[black, very thick] (0.26648526291374264,0.4998123096218199)--(0.293302842223799,0.5346545129197736);
\draw[black, very thick] (0.26648526291374264,0.4998123096218199)--(0.2919655586376268,0.49690587581737233);
\draw[black, very thick] (0.293302842223799,0.5346545129197736)--(0.2919655586376268,0.49690587581737233);
\draw[black, very thick] (0.293302842223799,0.5346545129197736)--(0.34045815137782165,0.574156595807007);
\draw[black, very thick] (0.3498485323164778,0.5891158637590141)--(0.34045815137782165,0.574156595807007);
\draw[black, very thick] (0.5084688879484911,0.5290406992216823)--(0.49024128124118205,0.4678639258176544);
\draw[black, very thick] (0.5084688879484911,0.5290406992216823)--(0.5800834392674282,0.554694787525681);
\draw[black, very thick] (0.49024128124118205,0.4678639258176544)--(0.5326754741324017,0.40874338970209434);
\draw[black, very thick] (0.38355244884048856,0.06922894316707184)--(0.44145691925859737,0.03871333161972445);
\draw[black, very thick] (0.38355244884048856,0.06922894316707184)--(0.37508174221498614,0.07369299906895496);
\draw[black, very thick] (0.38355244884048856,0.06922894316707184)--(0.35676393151346963,0.08334647127462609);
\draw[black, very thick] (0.38355244884048856,0.06922894316707184)--(0.33402090497369247,0.09314723041436684);
\draw[black, very thick] (0.44145691925859737,0.03871333161972445)--(0.37508174221498614,0.07369299906895496);
\draw[black, very thick] (0.44145691925859737,0.03871333161972445)--(0.4744519247018341,0.0778293073440714);
\draw[black, very thick] (0.37508174221498614,0.07369299906895496)--(0.35676393151346963,0.08334647127462609);
\draw[black, very thick] (0.37508174221498614,0.07369299906895496)--(0.33402090497369247,0.09314723041436684);
\draw[black, very thick] (0.35676393151346963,0.08334647127462609)--(0.33402090497369247,0.09314723041436684);
\draw[black, very thick] (0.27052606386950184,0.04422519848677735)--(0.33402090497369247,0.09314723041436684);
\draw[black, very thick] (0.27052606386950184,0.04422519848677735)--(0.25415869288763704,0.031614330708849434);
\draw[black, very thick] (0.19012226347068706,-0.017724988674975257)--(0.25415869288763704,0.031614330708849434);
\draw[black, very thick] (0.19012226347068706,-0.017724988674975257)--(0.1872786858015116,0.01718136555514127);
\draw[black, very thick] (0.19012226347068706,-0.017724988674975257)--(0.18417778144474048,0.00928081221598432);
\draw[black, very thick] (0.19012226347068706,-0.017724988674975257)--(0.17858300600531926,-0.004973681536212152);
\draw[black, very thick] (0.25415869288763704,0.031614330708849434)--(0.1872786858015116,0.01718136555514127);
\draw[black, very thick] (0.25415869288763704,0.031614330708849434)--(0.18417778144474048,0.00928081221598432);
\draw[black, very thick] (0.19976901885816142,0.14065706287384386)--(0.15362496192050873,0.1841013922061019);
\draw[black, very thick] (0.19976901885816142,0.14065706287384386)--(0.153557661738525,0.18416475488949516);
\draw[black, very thick] (0.19976901885816142,0.14065706287384386)--(0.14846174817708205,0.18896252434680957);
\draw[black, very thick] (0.15362496192050873,0.1841013922061019)--(0.153557661738525,0.18416475488949516);
\draw[black, very thick] (0.15362496192050873,0.1841013922061019)--(0.14846174817708205,0.18896252434680957);
\draw[black, very thick] (0.153557661738525,0.18416475488949516)--(0.14846174817708205,0.18896252434680957);
\draw[black, very thick] (0.1872786858015116,0.01718136555514127)--(0.18417778144474048,0.00928081221598432);
\draw[black, very thick] (0.1872786858015116,0.01718136555514127)--(0.17858300600531926,-0.004973681536212152);
\draw[black, very thick] (0.18417778144474048,0.00928081221598432)--(0.17858300600531926,-0.004973681536212152);
\draw[black, very thick] (0.694675762855622,0.11435450611503778)--(0.6796560855289289,0.2444362909933815);
\draw[black, very thick] (0.694675762855622,0.11435450611503778)--(0.8118874384593696,0.1633564836689559);
\draw[black, very thick] (0.6796560855289289,0.2444362909933815)--(0.8118874384593696,0.1633564836689559);
\draw[black, very thick] (0.4744519247018341,0.0778293073440714)--(0.5263999923459467,0.10139917795079811);
\draw[black, very thick] (0.4744519247018341,0.0778293073440714)--(0.45877901890888895,0.13346796697448662);
\draw[black, very thick] (0.44625527759791844,0.20305910066999244)--(0.43025304052789864,0.18259516696622943);
\draw[black, very thick] (0.44625527759791844,0.20305910066999244)--(0.48139088714400735,0.22947900357015216);
\draw[black, very thick] (0.44625527759791844,0.20305910066999244)--(0.45162040934298203,0.19736030299255442);
\draw[black, very thick] (0.44625527759791844,0.20305910066999244)--(0.44500451138548675,0.17105528445839943);
\draw[black, very thick] (0.35675114120119855,0.3082552180370413)--(0.37932348476575223,0.34083217785400804);
\draw[black, very thick] (0.35675114120119855,0.3082552180370413)--(0.3849643126762702,0.3563026374685184);
\draw[black, very thick] (0.37932348476575223,0.34083217785400804)--(0.3849643126762702,0.3563026374685184);
\draw[black, very thick] (0.3781940958333602,0.1927108234591293)--(0.43025304052789864,0.18259516696622943);
\draw[black, very thick] (0.43025304052789864,0.18259516696622943)--(0.45162040934298203,0.19736030299255442);
\draw[black, very thick] (0.43025304052789864,0.18259516696622943)--(0.44500451138548675,0.17105528445839943);
\draw[black, very thick] (0.5646392290876473,0.35459424606590373)--(0.5570056583155196,0.3675261222104846);
\draw[black, very thick] (0.5646392290876473,0.35459424606590373)--(0.5326754741324017,0.40874338970209434);
\draw[black, very thick] (0.5646392290876473,0.35459424606590373)--(0.5464741713477421,0.2996958987076048);
\draw[black, very thick] (0.5570056583155196,0.3675261222104846)--(0.5326754741324017,0.40874338970209434);
\draw[black, very thick] (0.5570056583155196,0.3675261222104846)--(0.5464741713477421,0.2996958987076048);
\draw[black, very thick] (0.48139088714400735,0.22947900357015216)--(0.4952009553583087,0.24437837632781984);
\draw[black, very thick] (0.48139088714400735,0.22947900357015216)--(0.45162040934298203,0.19736030299255442);
\draw[black, very thick] (0.4952009553583087,0.24437837632781984)--(0.5464741713477421,0.2996958987076048);
\draw[black, very thick] (0.45162040934298203,0.19736030299255442)--(0.44500451138548675,0.17105528445839943);
\draw[black, very thick] (0.44500451138548675,0.17105528445839943)--(0.45877901890888895,0.13346796697448662);
\draw[black, very thick] (0.7221119059162853,0.4758477701896109)--(0.8272393746070328,0.4730369662060027);
\draw[black, very thick] (0.7221119059162853,0.4758477701896109)--(0.7711711320631728,0.41318386283041797);
\draw[black, very thick] (0.8272393746070328,0.4730369662060027)--(0.7711711320631728,0.41318386283041797);
\draw[black, very thick] (0.8272393746070328,0.4730369662060027)--(0.8851632825205141,0.5493531197369946);
\draw[black, very thick] (0.5800834392674282,0.554694787525681)--(0.6188652259059673,0.5462977123528099);
\draw[black, very thick] (0.5800834392674282,0.554694787525681)--(0.6361928686502474,0.5410765590383083);
\draw[black, very thick] (0.5800834392674282,0.554694787525681)--(0.6273026388087081,0.6064132747783626);
\draw[black, very thick] (0.6188652259059673,0.5462977123528099)--(0.6361928686502474,0.5410765590383083);
\draw[black, very thick] (0.6188652259059673,0.5462977123528099)--(0.6273026388087081,0.6064132747783626);
\draw[black, very thick] (0.6361928686502474,0.5410765590383083)--(0.6273026388087081,0.6064132747783626);
\draw[black, very thick] (0.7880511358690263,0.6120556809230696)--(0.7843733421016718,0.6092230402432683);
\draw[black, very thick] (0.7880511358690263,0.6120556809230696)--(0.7940377748297177,0.6166665966317921);
\draw[black, very thick] (0.7880511358690263,0.6120556809230696)--(0.7851149531728918,0.6097942298808102);
\draw[black, very thick] (0.7880511358690263,0.6120556809230696)--(0.7860546895085898,0.6214997109685254);
\draw[black, very thick] (0.7843733421016718,0.6092230402432683)--(0.7940377748297177,0.6166665966317921);
\draw[black, very thick] (0.7843733421016718,0.6092230402432683)--(0.7851149531728918,0.6097942298808102);
\draw[black, very thick] (0.7843733421016718,0.6092230402432683)--(0.7860546895085898,0.6214997109685254);
\draw[black, very thick] (0.7940377748297177,0.6166665966317921)--(0.7851149531728918,0.6097942298808102);
\draw[black, very thick] (0.7940377748297177,0.6166665966317921)--(0.7860546895085898,0.6214997109685254);
\draw[black, very thick] (0.7851149531728918,0.6097942298808102)--(0.7860546895085898,0.6214997109685254);
\draw[black, very thick] (0.8851632825205141,0.5493531197369946)--(0.9421973439307684,0.5584544111718442);
\draw[black, very thick] (1.0241664694608577,0.57153475091379)--(1.0378006206639054,0.5737104399756852);
\draw[black, very thick] (1.0241664694608577,0.57153475091379)--(1.0418135837704467,0.6112240638453353);
\draw[black, very thick] (1.0378006206639054,0.5737104399756852)--(1.0605948498489768,0.5980528201790608);
\draw[black, very thick] (1.0378006206639054,0.5737104399756852)--(1.0418135837704467,0.6112240638453353);
\draw[black, very thick] (1.0378006206639054,0.5737104399756852)--(1.0588598410834567,0.5992695764466222);
\draw[black, very thick] (1.0823433206558173,0.580818404653238)--(1.0605948498489768,0.5980528201790608);
\draw[black, very thick] (1.0823433206558173,0.580818404653238)--(1.0588598410834567,0.5992695764466222);
\draw[black, very thick] (0.9421973439307684,0.5584544111718442)--(0.9625363197915181,0.6668209632607429);
\draw[black, very thick] (1.0605948498489768,0.5980528201790608)--(1.0418135837704467,0.6112240638453353);
\draw[black, very thick] (1.0605948498489768,0.5980528201790608)--(1.0588598410834567,0.5992695764466222);
\draw[black, very thick] (0.9625363197915181,0.6668209632607429)--(0.9168260577530807,0.6988774282956411);
\draw[black, very thick] (0.9625363197915181,0.6668209632607429)--(0.9338936187195976,0.6869080007107792);
\draw[black, very thick] (1.0418135837704467,0.6112240638453353)--(1.0588598410834567,0.5992695764466222);
\draw[black, very thick] (0.9168260577530807,0.6988774282956411)--(0.9338936187195976,0.6869080007107792);
\draw[black, very thick] (0.9168260577530807,0.6988774282956411)--(0.84050884537445,0.7523984527019172);
\draw[black, very thick] (0.84050884537445,0.7523984527019172)--(0.8032960357300718,0.6954122436021932);
\draw[black, very thick] (0.642114287806097,0.6304560849618224)--(0.6595268144771232,0.6905536626697233);
\draw[black, very thick] (0.642114287806097,0.6304560849618224)--(0.6273026388087081,0.6064132747783626);
\draw[black, very thick] (0.7630519194473339,0.6944316881617858)--(0.7580659591276837,0.6931913911438812);
\draw[black, very thick] (0.7630519194473339,0.6944316881617858)--(0.7097088513509334,0.6811621785700015);
\draw[black, very thick] (0.7630519194473339,0.6944316881617858)--(0.8032960357300718,0.6954122436021932);
\draw[black, very thick] (0.7580659591276837,0.6931913911438812)--(0.7097088513509334,0.6811621785700015);
\draw[black, very thick] (0.7580659591276837,0.6931913911438812)--(0.6984784752351365,0.6783685337927906);
\draw[black, very thick] (0.7580659591276837,0.6931913911438812)--(0.8032960357300718,0.6954122436021932);
\draw[black, very thick] (0.7097088513509334,0.6811621785700015)--(0.6984784752351365,0.6783685337927906);
\draw[black, very thick] (0.7097088513509334,0.6811621785700015)--(0.6595268144771232,0.6905536626697233);
\draw[black, very thick] (0.6984784752351365,0.6783685337927906)--(0.714648332376192,0.6171165250868014);
\draw[black, very thick] (0.6984784752351365,0.6783685337927906)--(0.6595268144771232,0.6905536626697233);
\fill[blue] (0.1813061501361128,0.42139697266953957) circle (0.1pt);\draw(0.1813061501361128,0.42139697266953957) circle (0.1pt);
\fill[blue] (0.14872563829292249,0.5541217325950919) circle (0.1pt);\draw(0.14872563829292249,0.5541217325950919) circle (0.1pt);
\fill[blue] (0.05789114625497006,0.7591447714783306) circle (0.1pt);\draw(0.05789114625497006,0.7591447714783306) circle (0.1pt);
\fill[blue] (0.025416319820836576,0.7229033871050434) circle (0.1pt);\draw(0.025416319820836576,0.7229033871050434) circle (0.1pt);
\fill[blue] (0.07059719372583992,0.7820235001551941) circle (0.1pt);\draw(0.07059719372583992,0.7820235001551941) circle (0.1pt);
\fill[blue] (0.816871669547315,0.8207650844535032) circle (0.1pt);\draw(0.816871669547315,0.8207650844535032) circle (0.1pt);
\fill[blue] (0.797244196123064,0.8483019422161482) circle (0.1pt);\draw(0.797244196123064,0.8483019422161482) circle (0.1pt);
\fill[blue] (0.5150163133378375,0.8659430761214151) circle (0.1pt);\draw(0.5150163133378375,0.8659430761214151) circle (0.1pt);
\fill[blue] (0.47164836244225017,0.9109390813944592) circle (0.1pt);\draw(0.47164836244225017,0.9109390813944592) circle (0.1pt);
\fill[blue] (0.38653207447931925,0.9992506770062238) circle (0.1pt);\draw(0.38653207447931925,0.9992506770062238) circle (0.1pt);
\fill[blue] (0.7024556644193323,0.853730955492912) circle (0.1pt);\draw(0.7024556644193323,0.853730955492912) circle (0.1pt);
\fill[blue] (0.29224699803261267,0.34151062925235404) circle (0.1pt);\draw(0.29224699803261267,0.34151062925235404) circle (0.1pt);
\fill[blue] (0.25616296056720894,0.325003496798576) circle (0.1pt);\draw(0.25616296056720894,0.325003496798576) circle (0.1pt);
\fill[blue] (0.3693375759616444,0.43463383746958717) circle (0.1pt);\draw(0.3693375759616444,0.43463383746958717) circle (0.1pt);
\fill[blue] (0.37181151999371226,0.4343189103317002) circle (0.1pt);\draw(0.37181151999371226,0.4343189103317002) circle (0.1pt);
\fill[blue] (0.4100600707523283,0.4294499616155625) circle (0.1pt);\draw(0.4100600707523283,0.4294499616155625) circle (0.1pt);
\fill[blue] (0.25532816389877266,0.4461882000106892) circle (0.1pt);\draw(0.25532816389877266,0.4461882000106892) circle (0.1pt);
\fill[blue] (0.15937905866236854,0.9047088031490487) circle (0.1pt);\draw(0.15937905866236854,0.9047088031490487) circle (0.1pt);
\fill[blue] (0.1590067661594326,0.9455878746166037) circle (0.1pt);\draw(0.1590067661594326,0.9455878746166037) circle (0.1pt);
\fill[blue] (0.17361789056520271,0.9419355930212361) circle (0.1pt);\draw(0.17361789056520271,0.9419355930212361) circle (0.1pt);
\fill[blue] (0.06284152846679965,0.9189337987707326) circle (0.1pt);\draw(0.06284152846679965,0.9189337987707326) circle (0.1pt);
\fill[blue] (0.2921345429365123,0.9608102969999692) circle (0.1pt);\draw(0.2921345429365123,0.9608102969999692) circle (0.1pt);
\fill[blue] (0.2526856156454257,0.7769398874936125) circle (0.1pt);\draw(0.2526856156454257,0.7769398874936125) circle (0.1pt);
\fill[blue] (0.27892638396565794,0.7990815186159101) circle (0.1pt);\draw(0.27892638396565794,0.7990815186159101) circle (0.1pt);
\fill[blue] (0.19387030084463816,0.7007344954508157) circle (0.1pt);\draw(0.19387030084463816,0.7007344954508157) circle (0.1pt);
\fill[blue] (0.37293603650876134,0.6988993954844914) circle (0.1pt);\draw(0.37293603650876134,0.6988993954844914) circle (0.1pt);
\fill[blue] (0.3982187058250188,0.6582318996866904) circle (0.1pt);\draw(0.3982187058250188,0.6582318996866904) circle (0.1pt);
\fill[blue] (0.18134071416677522,0.6074653564194161) circle (0.1pt);\draw(0.18134071416677522,0.6074653564194161) circle (0.1pt);
\fill[blue] (0.1726754312359303,0.6109192601286578) circle (0.1pt);\draw(0.1726754312359303,0.6109192601286578) circle (0.1pt);
\fill[blue] (0.34509732830012835,0.6165202358288133) circle (0.1pt);\draw(0.34509732830012835,0.6165202358288133) circle (0.1pt);
\fill[blue] (0.232149332632162,0.5956374112256204) circle (0.1pt);\draw(0.232149332632162,0.5956374112256204) circle (0.1pt);
\fill[blue] (0.36683993239497326,0.6217975775006084) circle (0.1pt);\draw(0.36683993239497326,0.6217975775006084) circle (0.1pt);
\fill[blue] (0.3010446550260852,0.8408546046323935) circle (0.1pt);\draw(0.3010446550260852,0.8408546046323935) circle (0.1pt);
\fill[blue] (0.30700725570527704,0.9512527117166875) circle (0.1pt);\draw(0.30700725570527704,0.9512527117166875) circle (0.1pt);
\fill[blue] (0.304643617111923,0.9074897229510144) circle (0.1pt);\draw(0.304643617111923,0.9074897229510144) circle (0.1pt);
\fill[blue] (0.29981130217595947,0.8180189617063004) circle (0.1pt);\draw(0.29981130217595947,0.8180189617063004) circle (0.1pt);
\fill[blue] (0.30055791219730454,0.8318425158098488) circle (0.1pt);\draw(0.30055791219730454,0.8318425158098488) circle (0.1pt);
\fill[blue] (0.18563030063821997,0.7444433737589387) circle (0.1pt);\draw(0.18563030063821997,0.7444433737589387) circle (0.1pt);
\fill[blue] (0.5337666871494476,0.7563578295951302) circle (0.1pt);\draw(0.5337666871494476,0.7563578295951302) circle (0.1pt);
\fill[blue] (0.4828113417610102,0.6605018615470278) circle (0.1pt);\draw(0.4828113417610102,0.6605018615470278) circle (0.1pt);
\fill[blue] (0.2585381570265892,0.4541014013513852) circle (0.1pt);\draw(0.2585381570265892,0.4541014013513852) circle (0.1pt);
\fill[blue] (0.23156648132676858,0.5534121318594152) circle (0.1pt);\draw(0.23156648132676858,0.5534121318594152) circle (0.1pt);
\fill[blue] (0.2277040583383429,0.5644397202114) circle (0.1pt);\draw(0.2277040583383429,0.5644397202114) circle (0.1pt);
\fill[blue] (0.26648526291374264,0.4998123096218199) circle (0.1pt);\draw(0.26648526291374264,0.4998123096218199) circle (0.1pt);
\fill[blue] (0.293302842223799,0.5346545129197736) circle (0.1pt);\draw(0.293302842223799,0.5346545129197736) circle (0.1pt);
\fill[blue] (0.3498485323164778,0.5891158637590141) circle (0.1pt);\draw(0.3498485323164778,0.5891158637590141) circle (0.1pt);
\fill[blue] (0.2919655586376268,0.49690587581737233) circle (0.1pt);\draw(0.2919655586376268,0.49690587581737233) circle (0.1pt);
\fill[blue] (0.34045815137782165,0.574156595807007) circle (0.1pt);\draw(0.34045815137782165,0.574156595807007) circle (0.1pt);
\fill[blue] (0.5084688879484911,0.5290406992216823) circle (0.1pt);\draw(0.5084688879484911,0.5290406992216823) circle (0.1pt);
\fill[blue] (0.4169904242473743,0.42937728712224627) circle (0.1pt);\draw(0.4169904242473743,0.42937728712224627) circle (0.1pt);
\fill[blue] (0.49024128124118205,0.4678639258176544) circle (0.1pt);\draw(0.49024128124118205,0.4678639258176544) circle (0.1pt);
\fill[blue] (0.38355244884048856,0.06922894316707184) circle (0.1pt);\draw(0.38355244884048856,0.06922894316707184) circle (0.1pt);
\fill[blue] (0.44145691925859737,0.03871333161972445) circle (0.1pt);\draw(0.44145691925859737,0.03871333161972445) circle (0.1pt);
\fill[blue] (0.37508174221498614,0.07369299906895496) circle (0.1pt);\draw(0.37508174221498614,0.07369299906895496) circle (0.1pt);
\fill[blue] (0.35676393151346963,0.08334647127462609) circle (0.1pt);\draw(0.35676393151346963,0.08334647127462609) circle (0.1pt);
\fill[blue] (0.27052606386950184,0.04422519848677735) circle (0.1pt);\draw(0.27052606386950184,0.04422519848677735) circle (0.1pt);
\fill[blue] (0.33402090497369247,0.09314723041436684) circle (0.1pt);\draw(0.33402090497369247,0.09314723041436684) circle (0.1pt);
\fill[blue] (0.19012226347068706,-0.017724988674975257) circle (0.1pt);\draw(0.19012226347068706,-0.017724988674975257) circle (0.1pt);
\fill[blue] (0.25415869288763704,0.031614330708849434) circle (0.1pt);\draw(0.25415869288763704,0.031614330708849434) circle (0.1pt);
\fill[blue] (0.19976901885816142,0.14065706287384386) circle (0.1pt);\draw(0.19976901885816142,0.14065706287384386) circle (0.1pt);
\fill[blue] (0.15362496192050873,0.1841013922061019) circle (0.1pt);\draw(0.15362496192050873,0.1841013922061019) circle (0.1pt);
\fill[blue] (0.153557661738525,0.18416475488949516) circle (0.1pt);\draw(0.153557661738525,0.18416475488949516) circle (0.1pt);
\fill[blue] (0.14846174817708205,0.18896252434680957) circle (0.1pt);\draw(0.14846174817708205,0.18896252434680957) circle (0.1pt);
\fill[blue] (0.1872786858015116,0.01718136555514127) circle (0.1pt);\draw(0.1872786858015116,0.01718136555514127) circle (0.1pt);
\fill[blue] (0.18417778144474048,0.00928081221598432) circle (0.1pt);\draw(0.18417778144474048,0.00928081221598432) circle (0.1pt);
\fill[blue] (0.17858300600531926,-0.004973681536212152) circle (0.1pt);\draw(0.17858300600531926,-0.004973681536212152) circle (0.1pt);
\fill[blue] (0.694675762855622,0.11435450611503778) circle (0.1pt);\draw(0.694675762855622,0.11435450611503778) circle (0.1pt);
\fill[blue] (0.6796560855289289,0.2444362909933815) circle (0.1pt);\draw(0.6796560855289289,0.2444362909933815) circle (0.1pt);
\fill[blue] (0.4744519247018341,0.0778293073440714) circle (0.1pt);\draw(0.4744519247018341,0.0778293073440714) circle (0.1pt);
\fill[blue] (0.5263999923459467,0.10139917795079811) circle (0.1pt);\draw(0.5263999923459467,0.10139917795079811) circle (0.1pt);
\fill[blue] (0.44625527759791844,0.20305910066999244) circle (0.1pt);\draw(0.44625527759791844,0.20305910066999244) circle (0.1pt);
\fill[blue] (0.35675114120119855,0.3082552180370413) circle (0.1pt);\draw(0.35675114120119855,0.3082552180370413) circle (0.1pt);
\fill[blue] (0.37932348476575223,0.34083217785400804) circle (0.1pt);\draw(0.37932348476575223,0.34083217785400804) circle (0.1pt);
\fill[blue] (0.3781940958333602,0.1927108234591293) circle (0.1pt);\draw(0.3781940958333602,0.1927108234591293) circle (0.1pt);
\fill[blue] (0.43025304052789864,0.18259516696622943) circle (0.1pt);\draw(0.43025304052789864,0.18259516696622943) circle (0.1pt);
\fill[blue] (0.5646392290876473,0.35459424606590373) circle (0.1pt);\draw(0.5646392290876473,0.35459424606590373) circle (0.1pt);
\fill[blue] (0.5570056583155196,0.3675261222104846) circle (0.1pt);\draw(0.5570056583155196,0.3675261222104846) circle (0.1pt);
\fill[blue] (0.5326754741324017,0.40874338970209434) circle (0.1pt);\draw(0.5326754741324017,0.40874338970209434) circle (0.1pt);
\fill[blue] (0.3849643126762702,0.3563026374685184) circle (0.1pt);\draw(0.3849643126762702,0.3563026374685184) circle (0.1pt);
\fill[blue] (0.48139088714400735,0.22947900357015216) circle (0.1pt);\draw(0.48139088714400735,0.22947900357015216) circle (0.1pt);
\fill[blue] (0.4952009553583087,0.24437837632781984) circle (0.1pt);\draw(0.4952009553583087,0.24437837632781984) circle (0.1pt);
\fill[blue] (0.5464741713477421,0.2996958987076048) circle (0.1pt);\draw(0.5464741713477421,0.2996958987076048) circle (0.1pt);
\fill[blue] (0.45162040934298203,0.19736030299255442) circle (0.1pt);\draw(0.45162040934298203,0.19736030299255442) circle (0.1pt);
\fill[blue] (0.44500451138548675,0.17105528445839943) circle (0.1pt);\draw(0.44500451138548675,0.17105528445839943) circle (0.1pt);
\fill[blue] (0.45877901890888895,0.13346796697448662) circle (0.1pt);\draw(0.45877901890888895,0.13346796697448662) circle (0.1pt);
\fill[blue] (0.7221119059162853,0.4758477701896109) circle (0.1pt);\draw(0.7221119059162853,0.4758477701896109) circle (0.1pt);
\fill[blue] (0.8272393746070328,0.4730369662060027) circle (0.1pt);\draw(0.8272393746070328,0.4730369662060027) circle (0.1pt);
\fill[blue] (0.7711711320631728,0.41318386283041797) circle (0.1pt);\draw(0.7711711320631728,0.41318386283041797) circle (0.1pt);
\fill[blue] (0.5800834392674282,0.554694787525681) circle (0.1pt);\draw(0.5800834392674282,0.554694787525681) circle (0.1pt);
\fill[blue] (0.6188652259059673,0.5462977123528099) circle (0.1pt);\draw(0.6188652259059673,0.5462977123528099) circle (0.1pt);
\fill[blue] (0.6361928686502474,0.5410765590383083) circle (0.1pt);\draw(0.6361928686502474,0.5410765590383083) circle (0.1pt);
\fill[blue] (0.7880511358690263,0.6120556809230696) circle (0.1pt);\draw(0.7880511358690263,0.6120556809230696) circle (0.1pt);
\fill[blue] (0.7843733421016718,0.6092230402432683) circle (0.1pt);\draw(0.7843733421016718,0.6092230402432683) circle (0.1pt);
\fill[blue] (0.7940377748297177,0.6166665966317921) circle (0.1pt);\draw(0.7940377748297177,0.6166665966317921) circle (0.1pt);
\fill[blue] (0.7851149531728918,0.6097942298808102) circle (0.1pt);\draw(0.7851149531728918,0.6097942298808102) circle (0.1pt);
\fill[blue] (0.8108312686319776,0.5941237704915503) circle (0.1pt);\draw(0.8108312686319776,0.5941237704915503) circle (0.1pt);
\fill[blue] (0.8118874384593696,0.1633564836689559) circle (0.1pt);\draw(0.8118874384593696,0.1633564836689559) circle (0.1pt);
\fill[blue] (0.8851632825205141,0.5493531197369946) circle (0.1pt);\draw(0.8851632825205141,0.5493531197369946) circle (0.1pt);
\fill[blue] (1.0241664694608577,0.57153475091379) circle (0.1pt);\draw(1.0241664694608577,0.57153475091379) circle (0.1pt);
\fill[blue] (1.0378006206639054,0.5737104399756852) circle (0.1pt);\draw(1.0378006206639054,0.5737104399756852) circle (0.1pt);
\fill[blue] (1.0823433206558173,0.580818404653238) circle (0.1pt);\draw(1.0823433206558173,0.580818404653238) circle (0.1pt);
\fill[blue] (0.9421973439307684,0.5584544111718442) circle (0.1pt);\draw(0.9421973439307684,0.5584544111718442) circle (0.1pt);
\fill[blue] (1.0605948498489768,0.5980528201790608) circle (0.1pt);\draw(1.0605948498489768,0.5980528201790608) circle (0.1pt);
\fill[blue] (0.9625363197915181,0.6668209632607429) circle (0.1pt);\draw(0.9625363197915181,0.6668209632607429) circle (0.1pt);
\fill[blue] (1.0418135837704467,0.6112240638453353) circle (0.1pt);\draw(1.0418135837704467,0.6112240638453353) circle (0.1pt);
\fill[blue] (0.9168260577530807,0.6988774282956411) circle (0.1pt);\draw(0.9168260577530807,0.6988774282956411) circle (0.1pt);
\fill[blue] (1.0588598410834567,0.5992695764466222) circle (0.1pt);\draw(1.0588598410834567,0.5992695764466222) circle (0.1pt);
\fill[blue] (0.9338936187195976,0.6869080007107792) circle (0.1pt);\draw(0.9338936187195976,0.6869080007107792) circle (0.1pt);
\fill[blue] (0.84050884537445,0.7523984527019172) circle (0.1pt);\draw(0.84050884537445,0.7523984527019172) circle (0.1pt);
\fill[blue] (0.642114287806097,0.6304560849618224) circle (0.1pt);\draw(0.642114287806097,0.6304560849618224) circle (0.1pt);
\fill[blue] (0.7630519194473339,0.6944316881617858) circle (0.1pt);\draw(0.7630519194473339,0.6944316881617858) circle (0.1pt);
\fill[blue] (0.7580659591276837,0.6931913911438812) circle (0.1pt);\draw(0.7580659591276837,0.6931913911438812) circle (0.1pt);
\fill[blue] (0.7097088513509334,0.6811621785700015) circle (0.1pt);\draw(0.7097088513509334,0.6811621785700015) circle (0.1pt);
\fill[blue] (0.6984784752351365,0.6783685337927906) circle (0.1pt);\draw(0.6984784752351365,0.6783685337927906) circle (0.1pt);
\fill[blue] (0.714648332376192,0.6171165250868014) circle (0.1pt);\draw(0.714648332376192,0.6171165250868014) circle (0.1pt);
\fill[blue] (0.7860546895085898,0.6214997109685254) circle (0.1pt);\draw(0.7860546895085898,0.6214997109685254) circle (0.1pt);
\fill[blue] (0.8032960357300718,0.6954122436021932) circle (0.1pt);\draw(0.8032960357300718,0.6954122436021932) circle (0.1pt);
\fill[blue] (0.6595268144771232,0.6905536626697233) circle (0.1pt);\draw(0.6595268144771232,0.6905536626697233) circle (0.1pt);
\fill[blue] (0.6273026388087081,0.6064132747783626) circle (0.1pt);\draw(0.6273026388087081,0.6064132747783626) circle (0.1pt);
 \end{scope} 
 \draw (0.2,0.2) rectangle (0.8,0.8); 
 \end{tikzpicture}

%% file: VoronoiB5.tex
\begin{tikzpicture}[scale=8] 
 \begin{scope} 
\clip(0.2,0.2) rectangle (0.8,0.8);
\draw[red, thick] (0.8340383503796273,0.7758938836779715)--(0.822841928213028,0.8123889598178787);
\draw[red, thick] (0.16402452764015735,0.3738974712746203)--(0.17418575103719508,0.3127497753447539);
\draw[red, thick] (0.148397993788087,0.28793013543623097)--(0.17418575103719508,0.3127497753447539);
\draw[red, thick] (0.16402452764015735,0.3738974712746203)--(0.18588103912022907,0.40276003273976374);
\draw[red, thick] (0.18588103912022907,0.40276003273976374)--(0.13638721173493218,0.6043853598410402);
\draw[red, thick] (0.0072051335010195755,0.8663122452666777)--(0.07225031953352958,0.7988984371260657);
\draw[red, thick] (0.06964041109724696,0.7722567622110954)--(0.07225031953352958,0.7988984371260657);
\draw[red, thick] (0.822841928213028,0.8123889598178787)--(0.767274141410929,0.8903491846659611);
\draw[red, thick] (0.767274141410929,0.8903491846659611)--(0.6347894819842471,0.8155039577188894);
\draw[red, thick] (0.6347894819842471,0.8155039577188894)--(0.5646470338381138,0.8144491937677845);
\draw[red, thick] (0.148397993788087,0.28793013543623097)--(0.13656917844354743,0.2001593013620731);
\draw[red, thick] (0.19799452036188725,0.2983935547149186)--(0.17418575103719508,0.3127497753447539);
\draw[red, thick] (0.19799452036188725,0.2983935547149186)--(0.39386035741641695,0.3879950389869909);
\draw[red, thick] (0.2614932165060834,0.4483621656048845)--(0.25909724262654277,0.44854516125233673);
\draw[red, thick] (0.2614932165060834,0.4483621656048845)--(0.4157501376828479,0.4287256297513061);
\draw[red, thick] (0.18588103912022907,0.40276003273976374)--(0.25909724262654277,0.44854516125233673);
\draw[red, thick] (0.4157501376828479,0.4287256297513061)--(0.39386035741641695,0.3879950389869909);
\draw[red, thick] (0.2125410814435029,0.932206125956973)--(0.1503709212537891,0.9000494683595954);
\draw[red, thick] (0.2125410814435029,0.932206125956973)--(0.049259228651907655,0.9730210070877323);
\draw[red, thick] (0.1503709212537891,0.9000494683595954)--(0.07950484344558738,0.9060196908690696);
\draw[red, thick] (0.07950484344558738,0.9060196908690696)--(0.03418643383988759,0.9411415638172064);
\draw[red, thick] (0.049259228651907655,0.9730210070877323)--(0.03418643383988759,0.9411415638172064);
\draw[red, thick] (0.3078280725794059,0.9664502126130512)--(0.2125410814435029,0.932206125956973);
\draw[red, thick] (0.0072051335010195755,0.8663122452666777)--(0.07950484344558738,0.9060196908690696);
\draw[red, thick] (0.07225031953352958,0.7988984371260657)--(0.1503709212537891,0.9000494683595954);
\draw[red, thick] (0.20785657217564413,0.7391137029871157)--(0.299736887125702,0.8166411601178638);
\draw[red, thick] (0.20785657217564413,0.7391137029871157)--(0.1626026765847196,0.6149341693854876);
\draw[red, thick] (0.299736887125702,0.8166411601178638)--(0.410824171102973,0.6379558481429204);
\draw[red, thick] (0.1626026765847196,0.6149341693854876)--(0.21771196590924116,0.5929681031496008);
\draw[red, thick] (0.21771196590924116,0.5929681031496008)--(0.3581748076390561,0.6189381158954018);
\draw[red, thick] (0.3581748076390561,0.6189381158954018)--(0.37128451309713617,0.6232642739680077);
\draw[red, thick] (0.410824171102973,0.6379558481429204)--(0.37128451309713617,0.6232642739680077);
\draw[red, thick] (0.3078280725794059,0.9664502126130512)--(0.299736887125702,0.8166411601178638);
\draw[red, thick] (0.06964041109724696,0.7722567622110954)--(0.20785657217564413,0.7391137029871157);
\draw[red, thick] (0.13638721173493218,0.6043853598410402)--(0.1626026765847196,0.6149341693854876);
\draw[red, thick] (0.5646470338381138,0.8144491937677845)--(0.4573568415412147,0.6126174690235537);
\draw[red, thick] (0.4573568415412147,0.6126174690235537)--(0.410824171102973,0.6379558481429204);
\draw[red, thick] (0.25909724262654277,0.44854516125233673)--(0.2553892027923351,0.48539597554160624);
\draw[red, thick] (0.2553892027923351,0.48539597554160624)--(0.21771196590924116,0.5929681031496008);
\draw[red, thick] (0.2553892027923351,0.48539597554160624)--(0.3581748076390561,0.6189381158954018);
\draw[red, thick] (0.2614932165060834,0.4483621656048845)--(0.37128451309713617,0.6232642739680077);
\draw[red, thick] (0.5092587354725757,0.5327225311464738)--(0.4573568415412147,0.6126174690235537);
\draw[red, thick] (0.5092587354725757,0.5327225311464738)--(0.4959932207799768,0.47088604498338843);
\draw[red, thick] (0.4959932207799768,0.47088604498338843)--(0.4157501376828479,0.4287256297513061);
\draw[red, thick] (0.45995291048752523,0.028965958412618953)--(0.33570477350580863,0.0944446314103765);
\draw[red, thick] (0.13656917844354743,0.2001593013620731)--(0.22603564326488124,0.11592720692047373);
\draw[red, thick] (0.679280231333423,0.07685494809943949)--(0.7039360788752047,0.13691025596481718);
\draw[red, thick] (0.679280231333423,0.07685494809943949)--(0.5401162699392226,0.10592861733374764);
\draw[red, thick] (0.7039360788752047,0.13691025596481718)--(0.692984094105217,0.16657177838006781);
\draw[red, thick] (0.692984094105217,0.16657177838006781)--(0.5705030851078778,0.14956730017396283);
\draw[red, thick] (0.5705030851078778,0.14956730017396283)--(0.5401162699392226,0.10592861733374764);
\draw[red, thick] (0.6704261078411188,0.2983594041110991)--(0.692984094105217,0.16657177838006781);
\draw[red, thick] (0.6704261078411188,0.2983594041110991)--(0.609626416494816,0.3156521591180685);
\draw[red, thick] (0.5705030851078778,0.14956730017396283)--(0.609626416494816,0.3156521591180685);
\draw[red, thick] (0.45995291048752523,0.028965958412618953)--(0.4765623254864447,0.08494160121449679);
\draw[red, thick] (0.4765623254864447,0.08494160121449679)--(0.5401162699392226,0.10592861733374764);
\draw[red, thick] (0.38150587398438834,0.343981855497234)--(0.4497840091409947,0.1953790520155802);
\draw[red, thick] (0.38150587398438834,0.343981855497234)--(0.29113000848283804,0.21354918994579553);
\draw[red, thick] (0.297180677272192,0.20845266990432376)--(0.2922904993448189,0.21081962561616713);
\draw[red, thick] (0.297180677272192,0.20845266990432376)--(0.4415822666940058,0.18039376697821286);
\draw[red, thick] (0.4415822666940058,0.18039376697821286)--(0.4497840091409947,0.1953790520155802);
\draw[red, thick] (0.29113000848283804,0.21354918994579553)--(0.2922904993448189,0.21081962561616713);
\draw[red, thick] (0.4959932207799768,0.47088604498338843)--(0.5773694923567108,0.33302816738559193);
\draw[red, thick] (0.39386035741641695,0.3879950389869909)--(0.38150587398438834,0.343981855497234);
\draw[red, thick] (0.5773694923567108,0.33302816738559193)--(0.4497840091409947,0.1953790520155802);
\draw[red, thick] (0.19799452036188725,0.2983935547149186)--(0.29113000848283804,0.21354918994579553);
\draw[red, thick] (0.33570477350580863,0.0944446314103765)--(0.297180677272192,0.20845266990432376);
\draw[red, thick] (0.22603564326488124,0.11592720692047373)--(0.2922904993448189,0.21081962561616713);
\draw[red, thick] (0.4765623254864447,0.08494160121449679)--(0.4415822666940058,0.18039376697821286);
\draw[red, thick] (0.609626416494816,0.3156521591180685)--(0.5773694923567108,0.33302816738559193);
\draw[red, thick] (0.6966078171943033,0.517084675521012)--(0.8394170635303543,0.5386328082928883);
\draw[red, thick] (0.6966078171943033,0.517084675521012)--(0.7629987978835259,0.4097388057906942);
\draw[red, thick] (0.8394170635303543,0.5386328082928883)--(0.8199448304662863,0.43374447248530434);
\draw[red, thick] (0.7629987978835259,0.4097388057906942)--(0.8199448304662863,0.43374447248530434);
\draw[red, thick] (0.5092587354725757,0.5327225311464738)--(0.5854609328882098,0.5563630707966651);
\draw[red, thick] (0.6704261078411188,0.2983594041110991)--(0.706255398851686,0.31142037155747987);
\draw[red, thick] (0.706255398851686,0.31142037155747987)--(0.7629987978835259,0.4097388057906942);
\draw[red, thick] (0.6966078171943033,0.517084675521012)--(0.6791059840045244,0.528146008831672);
\draw[red, thick] (0.5854609328882098,0.5563630707966651)--(0.6791059840045244,0.528146008831672);
\draw[red, thick] (0.8394170635303543,0.5386328082928883)--(0.8421316549547286,0.5424862865901321);
\draw[red, thick] (0.6791059840045244,0.528146008831672)--(0.7961709115289796,0.6183095407912882);
\draw[red, thick] (0.8421316549547286,0.5424862865901321)--(0.7961709115289796,0.6183095407912882);
\draw[red, thick] (0.8180743974620243,0.1360181565046783)--(0.7847271846057482,0.2833695475015336);
\draw[red, thick] (0.7039360788752047,0.13691025596481718)--(0.8180743974620243,0.1360181565046783);
\draw[red, thick] (0.706255398851686,0.31142037155747987)--(0.7847271846057482,0.2833695475015336);
\draw[red, thick] (0.8340383503796273,0.7758938836779715)--(0.8378339831936373,0.7542743252554434);
\draw[red, thick] (0.663654233906528,0.6697057286556454)--(0.6320097831139766,0.6120438657752406);
\draw[red, thick] (0.663654233906528,0.6697057286556454)--(0.8042869452731249,0.7046892265869095);
\draw[red, thick] (0.6320097831139766,0.6120438657752406)--(0.7954628450760505,0.6220772182985786);
\draw[red, thick] (0.8042869452731249,0.7046892265869095)--(0.7954628450760505,0.6220772182985786);
\draw[red, thick] (0.6347894819842471,0.8155039577188894)--(0.663654233906528,0.6697057286556454);
\draw[red, thick] (0.5854609328882098,0.5563630707966651)--(0.6320097831139766,0.6120438657752406);
\draw[red, thick] (0.8378339831936373,0.7542743252554434)--(0.8042869452731249,0.7046892265869095);
\draw[red, thick] (0.7961709115289796,0.6183095407912882)--(0.7954628450760505,0.6220772182985786);
\draw[black, very thick] (0.1813061501361128,0.42139697266953957)--(0.25616296056720894,0.325003496798576);
\draw[black, very thick] (0.1813061501361128,0.42139697266953957)--(0.25532816389877266,0.4461882000106892);
\draw[black, very thick] (0.1813061501361128,0.42139697266953957)--(0.2585381570265892,0.4541014013513852);
\draw[black, very thick] (0.14872563829292249,0.5541217325950919)--(0.18134071416677522,0.6074653564194161);
\draw[black, very thick] (0.14872563829292249,0.5541217325950919)--(0.1726754312359303,0.6109192601286578);
\draw[black, very thick] (0.05789114625497006,0.7591447714783306)--(0.025416319820836576,0.7229033871050434);
\draw[black, very thick] (0.05789114625497006,0.7591447714783306)--(0.07059719372583992,0.7820235001551941);
\draw[black, very thick] (0.05789114625497006,0.7591447714783306)--(0.06284152846679965,0.9189337987707326);
\draw[black, very thick] (0.05789114625497006,0.7591447714783306)--(0.18563030063821997,0.7444433737589387);
\draw[black, very thick] (0.025416319820836576,0.7229033871050434)--(0.07059719372583992,0.7820235001551941);
\draw[black, very thick] (0.07059719372583992,0.7820235001551941)--(0.06284152846679965,0.9189337987707326);
\draw[black, very thick] (0.07059719372583992,0.7820235001551941)--(0.18563030063821997,0.7444433737589387);
\draw[black, very thick] (0.816871669547315,0.8207650844535032)--(0.797244196123064,0.8483019422161482);
\draw[black, very thick] (0.816871669547315,0.8207650844535032)--(0.7024556644193323,0.853730955492912);
\draw[black, very thick] (0.816871669547315,0.8207650844535032)--(0.84050884537445,0.7523984527019172);
\draw[black, very thick] (0.797244196123064,0.8483019422161482)--(0.7024556644193323,0.853730955492912);
\draw[black, very thick] (0.5150163133378375,0.8659430761214151)--(0.47164836244225017,0.9109390813944592);
\draw[black, very thick] (0.5150163133378375,0.8659430761214151)--(0.5337666871494476,0.7563578295951302);
\draw[black, very thick] (0.47164836244225017,0.9109390813944592)--(0.38653207447931925,0.9992506770062238);
\draw[black, very thick] (0.47164836244225017,0.9109390813944592)--(0.5337666871494476,0.7563578295951302);
\draw[black, very thick] (0.38653207447931925,0.9992506770062238)--(0.2921345429365123,0.9608102969999692);
\draw[black, very thick] (0.38653207447931925,0.9992506770062238)--(0.30700725570527704,0.9512527117166875);
\draw[black, very thick] (0.29224699803261267,0.34151062925235404)--(0.25616296056720894,0.325003496798576);
\draw[black, very thick] (0.29224699803261267,0.34151062925235404)--(0.35675114120119855,0.3082552180370413);
\draw[black, very thick] (0.29224699803261267,0.34151062925235404)--(0.37932348476575223,0.34083217785400804);
\draw[black, very thick] (0.25616296056720894,0.325003496798576)--(0.35675114120119855,0.3082552180370413);
\draw[black, very thick] (0.3693375759616444,0.43463383746958717)--(0.37181151999371226,0.4343189103317002);
\draw[black, very thick] (0.3693375759616444,0.43463383746958717)--(0.4100600707523283,0.4294499616155625);
\draw[black, very thick] (0.3693375759616444,0.43463383746958717)--(0.4169904242473743,0.42937728712224627);
\draw[black, very thick] (0.37181151999371226,0.4343189103317002)--(0.4100600707523283,0.4294499616155625);
\draw[black, very thick] (0.37181151999371226,0.4343189103317002)--(0.4169904242473743,0.42937728712224627);
\draw[black, very thick] (0.37181151999371226,0.4343189103317002)--(0.37932348476575223,0.34083217785400804);
\draw[black, very thick] (0.37181151999371226,0.4343189103317002)--(0.3849643126762702,0.3563026374685184);
\draw[black, very thick] (0.4100600707523283,0.4294499616155625)--(0.4169904242473743,0.42937728712224627);
\draw[black, very thick] (0.4100600707523283,0.4294499616155625)--(0.49024128124118205,0.4678639258176544);
\draw[black, very thick] (0.4100600707523283,0.4294499616155625)--(0.3849643126762702,0.3563026374685184);
\draw[black, very thick] (0.25532816389877266,0.4461882000106892)--(0.2585381570265892,0.4541014013513852);
\draw[black, very thick] (0.25532816389877266,0.4461882000106892)--(0.26648526291374264,0.4998123096218199);
\draw[black, very thick] (0.25532816389877266,0.4461882000106892)--(0.2919655586376268,0.49690587581737233);
\draw[black, very thick] (0.15937905866236854,0.9047088031490487)--(0.1590067661594326,0.9455878746166037);
\draw[black, very thick] (0.15937905866236854,0.9047088031490487)--(0.17361789056520271,0.9419355930212361);
\draw[black, very thick] (0.15937905866236854,0.9047088031490487)--(0.06284152846679965,0.9189337987707326);
\draw[black, very thick] (0.1590067661594326,0.9455878746166037)--(0.17361789056520271,0.9419355930212361);
\draw[black, very thick] (0.1590067661594326,0.9455878746166037)--(0.06284152846679965,0.9189337987707326);
\draw[black, very thick] (0.17361789056520271,0.9419355930212361)--(0.06284152846679965,0.9189337987707326);
\draw[black, very thick] (0.17361789056520271,0.9419355930212361)--(0.2921345429365123,0.9608102969999692);
\draw[black, very thick] (0.2921345429365123,0.9608102969999692)--(0.30700725570527704,0.9512527117166875);
\draw[black, very thick] (0.2921345429365123,0.9608102969999692)--(0.304643617111923,0.9074897229510144);
\draw[black, very thick] (0.2526856156454257,0.7769398874936125)--(0.27892638396565794,0.7990815186159101);
\draw[black, very thick] (0.2526856156454257,0.7769398874936125)--(0.3010446550260852,0.8408546046323935);
\draw[black, very thick] (0.2526856156454257,0.7769398874936125)--(0.29981130217595947,0.8180189617063004);
\draw[black, very thick] (0.2526856156454257,0.7769398874936125)--(0.30055791219730454,0.8318425158098488);
\draw[black, very thick] (0.2526856156454257,0.7769398874936125)--(0.18563030063821997,0.7444433737589387);
\draw[black, very thick] (0.27892638396565794,0.7990815186159101)--(0.3010446550260852,0.8408546046323935);
\draw[black, very thick] (0.27892638396565794,0.7990815186159101)--(0.29981130217595947,0.8180189617063004);
\draw[black, very thick] (0.27892638396565794,0.7990815186159101)--(0.30055791219730454,0.8318425158098488);
\draw[black, very thick] (0.27892638396565794,0.7990815186159101)--(0.18563030063821997,0.7444433737589387);
\draw[black, very thick] (0.19387030084463816,0.7007344954508157)--(0.18563030063821997,0.7444433737589387);
\draw[black, very thick] (0.37293603650876134,0.6988993954844914)--(0.3982187058250188,0.6582318996866904);
\draw[black, very thick] (0.37293603650876134,0.6988993954844914)--(0.34509732830012835,0.6165202358288133);
\draw[black, very thick] (0.37293603650876134,0.6988993954844914)--(0.36683993239497326,0.6217975775006084);
\draw[black, very thick] (0.37293603650876134,0.6988993954844914)--(0.4828113417610102,0.6605018615470278);
\draw[black, very thick] (0.3982187058250188,0.6582318996866904)--(0.34509732830012835,0.6165202358288133);
\draw[black, very thick] (0.3982187058250188,0.6582318996866904)--(0.36683993239497326,0.6217975775006084);
\draw[black, very thick] (0.3982187058250188,0.6582318996866904)--(0.4828113417610102,0.6605018615470278);
\draw[black, very thick] (0.3982187058250188,0.6582318996866904)--(0.3498485323164778,0.5891158637590141);
\draw[black, very thick] (0.18134071416677522,0.6074653564194161)--(0.1726754312359303,0.6109192601286578);
\draw[black, very thick] (0.18134071416677522,0.6074653564194161)--(0.232149332632162,0.5956374112256204);
\draw[black, very thick] (0.18134071416677522,0.6074653564194161)--(0.23156648132676858,0.5534121318594152);
\draw[black, very thick] (0.18134071416677522,0.6074653564194161)--(0.2277040583383429,0.5644397202114);
\draw[black, very thick] (0.1726754312359303,0.6109192601286578)--(0.232149332632162,0.5956374112256204);
\draw[black, very thick] (0.1726754312359303,0.6109192601286578)--(0.2277040583383429,0.5644397202114);
\draw[black, very thick] (0.34509732830012835,0.6165202358288133)--(0.36683993239497326,0.6217975775006084);
\draw[black, very thick] (0.34509732830012835,0.6165202358288133)--(0.3498485323164778,0.5891158637590141);
\draw[black, very thick] (0.34509732830012835,0.6165202358288133)--(0.34045815137782165,0.574156595807007);
\draw[black, very thick] (0.232149332632162,0.5956374112256204)--(0.23156648132676858,0.5534121318594152);
\draw[black, very thick] (0.232149332632162,0.5956374112256204)--(0.2277040583383429,0.5644397202114);
\draw[black, very thick] (0.36683993239497326,0.6217975775006084)--(0.3498485323164778,0.5891158637590141);
\draw[black, very thick] (0.36683993239497326,0.6217975775006084)--(0.34045815137782165,0.574156595807007);
\draw[black, very thick] (0.3010446550260852,0.8408546046323935)--(0.304643617111923,0.9074897229510144);
\draw[black, very thick] (0.3010446550260852,0.8408546046323935)--(0.29981130217595947,0.8180189617063004);
\draw[black, very thick] (0.3010446550260852,0.8408546046323935)--(0.30055791219730454,0.8318425158098488);
\draw[black, very thick] (0.30700725570527704,0.9512527117166875)--(0.304643617111923,0.9074897229510144);
\draw[black, very thick] (0.304643617111923,0.9074897229510144)--(0.29981130217595947,0.8180189617063004);
\draw[black, very thick] (0.304643617111923,0.9074897229510144)--(0.30055791219730454,0.8318425158098488);
\draw[black, very thick] (0.29981130217595947,0.8180189617063004)--(0.30055791219730454,0.8318425158098488);
\draw[black, very thick] (0.5337666871494476,0.7563578295951302)--(0.4828113417610102,0.6605018615470278);
\draw[black, very thick] (0.2585381570265892,0.4541014013513852)--(0.26648526291374264,0.4998123096218199);
\draw[black, very thick] (0.2585381570265892,0.4541014013513852)--(0.2919655586376268,0.49690587581737233);
\draw[black, very thick] (0.23156648132676858,0.5534121318594152)--(0.2277040583383429,0.5644397202114);
\draw[black, very thick] (0.23156648132676858,0.5534121318594152)--(0.26648526291374264,0.4998123096218199);
\draw[black, very thick] (0.23156648132676858,0.5534121318594152)--(0.293302842223799,0.5346545129197736);
\draw[black, very thick] (0.2277040583383429,0.5644397202114)--(0.293302842223799,0.5346545129197736);
\draw[black, very thick] (0.26648526291374264,0.4998123096218199)--(0.293302842223799,0.5346545129197736);
\draw[black, very thick] (0.26648526291374264,0.4998123096218199)--(0.2919655586376268,0.49690587581737233);
\draw[black, very thick] (0.293302842223799,0.5346545129197736)--(0.2919655586376268,0.49690587581737233);
\draw[black, very thick] (0.293302842223799,0.5346545129197736)--(0.34045815137782165,0.574156595807007);
\draw[black, very thick] (0.3498485323164778,0.5891158637590141)--(0.34045815137782165,0.574156595807007);
\draw[black, very thick] (0.5084688879484911,0.5290406992216823)--(0.49024128124118205,0.4678639258176544);
\draw[black, very thick] (0.5084688879484911,0.5290406992216823)--(0.5800834392674282,0.554694787525681);
\draw[black, very thick] (0.5084688879484911,0.5290406992216823)--(0.6188652259059673,0.5462977123528099);
\draw[black, very thick] (0.4169904242473743,0.42937728712224627)--(0.49024128124118205,0.4678639258176544);
\draw[black, very thick] (0.4169904242473743,0.42937728712224627)--(0.3849643126762702,0.3563026374685184);
\draw[black, very thick] (0.49024128124118205,0.4678639258176544)--(0.5570056583155196,0.3675261222104846);
\draw[black, very thick] (0.49024128124118205,0.4678639258176544)--(0.5326754741324017,0.40874338970209434);
\draw[black, very thick] (0.38355244884048856,0.06922894316707184)--(0.44145691925859737,0.03871333161972445);
\draw[black, very thick] (0.38355244884048856,0.06922894316707184)--(0.37508174221498614,0.07369299906895496);
\draw[black, very thick] (0.38355244884048856,0.06922894316707184)--(0.35676393151346963,0.08334647127462609);
\draw[black, very thick] (0.38355244884048856,0.06922894316707184)--(0.33402090497369247,0.09314723041436684);
\draw[black, very thick] (0.38355244884048856,0.06922894316707184)--(0.4744519247018341,0.0778293073440714);
\draw[black, very thick] (0.44145691925859737,0.03871333161972445)--(0.37508174221498614,0.07369299906895496);
\draw[black, very thick] (0.44145691925859737,0.03871333161972445)--(0.35676393151346963,0.08334647127462609);
\draw[black, very thick] (0.44145691925859737,0.03871333161972445)--(0.4744519247018341,0.0778293073440714);
\draw[black, very thick] (0.37508174221498614,0.07369299906895496)--(0.35676393151346963,0.08334647127462609);
\draw[black, very thick] (0.37508174221498614,0.07369299906895496)--(0.33402090497369247,0.09314723041436684);
\draw[black, very thick] (0.35676393151346963,0.08334647127462609)--(0.27052606386950184,0.04422519848677735);
\draw[black, very thick] (0.35676393151346963,0.08334647127462609)--(0.33402090497369247,0.09314723041436684);
\draw[black, very thick] (0.27052606386950184,0.04422519848677735)--(0.33402090497369247,0.09314723041436684);
\draw[black, very thick] (0.27052606386950184,0.04422519848677735)--(0.25415869288763704,0.031614330708849434);
\draw[black, very thick] (0.27052606386950184,0.04422519848677735)--(0.1872786858015116,0.01718136555514127);
\draw[black, very thick] (0.27052606386950184,0.04422519848677735)--(0.18417778144474048,0.00928081221598432);
\draw[black, very thick] (0.19012226347068706,-0.017724988674975257)--(0.25415869288763704,0.031614330708849434);
\draw[black, very thick] (0.19012226347068706,-0.017724988674975257)--(0.1872786858015116,0.01718136555514127);
\draw[black, very thick] (0.19012226347068706,-0.017724988674975257)--(0.18417778144474048,0.00928081221598432);
\draw[black, very thick] (0.19012226347068706,-0.017724988674975257)--(0.17858300600531926,-0.004973681536212152);
\draw[black, very thick] (0.25415869288763704,0.031614330708849434)--(0.1872786858015116,0.01718136555514127);
\draw[black, very thick] (0.25415869288763704,0.031614330708849434)--(0.18417778144474048,0.00928081221598432);
\draw[black, very thick] (0.25415869288763704,0.031614330708849434)--(0.17858300600531926,-0.004973681536212152);
\draw[black, very thick] (0.19976901885816142,0.14065706287384386)--(0.15362496192050873,0.1841013922061019);
\draw[black, very thick] (0.19976901885816142,0.14065706287384386)--(0.153557661738525,0.18416475488949516);
\draw[black, very thick] (0.19976901885816142,0.14065706287384386)--(0.14846174817708205,0.18896252434680957);
\draw[black, very thick] (0.15362496192050873,0.1841013922061019)--(0.153557661738525,0.18416475488949516);
\draw[black, very thick] (0.15362496192050873,0.1841013922061019)--(0.14846174817708205,0.18896252434680957);
\draw[black, very thick] (0.153557661738525,0.18416475488949516)--(0.14846174817708205,0.18896252434680957);
\draw[black, very thick] (0.1872786858015116,0.01718136555514127)--(0.18417778144474048,0.00928081221598432);
\draw[black, very thick] (0.1872786858015116,0.01718136555514127)--(0.17858300600531926,-0.004973681536212152);
\draw[black, very thick] (0.18417778144474048,0.00928081221598432)--(0.17858300600531926,-0.004973681536212152);
\draw[black, very thick] (0.694675762855622,0.11435450611503778)--(0.6796560855289289,0.2444362909933815);
\draw[black, very thick] (0.694675762855622,0.11435450611503778)--(0.8118874384593696,0.1633564836689559);
\draw[black, very thick] (0.6796560855289289,0.2444362909933815)--(0.8118874384593696,0.1633564836689559);
\draw[black, very thick] (0.4744519247018341,0.0778293073440714)--(0.5263999923459467,0.10139917795079811);
\draw[black, very thick] (0.4744519247018341,0.0778293073440714)--(0.45877901890888895,0.13346796697448662);
\draw[black, very thick] (0.44625527759791844,0.20305910066999244)--(0.43025304052789864,0.18259516696622943);
\draw[black, very thick] (0.44625527759791844,0.20305910066999244)--(0.48139088714400735,0.22947900357015216);
\draw[black, very thick] (0.44625527759791844,0.20305910066999244)--(0.4952009553583087,0.24437837632781984);
\draw[black, very thick] (0.44625527759791844,0.20305910066999244)--(0.45162040934298203,0.19736030299255442);
\draw[black, very thick] (0.44625527759791844,0.20305910066999244)--(0.44500451138548675,0.17105528445839943);
\draw[black, very thick] (0.35675114120119855,0.3082552180370413)--(0.37932348476575223,0.34083217785400804);
\draw[black, very thick] (0.35675114120119855,0.3082552180370413)--(0.3849643126762702,0.3563026374685184);
\draw[black, very thick] (0.37932348476575223,0.34083217785400804)--(0.3849643126762702,0.3563026374685184);
\draw[black, very thick] (0.3781940958333602,0.1927108234591293)--(0.43025304052789864,0.18259516696622943);
\draw[black, very thick] (0.43025304052789864,0.18259516696622943)--(0.45162040934298203,0.19736030299255442);
\draw[black, very thick] (0.43025304052789864,0.18259516696622943)--(0.44500451138548675,0.17105528445839943);
\draw[black, very thick] (0.43025304052789864,0.18259516696622943)--(0.45877901890888895,0.13346796697448662);
\draw[black, very thick] (0.5646392290876473,0.35459424606590373)--(0.5570056583155196,0.3675261222104846);
\draw[black, very thick] (0.5646392290876473,0.35459424606590373)--(0.5326754741324017,0.40874338970209434);
\draw[black, very thick] (0.5646392290876473,0.35459424606590373)--(0.5464741713477421,0.2996958987076048);
\draw[black, very thick] (0.5570056583155196,0.3675261222104846)--(0.5326754741324017,0.40874338970209434);
\draw[black, very thick] (0.5570056583155196,0.3675261222104846)--(0.5464741713477421,0.2996958987076048);
\draw[black, very thick] (0.5326754741324017,0.40874338970209434)--(0.5464741713477421,0.2996958987076048);
\draw[black, very thick] (0.48139088714400735,0.22947900357015216)--(0.4952009553583087,0.24437837632781984);
\draw[black, very thick] (0.48139088714400735,0.22947900357015216)--(0.45162040934298203,0.19736030299255442);
\draw[black, very thick] (0.48139088714400735,0.22947900357015216)--(0.44500451138548675,0.17105528445839943);
\draw[black, very thick] (0.4952009553583087,0.24437837632781984)--(0.5464741713477421,0.2996958987076048);
\draw[black, very thick] (0.4952009553583087,0.24437837632781984)--(0.45162040934298203,0.19736030299255442);
\draw[black, very thick] (0.45162040934298203,0.19736030299255442)--(0.44500451138548675,0.17105528445839943);
\draw[black, very thick] (0.44500451138548675,0.17105528445839943)--(0.45877901890888895,0.13346796697448662);
\draw[black, very thick] (0.7221119059162853,0.4758477701896109)--(0.8272393746070328,0.4730369662060027);
\draw[black, very thick] (0.7221119059162853,0.4758477701896109)--(0.7711711320631728,0.41318386283041797);
\draw[black, very thick] (0.7221119059162853,0.4758477701896109)--(0.6361928686502474,0.5410765590383083);
\draw[black, very thick] (0.8272393746070328,0.4730369662060027)--(0.7711711320631728,0.41318386283041797);
\draw[black, very thick] (0.8272393746070328,0.4730369662060027)--(0.8851632825205141,0.5493531197369946);
\draw[black, very thick] (0.5800834392674282,0.554694787525681)--(0.6188652259059673,0.5462977123528099);
\draw[black, very thick] (0.5800834392674282,0.554694787525681)--(0.6361928686502474,0.5410765590383083);
\draw[black, very thick] (0.5800834392674282,0.554694787525681)--(0.6273026388087081,0.6064132747783626);
\draw[black, very thick] (0.6188652259059673,0.5462977123528099)--(0.6361928686502474,0.5410765590383083);
\draw[black, very thick] (0.6188652259059673,0.5462977123528099)--(0.6273026388087081,0.6064132747783626);
\draw[black, very thick] (0.6361928686502474,0.5410765590383083)--(0.6273026388087081,0.6064132747783626);
\draw[black, very thick] (0.7880511358690263,0.6120556809230696)--(0.7843733421016718,0.6092230402432683);
\draw[black, very thick] (0.7880511358690263,0.6120556809230696)--(0.7940377748297177,0.6166665966317921);
\draw[black, very thick] (0.7880511358690263,0.6120556809230696)--(0.7851149531728918,0.6097942298808102);
\draw[black, very thick] (0.7880511358690263,0.6120556809230696)--(0.8108312686319776,0.5941237704915503);
\draw[black, very thick] (0.7880511358690263,0.6120556809230696)--(0.7860546895085898,0.6214997109685254);
\draw[black, very thick] (0.7843733421016718,0.6092230402432683)--(0.7940377748297177,0.6166665966317921);
\draw[black, very thick] (0.7843733421016718,0.6092230402432683)--(0.7851149531728918,0.6097942298808102);
\draw[black, very thick] (0.7843733421016718,0.6092230402432683)--(0.8108312686319776,0.5941237704915503);
\draw[black, very thick] (0.7843733421016718,0.6092230402432683)--(0.7860546895085898,0.6214997109685254);
\draw[black, very thick] (0.7940377748297177,0.6166665966317921)--(0.7851149531728918,0.6097942298808102);
\draw[black, very thick] (0.7940377748297177,0.6166665966317921)--(0.8108312686319776,0.5941237704915503);
\draw[black, very thick] (0.7940377748297177,0.6166665966317921)--(0.7860546895085898,0.6214997109685254);
\draw[black, very thick] (0.7851149531728918,0.6097942298808102)--(0.8108312686319776,0.5941237704915503);
\draw[black, very thick] (0.7851149531728918,0.6097942298808102)--(0.7860546895085898,0.6214997109685254);
\draw[black, very thick] (0.8108312686319776,0.5941237704915503)--(0.7860546895085898,0.6214997109685254);
\draw[black, very thick] (0.8851632825205141,0.5493531197369946)--(0.9421973439307684,0.5584544111718442);
\draw[black, very thick] (1.0241664694608577,0.57153475091379)--(1.0378006206639054,0.5737104399756852);
\draw[black, very thick] (1.0241664694608577,0.57153475091379)--(1.0823433206558173,0.580818404653238);
\draw[black, very thick] (1.0241664694608577,0.57153475091379)--(1.0605948498489768,0.5980528201790608);
\draw[black, very thick] (1.0241664694608577,0.57153475091379)--(1.0418135837704467,0.6112240638453353);
\draw[black, very thick] (1.0241664694608577,0.57153475091379)--(1.0588598410834567,0.5992695764466222);
\draw[black, very thick] (1.0378006206639054,0.5737104399756852)--(1.0823433206558173,0.580818404653238);
\draw[black, very thick] (1.0378006206639054,0.5737104399756852)--(1.0605948498489768,0.5980528201790608);
\draw[black, very thick] (1.0378006206639054,0.5737104399756852)--(1.0418135837704467,0.6112240638453353);
\draw[black, very thick] (1.0378006206639054,0.5737104399756852)--(1.0588598410834567,0.5992695764466222);
\draw[black, very thick] (1.0823433206558173,0.580818404653238)--(1.0605948498489768,0.5980528201790608);
\draw[black, very thick] (1.0823433206558173,0.580818404653238)--(1.0418135837704467,0.6112240638453353);
\draw[black, very thick] (1.0823433206558173,0.580818404653238)--(1.0588598410834567,0.5992695764466222);
\draw[black, very thick] (0.9421973439307684,0.5584544111718442)--(0.9625363197915181,0.6668209632607429);
\draw[black, very thick] (1.0605948498489768,0.5980528201790608)--(1.0418135837704467,0.6112240638453353);
\draw[black, very thick] (1.0605948498489768,0.5980528201790608)--(1.0588598410834567,0.5992695764466222);
\draw[black, very thick] (0.9625363197915181,0.6668209632607429)--(0.9168260577530807,0.6988774282956411);
\draw[black, very thick] (0.9625363197915181,0.6668209632607429)--(0.9338936187195976,0.6869080007107792);
\draw[black, very thick] (1.0418135837704467,0.6112240638453353)--(1.0588598410834567,0.5992695764466222);
\draw[black, very thick] (0.9168260577530807,0.6988774282956411)--(0.9338936187195976,0.6869080007107792);
\draw[black, very thick] (0.9168260577530807,0.6988774282956411)--(0.84050884537445,0.7523984527019172);
\draw[black, very thick] (0.84050884537445,0.7523984527019172)--(0.8032960357300718,0.6954122436021932);
\draw[black, very thick] (0.642114287806097,0.6304560849618224)--(0.6595268144771232,0.6905536626697233);
\draw[black, very thick] (0.642114287806097,0.6304560849618224)--(0.6273026388087081,0.6064132747783626);
\draw[black, very thick] (0.7630519194473339,0.6944316881617858)--(0.7580659591276837,0.6931913911438812);
\draw[black, very thick] (0.7630519194473339,0.6944316881617858)--(0.7097088513509334,0.6811621785700015);
\draw[black, very thick] (0.7630519194473339,0.6944316881617858)--(0.6984784752351365,0.6783685337927906);
\draw[black, very thick] (0.7630519194473339,0.6944316881617858)--(0.8032960357300718,0.6954122436021932);
\draw[black, very thick] (0.7580659591276837,0.6931913911438812)--(0.7097088513509334,0.6811621785700015);
\draw[black, very thick] (0.7580659591276837,0.6931913911438812)--(0.6984784752351365,0.6783685337927906);
\draw[black, very thick] (0.7580659591276837,0.6931913911438812)--(0.8032960357300718,0.6954122436021932);
\draw[black, very thick] (0.7097088513509334,0.6811621785700015)--(0.6984784752351365,0.6783685337927906);
\draw[black, very thick] (0.7097088513509334,0.6811621785700015)--(0.714648332376192,0.6171165250868014);
\draw[black, very thick] (0.7097088513509334,0.6811621785700015)--(0.6595268144771232,0.6905536626697233);
\draw[black, very thick] (0.6984784752351365,0.6783685337927906)--(0.714648332376192,0.6171165250868014);
\draw[black, very thick] (0.6984784752351365,0.6783685337927906)--(0.6595268144771232,0.6905536626697233);
\fill[blue] (0.1813061501361128,0.42139697266953957) circle (0.1pt);\draw(0.1813061501361128,0.42139697266953957) circle (0.1pt);
\fill[blue] (0.14872563829292249,0.5541217325950919) circle (0.1pt);\draw(0.14872563829292249,0.5541217325950919) circle (0.1pt);
\fill[blue] (0.05789114625497006,0.7591447714783306) circle (0.1pt);\draw(0.05789114625497006,0.7591447714783306) circle (0.1pt);
\fill[blue] (0.025416319820836576,0.7229033871050434) circle (0.1pt);\draw(0.025416319820836576,0.7229033871050434) circle (0.1pt);
\fill[blue] (0.07059719372583992,0.7820235001551941) circle (0.1pt);\draw(0.07059719372583992,0.7820235001551941) circle (0.1pt);
\fill[blue] (0.816871669547315,0.8207650844535032) circle (0.1pt);\draw(0.816871669547315,0.8207650844535032) circle (0.1pt);
\fill[blue] (0.797244196123064,0.8483019422161482) circle (0.1pt);\draw(0.797244196123064,0.8483019422161482) circle (0.1pt);
\fill[blue] (0.5150163133378375,0.8659430761214151) circle (0.1pt);\draw(0.5150163133378375,0.8659430761214151) circle (0.1pt);
\fill[blue] (0.47164836244225017,0.9109390813944592) circle (0.1pt);\draw(0.47164836244225017,0.9109390813944592) circle (0.1pt);
\fill[blue] (0.38653207447931925,0.9992506770062238) circle (0.1pt);\draw(0.38653207447931925,0.9992506770062238) circle (0.1pt);
\fill[blue] (0.7024556644193323,0.853730955492912) circle (0.1pt);\draw(0.7024556644193323,0.853730955492912) circle (0.1pt);
\fill[blue] (0.29224699803261267,0.34151062925235404) circle (0.1pt);\draw(0.29224699803261267,0.34151062925235404) circle (0.1pt);
\fill[blue] (0.25616296056720894,0.325003496798576) circle (0.1pt);\draw(0.25616296056720894,0.325003496798576) circle (0.1pt);
\fill[blue] (0.3693375759616444,0.43463383746958717) circle (0.1pt);\draw(0.3693375759616444,0.43463383746958717) circle (0.1pt);
\fill[blue] (0.37181151999371226,0.4343189103317002) circle (0.1pt);\draw(0.37181151999371226,0.4343189103317002) circle (0.1pt);
\fill[blue] (0.4100600707523283,0.4294499616155625) circle (0.1pt);\draw(0.4100600707523283,0.4294499616155625) circle (0.1pt);
\fill[blue] (0.25532816389877266,0.4461882000106892) circle (0.1pt);\draw(0.25532816389877266,0.4461882000106892) circle (0.1pt);
\fill[blue] (0.15937905866236854,0.9047088031490487) circle (0.1pt);\draw(0.15937905866236854,0.9047088031490487) circle (0.1pt);
\fill[blue] (0.1590067661594326,0.9455878746166037) circle (0.1pt);\draw(0.1590067661594326,0.9455878746166037) circle (0.1pt);
\fill[blue] (0.17361789056520271,0.9419355930212361) circle (0.1pt);\draw(0.17361789056520271,0.9419355930212361) circle (0.1pt);
\fill[blue] (0.06284152846679965,0.9189337987707326) circle (0.1pt);\draw(0.06284152846679965,0.9189337987707326) circle (0.1pt);
\fill[blue] (0.2921345429365123,0.9608102969999692) circle (0.1pt);\draw(0.2921345429365123,0.9608102969999692) circle (0.1pt);
\fill[blue] (0.2526856156454257,0.7769398874936125) circle (0.1pt);\draw(0.2526856156454257,0.7769398874936125) circle (0.1pt);
\fill[blue] (0.27892638396565794,0.7990815186159101) circle (0.1pt);\draw(0.27892638396565794,0.7990815186159101) circle (0.1pt);
\fill[blue] (0.19387030084463816,0.7007344954508157) circle (0.1pt);\draw(0.19387030084463816,0.7007344954508157) circle (0.1pt);
\fill[blue] (0.37293603650876134,0.6988993954844914) circle (0.1pt);\draw(0.37293603650876134,0.6988993954844914) circle (0.1pt);
\fill[blue] (0.3982187058250188,0.6582318996866904) circle (0.1pt);\draw(0.3982187058250188,0.6582318996866904) circle (0.1pt);
\fill[blue] (0.18134071416677522,0.6074653564194161) circle (0.1pt);\draw(0.18134071416677522,0.6074653564194161) circle (0.1pt);
\fill[blue] (0.1726754312359303,0.6109192601286578) circle (0.1pt);\draw(0.1726754312359303,0.6109192601286578) circle (0.1pt);
\fill[blue] (0.34509732830012835,0.6165202358288133) circle (0.1pt);\draw(0.34509732830012835,0.6165202358288133) circle (0.1pt);
\fill[blue] (0.232149332632162,0.5956374112256204) circle (0.1pt);\draw(0.232149332632162,0.5956374112256204) circle (0.1pt);
\fill[blue] (0.36683993239497326,0.6217975775006084) circle (0.1pt);\draw(0.36683993239497326,0.6217975775006084) circle (0.1pt);
\fill[blue] (0.3010446550260852,0.8408546046323935) circle (0.1pt);\draw(0.3010446550260852,0.8408546046323935) circle (0.1pt);
\fill[blue] (0.30700725570527704,0.9512527117166875) circle (0.1pt);\draw(0.30700725570527704,0.9512527117166875) circle (0.1pt);
\fill[blue] (0.304643617111923,0.9074897229510144) circle (0.1pt);\draw(0.304643617111923,0.9074897229510144) circle (0.1pt);
\fill[blue] (0.29981130217595947,0.8180189617063004) circle (0.1pt);\draw(0.29981130217595947,0.8180189617063004) circle (0.1pt);
\fill[blue] (0.30055791219730454,0.8318425158098488) circle (0.1pt);\draw(0.30055791219730454,0.8318425158098488) circle (0.1pt);
\fill[blue] (0.18563030063821997,0.7444433737589387) circle (0.1pt);\draw(0.18563030063821997,0.7444433737589387) circle (0.1pt);
\fill[blue] (0.5337666871494476,0.7563578295951302) circle (0.1pt);\draw(0.5337666871494476,0.7563578295951302) circle (0.1pt);
\fill[blue] (0.4828113417610102,0.6605018615470278) circle (0.1pt);\draw(0.4828113417610102,0.6605018615470278) circle (0.1pt);
\fill[blue] (0.2585381570265892,0.4541014013513852) circle (0.1pt);\draw(0.2585381570265892,0.4541014013513852) circle (0.1pt);
\fill[blue] (0.23156648132676858,0.5534121318594152) circle (0.1pt);\draw(0.23156648132676858,0.5534121318594152) circle (0.1pt);
\fill[blue] (0.2277040583383429,0.5644397202114) circle (0.1pt);\draw(0.2277040583383429,0.5644397202114) circle (0.1pt);
\fill[blue] (0.26648526291374264,0.4998123096218199) circle (0.1pt);\draw(0.26648526291374264,0.4998123096218199) circle (0.1pt);
\fill[blue] (0.293302842223799,0.5346545129197736) circle (0.1pt);\draw(0.293302842223799,0.5346545129197736) circle (0.1pt);
\fill[blue] (0.3498485323164778,0.5891158637590141) circle (0.1pt);\draw(0.3498485323164778,0.5891158637590141) circle (0.1pt);
\fill[blue] (0.2919655586376268,0.49690587581737233) circle (0.1pt);\draw(0.2919655586376268,0.49690587581737233) circle (0.1pt);
\fill[blue] (0.34045815137782165,0.574156595807007) circle (0.1pt);\draw(0.34045815137782165,0.574156595807007) circle (0.1pt);
\fill[blue] (0.5084688879484911,0.5290406992216823) circle (0.1pt);\draw(0.5084688879484911,0.5290406992216823) circle (0.1pt);
\fill[blue] (0.4169904242473743,0.42937728712224627) circle (0.1pt);\draw(0.4169904242473743,0.42937728712224627) circle (0.1pt);
\fill[blue] (0.49024128124118205,0.4678639258176544) circle (0.1pt);\draw(0.49024128124118205,0.4678639258176544) circle (0.1pt);
\fill[blue] (0.38355244884048856,0.06922894316707184) circle (0.1pt);\draw(0.38355244884048856,0.06922894316707184) circle (0.1pt);
\fill[blue] (0.44145691925859737,0.03871333161972445) circle (0.1pt);\draw(0.44145691925859737,0.03871333161972445) circle (0.1pt);
\fill[blue] (0.37508174221498614,0.07369299906895496) circle (0.1pt);\draw(0.37508174221498614,0.07369299906895496) circle (0.1pt);
\fill[blue] (0.35676393151346963,0.08334647127462609) circle (0.1pt);\draw(0.35676393151346963,0.08334647127462609) circle (0.1pt);
\fill[blue] (0.27052606386950184,0.04422519848677735) circle (0.1pt);\draw(0.27052606386950184,0.04422519848677735) circle (0.1pt);
\fill[blue] (0.33402090497369247,0.09314723041436684) circle (0.1pt);\draw(0.33402090497369247,0.09314723041436684) circle (0.1pt);
\fill[blue] (0.19012226347068706,-0.017724988674975257) circle (0.1pt);\draw(0.19012226347068706,-0.017724988674975257) circle (0.1pt);
\fill[blue] (0.25415869288763704,0.031614330708849434) circle (0.1pt);\draw(0.25415869288763704,0.031614330708849434) circle (0.1pt);
\fill[blue] (0.19976901885816142,0.14065706287384386) circle (0.1pt);\draw(0.19976901885816142,0.14065706287384386) circle (0.1pt);
\fill[blue] (0.15362496192050873,0.1841013922061019) circle (0.1pt);\draw(0.15362496192050873,0.1841013922061019) circle (0.1pt);
\fill[blue] (0.153557661738525,0.18416475488949516) circle (0.1pt);\draw(0.153557661738525,0.18416475488949516) circle (0.1pt);
\fill[blue] (0.14846174817708205,0.18896252434680957) circle (0.1pt);\draw(0.14846174817708205,0.18896252434680957) circle (0.1pt);
\fill[blue] (0.1872786858015116,0.01718136555514127) circle (0.1pt);\draw(0.1872786858015116,0.01718136555514127) circle (0.1pt);
\fill[blue] (0.18417778144474048,0.00928081221598432) circle (0.1pt);\draw(0.18417778144474048,0.00928081221598432) circle (0.1pt);
\fill[blue] (0.17858300600531926,-0.004973681536212152) circle (0.1pt);\draw(0.17858300600531926,-0.004973681536212152) circle (0.1pt);
\fill[blue] (0.694675762855622,0.11435450611503778) circle (0.1pt);\draw(0.694675762855622,0.11435450611503778) circle (0.1pt);
\fill[blue] (0.6796560855289289,0.2444362909933815) circle (0.1pt);\draw(0.6796560855289289,0.2444362909933815) circle (0.1pt);
\fill[blue] (0.4744519247018341,0.0778293073440714) circle (0.1pt);\draw(0.4744519247018341,0.0778293073440714) circle (0.1pt);
\fill[blue] (0.5263999923459467,0.10139917795079811) circle (0.1pt);\draw(0.5263999923459467,0.10139917795079811) circle (0.1pt);
\fill[blue] (0.44625527759791844,0.20305910066999244) circle (0.1pt);\draw(0.44625527759791844,0.20305910066999244) circle (0.1pt);
\fill[blue] (0.35675114120119855,0.3082552180370413) circle (0.1pt);\draw(0.35675114120119855,0.3082552180370413) circle (0.1pt);
\fill[blue] (0.37932348476575223,0.34083217785400804) circle (0.1pt);\draw(0.37932348476575223,0.34083217785400804) circle (0.1pt);
\fill[blue] (0.3781940958333602,0.1927108234591293) circle (0.1pt);\draw(0.3781940958333602,0.1927108234591293) circle (0.1pt);
\fill[blue] (0.43025304052789864,0.18259516696622943) circle (0.1pt);\draw(0.43025304052789864,0.18259516696622943) circle (0.1pt);
\fill[blue] (0.5646392290876473,0.35459424606590373) circle (0.1pt);\draw(0.5646392290876473,0.35459424606590373) circle (0.1pt);
\fill[blue] (0.5570056583155196,0.3675261222104846) circle (0.1pt);\draw(0.5570056583155196,0.3675261222104846) circle (0.1pt);
\fill[blue] (0.5326754741324017,0.40874338970209434) circle (0.1pt);\draw(0.5326754741324017,0.40874338970209434) circle (0.1pt);
\fill[blue] (0.3849643126762702,0.3563026374685184) circle (0.1pt);\draw(0.3849643126762702,0.3563026374685184) circle (0.1pt);
\fill[blue] (0.48139088714400735,0.22947900357015216) circle (0.1pt);\draw(0.48139088714400735,0.22947900357015216) circle (0.1pt);
\fill[blue] (0.4952009553583087,0.24437837632781984) circle (0.1pt);\draw(0.4952009553583087,0.24437837632781984) circle (0.1pt);
\fill[blue] (0.5464741713477421,0.2996958987076048) circle (0.1pt);\draw(0.5464741713477421,0.2996958987076048) circle (0.1pt);
\fill[blue] (0.45162040934298203,0.19736030299255442) circle (0.1pt);\draw(0.45162040934298203,0.19736030299255442) circle (0.1pt);
\fill[blue] (0.44500451138548675,0.17105528445839943) circle (0.1pt);\draw(0.44500451138548675,0.17105528445839943) circle (0.1pt);
\fill[blue] (0.45877901890888895,0.13346796697448662) circle (0.1pt);\draw(0.45877901890888895,0.13346796697448662) circle (0.1pt);
\fill[blue] (0.7221119059162853,0.4758477701896109) circle (0.1pt);\draw(0.7221119059162853,0.4758477701896109) circle (0.1pt);
\fill[blue] (0.8272393746070328,0.4730369662060027) circle (0.1pt);\draw(0.8272393746070328,0.4730369662060027) circle (0.1pt);
\fill[blue] (0.7711711320631728,0.41318386283041797) circle (0.1pt);\draw(0.7711711320631728,0.41318386283041797) circle (0.1pt);
\fill[blue] (0.5800834392674282,0.554694787525681) circle (0.1pt);\draw(0.5800834392674282,0.554694787525681) circle (0.1pt);
\fill[blue] (0.6188652259059673,0.5462977123528099) circle (0.1pt);\draw(0.6188652259059673,0.5462977123528099) circle (0.1pt);
\fill[blue] (0.6361928686502474,0.5410765590383083) circle (0.1pt);\draw(0.6361928686502474,0.5410765590383083) circle (0.1pt);
\fill[blue] (0.7880511358690263,0.6120556809230696) circle (0.1pt);\draw(0.7880511358690263,0.6120556809230696) circle (0.1pt);
\fill[blue] (0.7843733421016718,0.6092230402432683) circle (0.1pt);\draw(0.7843733421016718,0.6092230402432683) circle (0.1pt);
\fill[blue] (0.7940377748297177,0.6166665966317921) circle (0.1pt);\draw(0.7940377748297177,0.6166665966317921) circle (0.1pt);
\fill[blue] (0.7851149531728918,0.6097942298808102) circle (0.1pt);\draw(0.7851149531728918,0.6097942298808102) circle (0.1pt);
\fill[blue] (0.8108312686319776,0.5941237704915503) circle (0.1pt);\draw(0.8108312686319776,0.5941237704915503) circle (0.1pt);
\fill[blue] (0.8118874384593696,0.1633564836689559) circle (0.1pt);\draw(0.8118874384593696,0.1633564836689559) circle (0.1pt);
\fill[blue] (0.8851632825205141,0.5493531197369946) circle (0.1pt);\draw(0.8851632825205141,0.5493531197369946) circle (0.1pt);
\fill[blue] (1.0241664694608577,0.57153475091379) circle (0.1pt);\draw(1.0241664694608577,0.57153475091379) circle (0.1pt);
\fill[blue] (1.0378006206639054,0.5737104399756852) circle (0.1pt);\draw(1.0378006206639054,0.5737104399756852) circle (0.1pt);
\fill[blue] (1.0823433206558173,0.580818404653238) circle (0.1pt);\draw(1.0823433206558173,0.580818404653238) circle (0.1pt);
\fill[blue] (0.9421973439307684,0.5584544111718442) circle (0.1pt);\draw(0.9421973439307684,0.5584544111718442) circle (0.1pt);
\fill[blue] (1.0605948498489768,0.5980528201790608) circle (0.1pt);\draw(1.0605948498489768,0.5980528201790608) circle (0.1pt);
\fill[blue] (0.9625363197915181,0.6668209632607429) circle (0.1pt);\draw(0.9625363197915181,0.6668209632607429) circle (0.1pt);
\fill[blue] (1.0418135837704467,0.6112240638453353) circle (0.1pt);\draw(1.0418135837704467,0.6112240638453353) circle (0.1pt);
\fill[blue] (0.9168260577530807,0.6988774282956411) circle (0.1pt);\draw(0.9168260577530807,0.6988774282956411) circle (0.1pt);
\fill[blue] (1.0588598410834567,0.5992695764466222) circle (0.1pt);\draw(1.0588598410834567,0.5992695764466222) circle (0.1pt);
\fill[blue] (0.9338936187195976,0.6869080007107792) circle (0.1pt);\draw(0.9338936187195976,0.6869080007107792) circle (0.1pt);
\fill[blue] (0.84050884537445,0.7523984527019172) circle (0.1pt);\draw(0.84050884537445,0.7523984527019172) circle (0.1pt);
\fill[blue] (0.642114287806097,0.6304560849618224) circle (0.1pt);\draw(0.642114287806097,0.6304560849618224) circle (0.1pt);
\fill[blue] (0.7630519194473339,0.6944316881617858) circle (0.1pt);\draw(0.7630519194473339,0.6944316881617858) circle (0.1pt);
\fill[blue] (0.7580659591276837,0.6931913911438812) circle (0.1pt);\draw(0.7580659591276837,0.6931913911438812) circle (0.1pt);
\fill[blue] (0.7097088513509334,0.6811621785700015) circle (0.1pt);\draw(0.7097088513509334,0.6811621785700015) circle (0.1pt);
\fill[blue] (0.6984784752351365,0.6783685337927906) circle (0.1pt);\draw(0.6984784752351365,0.6783685337927906) circle (0.1pt);
\fill[blue] (0.714648332376192,0.6171165250868014) circle (0.1pt);\draw(0.714648332376192,0.6171165250868014) circle (0.1pt);
\fill[blue] (0.7860546895085898,0.6214997109685254) circle (0.1pt);\draw(0.7860546895085898,0.6214997109685254) circle (0.1pt);
\fill[blue] (0.8032960357300718,0.6954122436021932) circle (0.1pt);\draw(0.8032960357300718,0.6954122436021932) circle (0.1pt);
\fill[blue] (0.6595268144771232,0.6905536626697233) circle (0.1pt);\draw(0.6595268144771232,0.6905536626697233) circle (0.1pt);
\fill[blue] (0.6273026388087081,0.6064132747783626) circle (0.1pt);\draw(0.6273026388087081,0.6064132747783626) circle (0.1pt);
 \end{scope} 
 \draw (0.2,0.2) rectangle (0.8,0.8); 
 \end{tikzpicture}

%% file: VoronoiSINR2.tex
\begin{tikzpicture}[scale=8] 
 \begin{scope} 
\clip(0.2,0.2) rectangle (0.8,0.8);
\draw[red, thick] (0.8340383503796273,0.7758938836779715)--(0.822841928213028,0.8123889598178787);
\draw[red, thick] (0.16402452764015735,0.3738974712746203)--(0.17418575103719508,0.3127497753447539);
\draw[red, thick] (0.148397993788087,0.28793013543623097)--(0.17418575103719508,0.3127497753447539);
\draw[red, thick] (0.16402452764015735,0.3738974712746203)--(0.18588103912022907,0.40276003273976374);
\draw[red, thick] (0.18588103912022907,0.40276003273976374)--(0.13638721173493218,0.6043853598410402);
\draw[red, thick] (0.0072051335010195755,0.8663122452666777)--(0.07225031953352958,0.7988984371260657);
\draw[red, thick] (0.06964041109724696,0.7722567622110954)--(0.07225031953352958,0.7988984371260657);
\draw[red, thick] (0.822841928213028,0.8123889598178787)--(0.767274141410929,0.8903491846659611);
\draw[red, thick] (0.767274141410929,0.8903491846659611)--(0.6347894819842471,0.8155039577188894);
\draw[red, thick] (0.6347894819842471,0.8155039577188894)--(0.5646470338381138,0.8144491937677845);
\draw[red, thick] (0.148397993788087,0.28793013543623097)--(0.13656917844354743,0.2001593013620731);
\draw[red, thick] (0.19799452036188725,0.2983935547149186)--(0.17418575103719508,0.3127497753447539);
\draw[red, thick] (0.19799452036188725,0.2983935547149186)--(0.39386035741641695,0.3879950389869909);
\draw[red, thick] (0.2614932165060834,0.4483621656048845)--(0.25909724262654277,0.44854516125233673);
\draw[red, thick] (0.2614932165060834,0.4483621656048845)--(0.4157501376828479,0.4287256297513061);
\draw[red, thick] (0.18588103912022907,0.40276003273976374)--(0.25909724262654277,0.44854516125233673);
\draw[red, thick] (0.4157501376828479,0.4287256297513061)--(0.39386035741641695,0.3879950389869909);
\draw[red, thick] (0.2125410814435029,0.932206125956973)--(0.1503709212537891,0.9000494683595954);
\draw[red, thick] (0.2125410814435029,0.932206125956973)--(0.049259228651907655,0.9730210070877323);
\draw[red, thick] (0.1503709212537891,0.9000494683595954)--(0.07950484344558738,0.9060196908690696);
\draw[red, thick] (0.07950484344558738,0.9060196908690696)--(0.03418643383988759,0.9411415638172064);
\draw[red, thick] (0.049259228651907655,0.9730210070877323)--(0.03418643383988759,0.9411415638172064);
\draw[red, thick] (0.3078280725794059,0.9664502126130512)--(0.2125410814435029,0.932206125956973);
\draw[red, thick] (0.0072051335010195755,0.8663122452666777)--(0.07950484344558738,0.9060196908690696);
\draw[red, thick] (0.07225031953352958,0.7988984371260657)--(0.1503709212537891,0.9000494683595954);
\draw[red, thick] (0.20785657217564413,0.7391137029871157)--(0.299736887125702,0.8166411601178638);
\draw[red, thick] (0.20785657217564413,0.7391137029871157)--(0.1626026765847196,0.6149341693854876);
\draw[red, thick] (0.299736887125702,0.8166411601178638)--(0.410824171102973,0.6379558481429204);
\draw[red, thick] (0.1626026765847196,0.6149341693854876)--(0.21771196590924116,0.5929681031496008);
\draw[red, thick] (0.21771196590924116,0.5929681031496008)--(0.3581748076390561,0.6189381158954018);
\draw[red, thick] (0.3581748076390561,0.6189381158954018)--(0.37128451309713617,0.6232642739680077);
\draw[red, thick] (0.410824171102973,0.6379558481429204)--(0.37128451309713617,0.6232642739680077);
\draw[red, thick] (0.3078280725794059,0.9664502126130512)--(0.299736887125702,0.8166411601178638);
\draw[red, thick] (0.06964041109724696,0.7722567622110954)--(0.20785657217564413,0.7391137029871157);
\draw[red, thick] (0.13638721173493218,0.6043853598410402)--(0.1626026765847196,0.6149341693854876);
\draw[red, thick] (0.5646470338381138,0.8144491937677845)--(0.4573568415412147,0.6126174690235537);
\draw[red, thick] (0.4573568415412147,0.6126174690235537)--(0.410824171102973,0.6379558481429204);
\draw[red, thick] (0.25909724262654277,0.44854516125233673)--(0.2553892027923351,0.48539597554160624);
\draw[red, thick] (0.2553892027923351,0.48539597554160624)--(0.21771196590924116,0.5929681031496008);
\draw[red, thick] (0.2553892027923351,0.48539597554160624)--(0.3581748076390561,0.6189381158954018);
\draw[red, thick] (0.2614932165060834,0.4483621656048845)--(0.37128451309713617,0.6232642739680077);
\draw[red, thick] (0.5092587354725757,0.5327225311464738)--(0.4573568415412147,0.6126174690235537);
\draw[red, thick] (0.5092587354725757,0.5327225311464738)--(0.4959932207799768,0.47088604498338843);
\draw[red, thick] (0.4959932207799768,0.47088604498338843)--(0.4157501376828479,0.4287256297513061);
\draw[red, thick] (0.45995291048752523,0.028965958412618953)--(0.33570477350580863,0.0944446314103765);
\draw[red, thick] (0.13656917844354743,0.2001593013620731)--(0.22603564326488124,0.11592720692047373);
\draw[red, thick] (0.679280231333423,0.07685494809943949)--(0.7039360788752047,0.13691025596481718);
\draw[red, thick] (0.679280231333423,0.07685494809943949)--(0.5401162699392226,0.10592861733374764);
\draw[red, thick] (0.7039360788752047,0.13691025596481718)--(0.692984094105217,0.16657177838006781);
\draw[red, thick] (0.692984094105217,0.16657177838006781)--(0.5705030851078778,0.14956730017396283);
\draw[red, thick] (0.5705030851078778,0.14956730017396283)--(0.5401162699392226,0.10592861733374764);
\draw[red, thick] (0.6704261078411188,0.2983594041110991)--(0.692984094105217,0.16657177838006781);
\draw[red, thick] (0.6704261078411188,0.2983594041110991)--(0.609626416494816,0.3156521591180685);
\draw[red, thick] (0.5705030851078778,0.14956730017396283)--(0.609626416494816,0.3156521591180685);
\draw[red, thick] (0.45995291048752523,0.028965958412618953)--(0.4765623254864447,0.08494160121449679);
\draw[red, thick] (0.4765623254864447,0.08494160121449679)--(0.5401162699392226,0.10592861733374764);
\draw[red, thick] (0.38150587398438834,0.343981855497234)--(0.4497840091409947,0.1953790520155802);
\draw[red, thick] (0.38150587398438834,0.343981855497234)--(0.29113000848283804,0.21354918994579553);
\draw[red, thick] (0.297180677272192,0.20845266990432376)--(0.2922904993448189,0.21081962561616713);
\draw[red, thick] (0.297180677272192,0.20845266990432376)--(0.4415822666940058,0.18039376697821286);
\draw[red, thick] (0.4415822666940058,0.18039376697821286)--(0.4497840091409947,0.1953790520155802);
\draw[red, thick] (0.29113000848283804,0.21354918994579553)--(0.2922904993448189,0.21081962561616713);
\draw[red, thick] (0.4959932207799768,0.47088604498338843)--(0.5773694923567108,0.33302816738559193);
\draw[red, thick] (0.39386035741641695,0.3879950389869909)--(0.38150587398438834,0.343981855497234);
\draw[red, thick] (0.5773694923567108,0.33302816738559193)--(0.4497840091409947,0.1953790520155802);
\draw[red, thick] (0.19799452036188725,0.2983935547149186)--(0.29113000848283804,0.21354918994579553);
\draw[red, thick] (0.33570477350580863,0.0944446314103765)--(0.297180677272192,0.20845266990432376);
\draw[red, thick] (0.22603564326488124,0.11592720692047373)--(0.2922904993448189,0.21081962561616713);
\draw[red, thick] (0.4765623254864447,0.08494160121449679)--(0.4415822666940058,0.18039376697821286);
\draw[red, thick] (0.609626416494816,0.3156521591180685)--(0.5773694923567108,0.33302816738559193);
\draw[red, thick] (0.6966078171943033,0.517084675521012)--(0.8394170635303543,0.5386328082928883);
\draw[red, thick] (0.6966078171943033,0.517084675521012)--(0.7629987978835259,0.4097388057906942);
\draw[red, thick] (0.8394170635303543,0.5386328082928883)--(0.8199448304662863,0.43374447248530434);
\draw[red, thick] (0.7629987978835259,0.4097388057906942)--(0.8199448304662863,0.43374447248530434);
\draw[red, thick] (0.5092587354725757,0.5327225311464738)--(0.5854609328882098,0.5563630707966651);
\draw[red, thick] (0.6704261078411188,0.2983594041110991)--(0.706255398851686,0.31142037155747987);
\draw[red, thick] (0.706255398851686,0.31142037155747987)--(0.7629987978835259,0.4097388057906942);
\draw[red, thick] (0.6966078171943033,0.517084675521012)--(0.6791059840045244,0.528146008831672);
\draw[red, thick] (0.5854609328882098,0.5563630707966651)--(0.6791059840045244,0.528146008831672);
\draw[red, thick] (0.8394170635303543,0.5386328082928883)--(0.8421316549547286,0.5424862865901321);
\draw[red, thick] (0.6791059840045244,0.528146008831672)--(0.7961709115289796,0.6183095407912882);
\draw[red, thick] (0.8421316549547286,0.5424862865901321)--(0.7961709115289796,0.6183095407912882);
\draw[red, thick] (0.8180743974620243,0.1360181565046783)--(0.7847271846057482,0.2833695475015336);
\draw[red, thick] (0.7039360788752047,0.13691025596481718)--(0.8180743974620243,0.1360181565046783);
\draw[red, thick] (0.706255398851686,0.31142037155747987)--(0.7847271846057482,0.2833695475015336);
\draw[red, thick] (0.8340383503796273,0.7758938836779715)--(0.8378339831936373,0.7542743252554434);
\draw[red, thick] (0.663654233906528,0.6697057286556454)--(0.6320097831139766,0.6120438657752406);
\draw[red, thick] (0.663654233906528,0.6697057286556454)--(0.8042869452731249,0.7046892265869095);
\draw[red, thick] (0.6320097831139766,0.6120438657752406)--(0.7954628450760505,0.6220772182985786);
\draw[red, thick] (0.8042869452731249,0.7046892265869095)--(0.7954628450760505,0.6220772182985786);
\draw[red, thick] (0.6347894819842471,0.8155039577188894)--(0.663654233906528,0.6697057286556454);
\draw[red, thick] (0.5854609328882098,0.5563630707966651)--(0.6320097831139766,0.6120438657752406);
\draw[red, thick] (0.8378339831936373,0.7542743252554434)--(0.8042869452731249,0.7046892265869095);
\draw[red, thick] (0.7961709115289796,0.6183095407912882)--(0.7954628450760505,0.6220772182985786);
\draw[black, very thick] (0.05789114625497006,0.7591447714783306)--(0.07059719372583992,0.7820235001551941);
\draw[black, very thick] (0.05789114625497006,0.7591447714783306)--(0.025416319820836576,0.7229033871050434);
\draw[black, very thick] (0.816871669547315,0.8207650844535032)--(0.797244196123064,0.8483019422161482);
\draw[black, very thick] (0.5150163133378375,0.8659430761214151)--(0.47164836244225017,0.9109390813944592);
\draw[black, very thick] (0.29224699803261267,0.34151062925235404)--(0.25616296056720894,0.325003496798576);
\draw[black, very thick] (0.3693375759616444,0.43463383746958717)--(0.37181151999371226,0.4343189103317002);
\draw[black, very thick] (0.4100600707523283,0.4294499616155625)--(0.4169904242473743,0.42937728712224627);
\draw[black, very thick] (0.25532816389877266,0.4461882000106892)--(0.2585381570265892,0.4541014013513852);
\draw[black, very thick] (0.15937905866236854,0.9047088031490487)--(0.1590067661594326,0.9455878746166037);
\draw[black, very thick] (0.15937905866236854,0.9047088031490487)--(0.17361789056520271,0.9419355930212361);
\draw[black, very thick] (0.1590067661594326,0.9455878746166037)--(0.17361789056520271,0.9419355930212361);
\draw[black, very thick] (0.2921345429365123,0.9608102969999692)--(0.30700725570527704,0.9512527117166875);
\draw[black, very thick] (0.2526856156454257,0.7769398874936125)--(0.27892638396565794,0.7990815186159101);
\draw[black, very thick] (0.27892638396565794,0.7990815186159101)--(0.29981130217595947,0.8180189617063004);
\draw[black, very thick] (0.19387030084463816,0.7007344954508157)--(0.18563030063821997,0.7444433737589387);
\draw[black, very thick] (0.37293603650876134,0.6988993954844914)--(0.3982187058250188,0.6582318996866904);
\draw[black, very thick] (0.18134071416677522,0.6074653564194161)--(0.1726754312359303,0.6109192601286578);
\draw[black, very thick] (0.34509732830012835,0.6165202358288133)--(0.36683993239497326,0.6217975775006084);
\draw[black, very thick] (0.232149332632162,0.5956374112256204)--(0.2277040583383429,0.5644397202114);
\draw[black, very thick] (0.3010446550260852,0.8408546046323935)--(0.30055791219730454,0.8318425158098488);
\draw[black, very thick] (0.30700725570527704,0.9512527117166875)--(0.304643617111923,0.9074897229510144);
\draw[black, very thick] (0.29981130217595947,0.8180189617063004)--(0.30055791219730454,0.8318425158098488);
\draw[black, very thick] (0.23156648132676858,0.5534121318594152)--(0.2277040583383429,0.5644397202114);
\draw[black, very thick] (0.26648526291374264,0.4998123096218199)--(0.2919655586376268,0.49690587581737233);
\draw[black, very thick] (0.293302842223799,0.5346545129197736)--(0.2919655586376268,0.49690587581737233);
\draw[black, very thick] (0.3498485323164778,0.5891158637590141)--(0.34045815137782165,0.574156595807007);
\draw[black, very thick] (0.38355244884048856,0.06922894316707184)--(0.37508174221498614,0.07369299906895496);
\draw[black, very thick] (0.44145691925859737,0.03871333161972445)--(0.4744519247018341,0.0778293073440714);
\draw[black, very thick] (0.37508174221498614,0.07369299906895496)--(0.35676393151346963,0.08334647127462609);
\draw[black, very thick] (0.35676393151346963,0.08334647127462609)--(0.33402090497369247,0.09314723041436684);
\draw[black, very thick] (0.27052606386950184,0.04422519848677735)--(0.25415869288763704,0.031614330708849434);
\draw[black, very thick] (0.19012226347068706,-0.017724988674975257)--(0.17858300600531926,-0.004973681536212152);
\draw[black, very thick] (0.15362496192050873,0.1841013922061019)--(0.153557661738525,0.18416475488949516);
\draw[black, very thick] (0.15362496192050873,0.1841013922061019)--(0.14846174817708205,0.18896252434680957);
\draw[black, very thick] (0.153557661738525,0.18416475488949516)--(0.14846174817708205,0.18896252434680957);
\draw[black, very thick] (0.1872786858015116,0.01718136555514127)--(0.18417778144474048,0.00928081221598432);
\draw[black, very thick] (0.18417778144474048,0.00928081221598432)--(0.17858300600531926,-0.004973681536212152);
\draw[black, very thick] (0.4744519247018341,0.0778293073440714)--(0.5263999923459467,0.10139917795079811);
\draw[black, very thick] (0.44625527759791844,0.20305910066999244)--(0.45162040934298203,0.19736030299255442);
\draw[black, very thick] (0.35675114120119855,0.3082552180370413)--(0.37932348476575223,0.34083217785400804);
\draw[black, very thick] (0.37932348476575223,0.34083217785400804)--(0.3849643126762702,0.3563026374685184);
\draw[black, very thick] (0.43025304052789864,0.18259516696622943)--(0.44500451138548675,0.17105528445839943);
\draw[black, very thick] (0.5646392290876473,0.35459424606590373)--(0.5570056583155196,0.3675261222104846);
\draw[black, very thick] (0.5570056583155196,0.3675261222104846)--(0.5326754741324017,0.40874338970209434);
\draw[black, very thick] (0.48139088714400735,0.22947900357015216)--(0.4952009553583087,0.24437837632781984);
\draw[black, very thick] (0.44500451138548675,0.17105528445839943)--(0.45877901890888895,0.13346796697448662);
\draw[black, very thick] (0.5800834392674282,0.554694787525681)--(0.6188652259059673,0.5462977123528099);
\draw[black, very thick] (0.6188652259059673,0.5462977123528099)--(0.6361928686502474,0.5410765590383083);
\draw[black, very thick] (0.7880511358690263,0.6120556809230696)--(0.7851149531728918,0.6097942298808102);
\draw[black, very thick] (0.7880511358690263,0.6120556809230696)--(0.7940377748297177,0.6166665966317921);
\draw[black, very thick] (0.7843733421016718,0.6092230402432683)--(0.7851149531728918,0.6097942298808102);
\draw[black, very thick] (0.7940377748297177,0.6166665966317921)--(0.7860546895085898,0.6214997109685254);
\draw[black, very thick] (0.8851632825205141,0.5493531197369946)--(0.9421973439307684,0.5584544111718442);
\draw[black, very thick] (1.0241664694608577,0.57153475091379)--(1.0378006206639054,0.5737104399756852);
\draw[black, very thick] (1.0823433206558173,0.580818404653238)--(1.0605948498489768,0.5980528201790608);
\draw[black, very thick] (1.0823433206558173,0.580818404653238)--(1.0588598410834567,0.5992695764466222);
\draw[black, very thick] (1.0605948498489768,0.5980528201790608)--(1.0588598410834567,0.5992695764466222);
\draw[black, very thick] (1.0605948498489768,0.5980528201790608)--(1.0418135837704467,0.6112240638453353);
\draw[black, very thick] (0.9625363197915181,0.6668209632607429)--(0.9338936187195976,0.6869080007107792);
\draw[black, very thick] (1.0418135837704467,0.6112240638453353)--(1.0588598410834567,0.5992695764466222);
\draw[black, very thick] (0.9168260577530807,0.6988774282956411)--(0.9338936187195976,0.6869080007107792);
\draw[black, very thick] (0.642114287806097,0.6304560849618224)--(0.6273026388087081,0.6064132747783626);
\draw[black, very thick] (0.7630519194473339,0.6944316881617858)--(0.7580659591276837,0.6931913911438812);
\draw[black, very thick] (0.7097088513509334,0.6811621785700015)--(0.6984784752351365,0.6783685337927906);
\draw[black, very thick] (0.6984784752351365,0.6783685337927906)--(0.6595268144771232,0.6905536626697233);
\fill[blue] (0.1813061501361128,0.42139697266953957) circle (0.1pt);\draw(0.1813061501361128,0.42139697266953957) circle (0.1pt);
\fill[blue] (0.14872563829292249,0.5541217325950919) circle (0.1pt);\draw(0.14872563829292249,0.5541217325950919) circle (0.1pt);
\fill[blue] (0.05789114625497006,0.7591447714783306) circle (0.1pt);\draw(0.05789114625497006,0.7591447714783306) circle (0.1pt);
\fill[blue] (0.025416319820836576,0.7229033871050434) circle (0.1pt);\draw(0.025416319820836576,0.7229033871050434) circle (0.1pt);
\fill[blue] (0.07059719372583992,0.7820235001551941) circle (0.1pt);\draw(0.07059719372583992,0.7820235001551941) circle (0.1pt);
\fill[blue] (0.816871669547315,0.8207650844535032) circle (0.1pt);\draw(0.816871669547315,0.8207650844535032) circle (0.1pt);
\fill[blue] (0.797244196123064,0.8483019422161482) circle (0.1pt);\draw(0.797244196123064,0.8483019422161482) circle (0.1pt);
\fill[blue] (0.5150163133378375,0.8659430761214151) circle (0.1pt);\draw(0.5150163133378375,0.8659430761214151) circle (0.1pt);
\fill[blue] (0.47164836244225017,0.9109390813944592) circle (0.1pt);\draw(0.47164836244225017,0.9109390813944592) circle (0.1pt);
\fill[blue] (0.38653207447931925,0.9992506770062238) circle (0.1pt);\draw(0.38653207447931925,0.9992506770062238) circle (0.1pt);
\fill[blue] (0.7024556644193323,0.853730955492912) circle (0.1pt);\draw(0.7024556644193323,0.853730955492912) circle (0.1pt);
\fill[blue] (0.29224699803261267,0.34151062925235404) circle (0.1pt);\draw(0.29224699803261267,0.34151062925235404) circle (0.1pt);
\fill[blue] (0.25616296056720894,0.325003496798576) circle (0.1pt);\draw(0.25616296056720894,0.325003496798576) circle (0.1pt);
\fill[blue] (0.3693375759616444,0.43463383746958717) circle (0.1pt);\draw(0.3693375759616444,0.43463383746958717) circle (0.1pt);
\fill[blue] (0.37181151999371226,0.4343189103317002) circle (0.1pt);\draw(0.37181151999371226,0.4343189103317002) circle (0.1pt);
\fill[blue] (0.4100600707523283,0.4294499616155625) circle (0.1pt);\draw(0.4100600707523283,0.4294499616155625) circle (0.1pt);
\fill[blue] (0.25532816389877266,0.4461882000106892) circle (0.1pt);\draw(0.25532816389877266,0.4461882000106892) circle (0.1pt);
\fill[blue] (0.15937905866236854,0.9047088031490487) circle (0.1pt);\draw(0.15937905866236854,0.9047088031490487) circle (0.1pt);
\fill[blue] (0.1590067661594326,0.9455878746166037) circle (0.1pt);\draw(0.1590067661594326,0.9455878746166037) circle (0.1pt);
\fill[blue] (0.17361789056520271,0.9419355930212361) circle (0.1pt);\draw(0.17361789056520271,0.9419355930212361) circle (0.1pt);
\fill[blue] (0.06284152846679965,0.9189337987707326) circle (0.1pt);\draw(0.06284152846679965,0.9189337987707326) circle (0.1pt);
\fill[blue] (0.2921345429365123,0.9608102969999692) circle (0.1pt);\draw(0.2921345429365123,0.9608102969999692) circle (0.1pt);
\fill[blue] (0.2526856156454257,0.7769398874936125) circle (0.1pt);\draw(0.2526856156454257,0.7769398874936125) circle (0.1pt);
\fill[blue] (0.27892638396565794,0.7990815186159101) circle (0.1pt);\draw(0.27892638396565794,0.7990815186159101) circle (0.1pt);
\fill[blue] (0.19387030084463816,0.7007344954508157) circle (0.1pt);\draw(0.19387030084463816,0.7007344954508157) circle (0.1pt);
\fill[blue] (0.37293603650876134,0.6988993954844914) circle (0.1pt);\draw(0.37293603650876134,0.6988993954844914) circle (0.1pt);
\fill[blue] (0.3982187058250188,0.6582318996866904) circle (0.1pt);\draw(0.3982187058250188,0.6582318996866904) circle (0.1pt);
\fill[blue] (0.18134071416677522,0.6074653564194161) circle (0.1pt);\draw(0.18134071416677522,0.6074653564194161) circle (0.1pt);
\fill[blue] (0.1726754312359303,0.6109192601286578) circle (0.1pt);\draw(0.1726754312359303,0.6109192601286578) circle (0.1pt);
\fill[blue] (0.34509732830012835,0.6165202358288133) circle (0.1pt);\draw(0.34509732830012835,0.6165202358288133) circle (0.1pt);
\fill[blue] (0.232149332632162,0.5956374112256204) circle (0.1pt);\draw(0.232149332632162,0.5956374112256204) circle (0.1pt);
\fill[blue] (0.36683993239497326,0.6217975775006084) circle (0.1pt);\draw(0.36683993239497326,0.6217975775006084) circle (0.1pt);
\fill[blue] (0.3010446550260852,0.8408546046323935) circle (0.1pt);\draw(0.3010446550260852,0.8408546046323935) circle (0.1pt);
\fill[blue] (0.30700725570527704,0.9512527117166875) circle (0.1pt);\draw(0.30700725570527704,0.9512527117166875) circle (0.1pt);
\fill[blue] (0.304643617111923,0.9074897229510144) circle (0.1pt);\draw(0.304643617111923,0.9074897229510144) circle (0.1pt);
\fill[blue] (0.29981130217595947,0.8180189617063004) circle (0.1pt);\draw(0.29981130217595947,0.8180189617063004) circle (0.1pt);
\fill[blue] (0.30055791219730454,0.8318425158098488) circle (0.1pt);\draw(0.30055791219730454,0.8318425158098488) circle (0.1pt);
\fill[blue] (0.18563030063821997,0.7444433737589387) circle (0.1pt);\draw(0.18563030063821997,0.7444433737589387) circle (0.1pt);
\fill[blue] (0.5337666871494476,0.7563578295951302) circle (0.1pt);\draw(0.5337666871494476,0.7563578295951302) circle (0.1pt);
\fill[blue] (0.4828113417610102,0.6605018615470278) circle (0.1pt);\draw(0.4828113417610102,0.6605018615470278) circle (0.1pt);
\fill[blue] (0.2585381570265892,0.4541014013513852) circle (0.1pt);\draw(0.2585381570265892,0.4541014013513852) circle (0.1pt);
\fill[blue] (0.23156648132676858,0.5534121318594152) circle (0.1pt);\draw(0.23156648132676858,0.5534121318594152) circle (0.1pt);
\fill[blue] (0.2277040583383429,0.5644397202114) circle (0.1pt);\draw(0.2277040583383429,0.5644397202114) circle (0.1pt);
\fill[blue] (0.26648526291374264,0.4998123096218199) circle (0.1pt);\draw(0.26648526291374264,0.4998123096218199) circle (0.1pt);
\fill[blue] (0.293302842223799,0.5346545129197736) circle (0.1pt);\draw(0.293302842223799,0.5346545129197736) circle (0.1pt);
\fill[blue] (0.3498485323164778,0.5891158637590141) circle (0.1pt);\draw(0.3498485323164778,0.5891158637590141) circle (0.1pt);
\fill[blue] (0.2919655586376268,0.49690587581737233) circle (0.1pt);\draw(0.2919655586376268,0.49690587581737233) circle (0.1pt);
\fill[blue] (0.34045815137782165,0.574156595807007) circle (0.1pt);\draw(0.34045815137782165,0.574156595807007) circle (0.1pt);
\fill[blue] (0.5084688879484911,0.5290406992216823) circle (0.1pt);\draw(0.5084688879484911,0.5290406992216823) circle (0.1pt);
\fill[blue] (0.4169904242473743,0.42937728712224627) circle (0.1pt);\draw(0.4169904242473743,0.42937728712224627) circle (0.1pt);
\fill[blue] (0.49024128124118205,0.4678639258176544) circle (0.1pt);\draw(0.49024128124118205,0.4678639258176544) circle (0.1pt);
\fill[blue] (0.38355244884048856,0.06922894316707184) circle (0.1pt);\draw(0.38355244884048856,0.06922894316707184) circle (0.1pt);
\fill[blue] (0.44145691925859737,0.03871333161972445) circle (0.1pt);\draw(0.44145691925859737,0.03871333161972445) circle (0.1pt);
\fill[blue] (0.37508174221498614,0.07369299906895496) circle (0.1pt);\draw(0.37508174221498614,0.07369299906895496) circle (0.1pt);
\fill[blue] (0.35676393151346963,0.08334647127462609) circle (0.1pt);\draw(0.35676393151346963,0.08334647127462609) circle (0.1pt);
\fill[blue] (0.27052606386950184,0.04422519848677735) circle (0.1pt);\draw(0.27052606386950184,0.04422519848677735) circle (0.1pt);
\fill[blue] (0.33402090497369247,0.09314723041436684) circle (0.1pt);\draw(0.33402090497369247,0.09314723041436684) circle (0.1pt);
\fill[blue] (0.19012226347068706,-0.017724988674975257) circle (0.1pt);\draw(0.19012226347068706,-0.017724988674975257) circle (0.1pt);
\fill[blue] (0.25415869288763704,0.031614330708849434) circle (0.1pt);\draw(0.25415869288763704,0.031614330708849434) circle (0.1pt);
\fill[blue] (0.19976901885816142,0.14065706287384386) circle (0.1pt);\draw(0.19976901885816142,0.14065706287384386) circle (0.1pt);
\fill[blue] (0.15362496192050873,0.1841013922061019) circle (0.1pt);\draw(0.15362496192050873,0.1841013922061019) circle (0.1pt);
\fill[blue] (0.153557661738525,0.18416475488949516) circle (0.1pt);\draw(0.153557661738525,0.18416475488949516) circle (0.1pt);
\fill[blue] (0.14846174817708205,0.18896252434680957) circle (0.1pt);\draw(0.14846174817708205,0.18896252434680957) circle (0.1pt);
\fill[blue] (0.1872786858015116,0.01718136555514127) circle (0.1pt);\draw(0.1872786858015116,0.01718136555514127) circle (0.1pt);
\fill[blue] (0.18417778144474048,0.00928081221598432) circle (0.1pt);\draw(0.18417778144474048,0.00928081221598432) circle (0.1pt);
\fill[blue] (0.17858300600531926,-0.004973681536212152) circle (0.1pt);\draw(0.17858300600531926,-0.004973681536212152) circle (0.1pt);
\fill[blue] (0.694675762855622,0.11435450611503778) circle (0.1pt);\draw(0.694675762855622,0.11435450611503778) circle (0.1pt);
\fill[blue] (0.6796560855289289,0.2444362909933815) circle (0.1pt);\draw(0.6796560855289289,0.2444362909933815) circle (0.1pt);
\fill[blue] (0.4744519247018341,0.0778293073440714) circle (0.1pt);\draw(0.4744519247018341,0.0778293073440714) circle (0.1pt);
\fill[blue] (0.5263999923459467,0.10139917795079811) circle (0.1pt);\draw(0.5263999923459467,0.10139917795079811) circle (0.1pt);
\fill[blue] (0.44625527759791844,0.20305910066999244) circle (0.1pt);\draw(0.44625527759791844,0.20305910066999244) circle (0.1pt);
\fill[blue] (0.35675114120119855,0.3082552180370413) circle (0.1pt);\draw(0.35675114120119855,0.3082552180370413) circle (0.1pt);
\fill[blue] (0.37932348476575223,0.34083217785400804) circle (0.1pt);\draw(0.37932348476575223,0.34083217785400804) circle (0.1pt);
\fill[blue] (0.3781940958333602,0.1927108234591293) circle (0.1pt);\draw(0.3781940958333602,0.1927108234591293) circle (0.1pt);
\fill[blue] (0.43025304052789864,0.18259516696622943) circle (0.1pt);\draw(0.43025304052789864,0.18259516696622943) circle (0.1pt);
\fill[blue] (0.5646392290876473,0.35459424606590373) circle (0.1pt);\draw(0.5646392290876473,0.35459424606590373) circle (0.1pt);
\fill[blue] (0.5570056583155196,0.3675261222104846) circle (0.1pt);\draw(0.5570056583155196,0.3675261222104846) circle (0.1pt);
\fill[blue] (0.5326754741324017,0.40874338970209434) circle (0.1pt);\draw(0.5326754741324017,0.40874338970209434) circle (0.1pt);
\fill[blue] (0.3849643126762702,0.3563026374685184) circle (0.1pt);\draw(0.3849643126762702,0.3563026374685184) circle (0.1pt);
\fill[blue] (0.48139088714400735,0.22947900357015216) circle (0.1pt);\draw(0.48139088714400735,0.22947900357015216) circle (0.1pt);
\fill[blue] (0.4952009553583087,0.24437837632781984) circle (0.1pt);\draw(0.4952009553583087,0.24437837632781984) circle (0.1pt);
\fill[blue] (0.5464741713477421,0.2996958987076048) circle (0.1pt);\draw(0.5464741713477421,0.2996958987076048) circle (0.1pt);
\fill[blue] (0.45162040934298203,0.19736030299255442) circle (0.1pt);\draw(0.45162040934298203,0.19736030299255442) circle (0.1pt);
\fill[blue] (0.44500451138548675,0.17105528445839943) circle (0.1pt);\draw(0.44500451138548675,0.17105528445839943) circle (0.1pt);
\fill[blue] (0.45877901890888895,0.13346796697448662) circle (0.1pt);\draw(0.45877901890888895,0.13346796697448662) circle (0.1pt);
\fill[blue] (0.7221119059162853,0.4758477701896109) circle (0.1pt);\draw(0.7221119059162853,0.4758477701896109) circle (0.1pt);
\fill[blue] (0.8272393746070328,0.4730369662060027) circle (0.1pt);\draw(0.8272393746070328,0.4730369662060027) circle (0.1pt);
\fill[blue] (0.7711711320631728,0.41318386283041797) circle (0.1pt);\draw(0.7711711320631728,0.41318386283041797) circle (0.1pt);
\fill[blue] (0.5800834392674282,0.554694787525681) circle (0.1pt);\draw(0.5800834392674282,0.554694787525681) circle (0.1pt);
\fill[blue] (0.6188652259059673,0.5462977123528099) circle (0.1pt);\draw(0.6188652259059673,0.5462977123528099) circle (0.1pt);
\fill[blue] (0.6361928686502474,0.5410765590383083) circle (0.1pt);\draw(0.6361928686502474,0.5410765590383083) circle (0.1pt);
\fill[blue] (0.7880511358690263,0.6120556809230696) circle (0.1pt);\draw(0.7880511358690263,0.6120556809230696) circle (0.1pt);
\fill[blue] (0.7843733421016718,0.6092230402432683) circle (0.1pt);\draw(0.7843733421016718,0.6092230402432683) circle (0.1pt);
\fill[blue] (0.7940377748297177,0.6166665966317921) circle (0.1pt);\draw(0.7940377748297177,0.6166665966317921) circle (0.1pt);
\fill[blue] (0.7851149531728918,0.6097942298808102) circle (0.1pt);\draw(0.7851149531728918,0.6097942298808102) circle (0.1pt);
\fill[blue] (0.8108312686319776,0.5941237704915503) circle (0.1pt);\draw(0.8108312686319776,0.5941237704915503) circle (0.1pt);
\fill[blue] (0.8118874384593696,0.1633564836689559) circle (0.1pt);\draw(0.8118874384593696,0.1633564836689559) circle (0.1pt);
\fill[blue] (0.8851632825205141,0.5493531197369946) circle (0.1pt);\draw(0.8851632825205141,0.5493531197369946) circle (0.1pt);
\fill[blue] (1.0241664694608577,0.57153475091379) circle (0.1pt);\draw(1.0241664694608577,0.57153475091379) circle (0.1pt);
\fill[blue] (1.0378006206639054,0.5737104399756852) circle (0.1pt);\draw(1.0378006206639054,0.5737104399756852) circle (0.1pt);
\fill[blue] (1.0823433206558173,0.580818404653238) circle (0.1pt);\draw(1.0823433206558173,0.580818404653238) circle (0.1pt);
\fill[blue] (0.9421973439307684,0.5584544111718442) circle (0.1pt);\draw(0.9421973439307684,0.5584544111718442) circle (0.1pt);
\fill[blue] (1.0605948498489768,0.5980528201790608) circle (0.1pt);\draw(1.0605948498489768,0.5980528201790608) circle (0.1pt);
\fill[blue] (0.9625363197915181,0.6668209632607429) circle (0.1pt);\draw(0.9625363197915181,0.6668209632607429) circle (0.1pt);
\fill[blue] (1.0418135837704467,0.6112240638453353) circle (0.1pt);\draw(1.0418135837704467,0.6112240638453353) circle (0.1pt);
\fill[blue] (0.9168260577530807,0.6988774282956411) circle (0.1pt);\draw(0.9168260577530807,0.6988774282956411) circle (0.1pt);
\fill[blue] (1.0588598410834567,0.5992695764466222) circle (0.1pt);\draw(1.0588598410834567,0.5992695764466222) circle (0.1pt);
\fill[blue] (0.9338936187195976,0.6869080007107792) circle (0.1pt);\draw(0.9338936187195976,0.6869080007107792) circle (0.1pt);
\fill[blue] (0.84050884537445,0.7523984527019172) circle (0.1pt);\draw(0.84050884537445,0.7523984527019172) circle (0.1pt);
\fill[blue] (0.642114287806097,0.6304560849618224) circle (0.1pt);\draw(0.642114287806097,0.6304560849618224) circle (0.1pt);
\fill[blue] (0.7630519194473339,0.6944316881617858) circle (0.1pt);\draw(0.7630519194473339,0.6944316881617858) circle (0.1pt);
\fill[blue] (0.7580659591276837,0.6931913911438812) circle (0.1pt);\draw(0.7580659591276837,0.6931913911438812) circle (0.1pt);
\fill[blue] (0.7097088513509334,0.6811621785700015) circle (0.1pt);\draw(0.7097088513509334,0.6811621785700015) circle (0.1pt);
\fill[blue] (0.6984784752351365,0.6783685337927906) circle (0.1pt);\draw(0.6984784752351365,0.6783685337927906) circle (0.1pt);
\fill[blue] (0.714648332376192,0.6171165250868014) circle (0.1pt);\draw(0.714648332376192,0.6171165250868014) circle (0.1pt);
\fill[blue] (0.7860546895085898,0.6214997109685254) circle (0.1pt);\draw(0.7860546895085898,0.6214997109685254) circle (0.1pt);
\fill[blue] (0.8032960357300718,0.6954122436021932) circle (0.1pt);\draw(0.8032960357300718,0.6954122436021932) circle (0.1pt);
\fill[blue] (0.6595268144771232,0.6905536626697233) circle (0.1pt);\draw(0.6595268144771232,0.6905536626697233) circle (0.1pt);
\fill[blue] (0.6273026388087081,0.6064132747783626) circle (0.1pt);\draw(0.6273026388087081,0.6064132747783626) circle (0.1pt);
 \end{scope} 
 \draw (0.2,0.2) rectangle (0.8,0.8); 
 \end{tikzpicture}

%% file: VoronoiSINR4.tex
\begin{tikzpicture}[scale=8] 
 \begin{scope} 
\clip(0.2,0.2) rectangle (0.8,0.8);
\draw[red, thick] (0.8340383503796273,0.7758938836779715)--(0.822841928213028,0.8123889598178787);
\draw[red, thick] (0.16402452764015735,0.3738974712746203)--(0.17418575103719508,0.3127497753447539);
\draw[red, thick] (0.148397993788087,0.28793013543623097)--(0.17418575103719508,0.3127497753447539);
\draw[red, thick] (0.16402452764015735,0.3738974712746203)--(0.18588103912022907,0.40276003273976374);
\draw[red, thick] (0.18588103912022907,0.40276003273976374)--(0.13638721173493218,0.6043853598410402);
\draw[red, thick] (0.0072051335010195755,0.8663122452666777)--(0.07225031953352958,0.7988984371260657);
\draw[red, thick] (0.06964041109724696,0.7722567622110954)--(0.07225031953352958,0.7988984371260657);
\draw[red, thick] (0.822841928213028,0.8123889598178787)--(0.767274141410929,0.8903491846659611);
\draw[red, thick] (0.767274141410929,0.8903491846659611)--(0.6347894819842471,0.8155039577188894);
\draw[red, thick] (0.6347894819842471,0.8155039577188894)--(0.5646470338381138,0.8144491937677845);
\draw[red, thick] (0.148397993788087,0.28793013543623097)--(0.13656917844354743,0.2001593013620731);
\draw[red, thick] (0.19799452036188725,0.2983935547149186)--(0.17418575103719508,0.3127497753447539);
\draw[red, thick] (0.19799452036188725,0.2983935547149186)--(0.39386035741641695,0.3879950389869909);
\draw[red, thick] (0.2614932165060834,0.4483621656048845)--(0.25909724262654277,0.44854516125233673);
\draw[red, thick] (0.2614932165060834,0.4483621656048845)--(0.4157501376828479,0.4287256297513061);
\draw[red, thick] (0.18588103912022907,0.40276003273976374)--(0.25909724262654277,0.44854516125233673);
\draw[red, thick] (0.4157501376828479,0.4287256297513061)--(0.39386035741641695,0.3879950389869909);
\draw[red, thick] (0.2125410814435029,0.932206125956973)--(0.1503709212537891,0.9000494683595954);
\draw[red, thick] (0.2125410814435029,0.932206125956973)--(0.049259228651907655,0.9730210070877323);
\draw[red, thick] (0.1503709212537891,0.9000494683595954)--(0.07950484344558738,0.9060196908690696);
\draw[red, thick] (0.07950484344558738,0.9060196908690696)--(0.03418643383988759,0.9411415638172064);
\draw[red, thick] (0.049259228651907655,0.9730210070877323)--(0.03418643383988759,0.9411415638172064);
\draw[red, thick] (0.3078280725794059,0.9664502126130512)--(0.2125410814435029,0.932206125956973);
\draw[red, thick] (0.0072051335010195755,0.8663122452666777)--(0.07950484344558738,0.9060196908690696);
\draw[red, thick] (0.07225031953352958,0.7988984371260657)--(0.1503709212537891,0.9000494683595954);
\draw[red, thick] (0.20785657217564413,0.7391137029871157)--(0.299736887125702,0.8166411601178638);
\draw[red, thick] (0.20785657217564413,0.7391137029871157)--(0.1626026765847196,0.6149341693854876);
\draw[red, thick] (0.299736887125702,0.8166411601178638)--(0.410824171102973,0.6379558481429204);
\draw[red, thick] (0.1626026765847196,0.6149341693854876)--(0.21771196590924116,0.5929681031496008);
\draw[red, thick] (0.21771196590924116,0.5929681031496008)--(0.3581748076390561,0.6189381158954018);
\draw[red, thick] (0.3581748076390561,0.6189381158954018)--(0.37128451309713617,0.6232642739680077);
\draw[red, thick] (0.410824171102973,0.6379558481429204)--(0.37128451309713617,0.6232642739680077);
\draw[red, thick] (0.3078280725794059,0.9664502126130512)--(0.299736887125702,0.8166411601178638);
\draw[red, thick] (0.06964041109724696,0.7722567622110954)--(0.20785657217564413,0.7391137029871157);
\draw[red, thick] (0.13638721173493218,0.6043853598410402)--(0.1626026765847196,0.6149341693854876);
\draw[red, thick] (0.5646470338381138,0.8144491937677845)--(0.4573568415412147,0.6126174690235537);
\draw[red, thick] (0.4573568415412147,0.6126174690235537)--(0.410824171102973,0.6379558481429204);
\draw[red, thick] (0.25909724262654277,0.44854516125233673)--(0.2553892027923351,0.48539597554160624);
\draw[red, thick] (0.2553892027923351,0.48539597554160624)--(0.21771196590924116,0.5929681031496008);
\draw[red, thick] (0.2553892027923351,0.48539597554160624)--(0.3581748076390561,0.6189381158954018);
\draw[red, thick] (0.2614932165060834,0.4483621656048845)--(0.37128451309713617,0.6232642739680077);
\draw[red, thick] (0.5092587354725757,0.5327225311464738)--(0.4573568415412147,0.6126174690235537);
\draw[red, thick] (0.5092587354725757,0.5327225311464738)--(0.4959932207799768,0.47088604498338843);
\draw[red, thick] (0.4959932207799768,0.47088604498338843)--(0.4157501376828479,0.4287256297513061);
\draw[red, thick] (0.45995291048752523,0.028965958412618953)--(0.33570477350580863,0.0944446314103765);
\draw[red, thick] (0.13656917844354743,0.2001593013620731)--(0.22603564326488124,0.11592720692047373);
\draw[red, thick] (0.679280231333423,0.07685494809943949)--(0.7039360788752047,0.13691025596481718);
\draw[red, thick] (0.679280231333423,0.07685494809943949)--(0.5401162699392226,0.10592861733374764);
\draw[red, thick] (0.7039360788752047,0.13691025596481718)--(0.692984094105217,0.16657177838006781);
\draw[red, thick] (0.692984094105217,0.16657177838006781)--(0.5705030851078778,0.14956730017396283);
\draw[red, thick] (0.5705030851078778,0.14956730017396283)--(0.5401162699392226,0.10592861733374764);
\draw[red, thick] (0.6704261078411188,0.2983594041110991)--(0.692984094105217,0.16657177838006781);
\draw[red, thick] (0.6704261078411188,0.2983594041110991)--(0.609626416494816,0.3156521591180685);
\draw[red, thick] (0.5705030851078778,0.14956730017396283)--(0.609626416494816,0.3156521591180685);
\draw[red, thick] (0.45995291048752523,0.028965958412618953)--(0.4765623254864447,0.08494160121449679);
\draw[red, thick] (0.4765623254864447,0.08494160121449679)--(0.5401162699392226,0.10592861733374764);
\draw[red, thick] (0.38150587398438834,0.343981855497234)--(0.4497840091409947,0.1953790520155802);
\draw[red, thick] (0.38150587398438834,0.343981855497234)--(0.29113000848283804,0.21354918994579553);
\draw[red, thick] (0.297180677272192,0.20845266990432376)--(0.2922904993448189,0.21081962561616713);
\draw[red, thick] (0.297180677272192,0.20845266990432376)--(0.4415822666940058,0.18039376697821286);
\draw[red, thick] (0.4415822666940058,0.18039376697821286)--(0.4497840091409947,0.1953790520155802);
\draw[red, thick] (0.29113000848283804,0.21354918994579553)--(0.2922904993448189,0.21081962561616713);
\draw[red, thick] (0.4959932207799768,0.47088604498338843)--(0.5773694923567108,0.33302816738559193);
\draw[red, thick] (0.39386035741641695,0.3879950389869909)--(0.38150587398438834,0.343981855497234);
\draw[red, thick] (0.5773694923567108,0.33302816738559193)--(0.4497840091409947,0.1953790520155802);
\draw[red, thick] (0.19799452036188725,0.2983935547149186)--(0.29113000848283804,0.21354918994579553);
\draw[red, thick] (0.33570477350580863,0.0944446314103765)--(0.297180677272192,0.20845266990432376);
\draw[red, thick] (0.22603564326488124,0.11592720692047373)--(0.2922904993448189,0.21081962561616713);
\draw[red, thick] (0.4765623254864447,0.08494160121449679)--(0.4415822666940058,0.18039376697821286);
\draw[red, thick] (0.609626416494816,0.3156521591180685)--(0.5773694923567108,0.33302816738559193);
\draw[red, thick] (0.6966078171943033,0.517084675521012)--(0.8394170635303543,0.5386328082928883);
\draw[red, thick] (0.6966078171943033,0.517084675521012)--(0.7629987978835259,0.4097388057906942);
\draw[red, thick] (0.8394170635303543,0.5386328082928883)--(0.8199448304662863,0.43374447248530434);
\draw[red, thick] (0.7629987978835259,0.4097388057906942)--(0.8199448304662863,0.43374447248530434);
\draw[red, thick] (0.5092587354725757,0.5327225311464738)--(0.5854609328882098,0.5563630707966651);
\draw[red, thick] (0.6704261078411188,0.2983594041110991)--(0.706255398851686,0.31142037155747987);
\draw[red, thick] (0.706255398851686,0.31142037155747987)--(0.7629987978835259,0.4097388057906942);
\draw[red, thick] (0.6966078171943033,0.517084675521012)--(0.6791059840045244,0.528146008831672);
\draw[red, thick] (0.5854609328882098,0.5563630707966651)--(0.6791059840045244,0.528146008831672);
\draw[red, thick] (0.8394170635303543,0.5386328082928883)--(0.8421316549547286,0.5424862865901321);
\draw[red, thick] (0.6791059840045244,0.528146008831672)--(0.7961709115289796,0.6183095407912882);
\draw[red, thick] (0.8421316549547286,0.5424862865901321)--(0.7961709115289796,0.6183095407912882);
\draw[red, thick] (0.8180743974620243,0.1360181565046783)--(0.7847271846057482,0.2833695475015336);
\draw[red, thick] (0.7039360788752047,0.13691025596481718)--(0.8180743974620243,0.1360181565046783);
\draw[red, thick] (0.706255398851686,0.31142037155747987)--(0.7847271846057482,0.2833695475015336);
\draw[red, thick] (0.8340383503796273,0.7758938836779715)--(0.8378339831936373,0.7542743252554434);
\draw[red, thick] (0.663654233906528,0.6697057286556454)--(0.6320097831139766,0.6120438657752406);
\draw[red, thick] (0.663654233906528,0.6697057286556454)--(0.8042869452731249,0.7046892265869095);
\draw[red, thick] (0.6320097831139766,0.6120438657752406)--(0.7954628450760505,0.6220772182985786);
\draw[red, thick] (0.8042869452731249,0.7046892265869095)--(0.7954628450760505,0.6220772182985786);
\draw[red, thick] (0.6347894819842471,0.8155039577188894)--(0.663654233906528,0.6697057286556454);
\draw[red, thick] (0.5854609328882098,0.5563630707966651)--(0.6320097831139766,0.6120438657752406);
\draw[red, thick] (0.8378339831936373,0.7542743252554434)--(0.8042869452731249,0.7046892265869095);
\draw[red, thick] (0.7961709115289796,0.6183095407912882)--(0.7954628450760505,0.6220772182985786);
\draw[black, very thick] (0.1813061501361128,0.42139697266953957)--(0.25532816389877266,0.4461882000106892);
\draw[black, very thick] (0.14872563829292249,0.5541217325950919)--(0.18134071416677522,0.6074653564194161);
\draw[black, very thick] (0.14872563829292249,0.5541217325950919)--(0.1726754312359303,0.6109192601286578);
\draw[black, very thick] (0.05789114625497006,0.7591447714783306)--(0.07059719372583992,0.7820235001551941);
\draw[black, very thick] (0.05789114625497006,0.7591447714783306)--(0.025416319820836576,0.7229033871050434);
\draw[black, very thick] (0.816871669547315,0.8207650844535032)--(0.797244196123064,0.8483019422161482);
\draw[black, very thick] (0.797244196123064,0.8483019422161482)--(0.7024556644193323,0.853730955492912);
\draw[black, very thick] (0.5150163133378375,0.8659430761214151)--(0.47164836244225017,0.9109390813944592);
\draw[black, very thick] (0.38653207447931925,0.9992506770062238)--(0.30700725570527704,0.9512527117166875);
\draw[black, very thick] (0.29224699803261267,0.34151062925235404)--(0.25616296056720894,0.325003496798576);
\draw[black, very thick] (0.3693375759616444,0.43463383746958717)--(0.37181151999371226,0.4343189103317002);
\draw[black, very thick] (0.4100600707523283,0.4294499616155625)--(0.4169904242473743,0.42937728712224627);
\draw[black, very thick] (0.25532816389877266,0.4461882000106892)--(0.2585381570265892,0.4541014013513852);
\draw[black, very thick] (0.15937905866236854,0.9047088031490487)--(0.1590067661594326,0.9455878746166037);
\draw[black, very thick] (0.15937905866236854,0.9047088031490487)--(0.17361789056520271,0.9419355930212361);
\draw[black, very thick] (0.15937905866236854,0.9047088031490487)--(0.06284152846679965,0.9189337987707326);
\draw[black, very thick] (0.1590067661594326,0.9455878746166037)--(0.17361789056520271,0.9419355930212361);
\draw[black, very thick] (0.1590067661594326,0.9455878746166037)--(0.06284152846679965,0.9189337987707326);
\draw[black, very thick] (0.2921345429365123,0.9608102969999692)--(0.30700725570527704,0.9512527117166875);
\draw[black, very thick] (0.2921345429365123,0.9608102969999692)--(0.304643617111923,0.9074897229510144);
\draw[black, very thick] (0.2526856156454257,0.7769398874936125)--(0.27892638396565794,0.7990815186159101);
\draw[black, very thick] (0.27892638396565794,0.7990815186159101)--(0.29981130217595947,0.8180189617063004);
\draw[black, very thick] (0.19387030084463816,0.7007344954508157)--(0.18563030063821997,0.7444433737589387);
\draw[black, very thick] (0.37293603650876134,0.6988993954844914)--(0.3982187058250188,0.6582318996866904);
\draw[black, very thick] (0.3982187058250188,0.6582318996866904)--(0.36683993239497326,0.6217975775006084);
\draw[black, very thick] (0.18134071416677522,0.6074653564194161)--(0.1726754312359303,0.6109192601286578);
\draw[black, very thick] (0.34509732830012835,0.6165202358288133)--(0.36683993239497326,0.6217975775006084);
\draw[black, very thick] (0.34509732830012835,0.6165202358288133)--(0.3498485323164778,0.5891158637590141);
\draw[black, very thick] (0.232149332632162,0.5956374112256204)--(0.2277040583383429,0.5644397202114);
\draw[black, very thick] (0.3010446550260852,0.8408546046323935)--(0.30055791219730454,0.8318425158098488);
\draw[black, very thick] (0.30700725570527704,0.9512527117166875)--(0.304643617111923,0.9074897229510144);
\draw[black, very thick] (0.29981130217595947,0.8180189617063004)--(0.30055791219730454,0.8318425158098488);
\draw[black, very thick] (0.23156648132676858,0.5534121318594152)--(0.2277040583383429,0.5644397202114);
\draw[black, very thick] (0.26648526291374264,0.4998123096218199)--(0.2919655586376268,0.49690587581737233);
\draw[black, very thick] (0.26648526291374264,0.4998123096218199)--(0.293302842223799,0.5346545129197736);
\draw[black, very thick] (0.293302842223799,0.5346545129197736)--(0.2919655586376268,0.49690587581737233);
\draw[black, very thick] (0.3498485323164778,0.5891158637590141)--(0.34045815137782165,0.574156595807007);
\draw[black, very thick] (0.5084688879484911,0.5290406992216823)--(0.49024128124118205,0.4678639258176544);
\draw[black, very thick] (0.38355244884048856,0.06922894316707184)--(0.37508174221498614,0.07369299906895496);
\draw[black, very thick] (0.44145691925859737,0.03871333161972445)--(0.4744519247018341,0.0778293073440714);
\draw[black, very thick] (0.37508174221498614,0.07369299906895496)--(0.35676393151346963,0.08334647127462609);
\draw[black, very thick] (0.35676393151346963,0.08334647127462609)--(0.33402090497369247,0.09314723041436684);
\draw[black, very thick] (0.27052606386950184,0.04422519848677735)--(0.25415869288763704,0.031614330708849434);
\draw[black, very thick] (0.19012226347068706,-0.017724988674975257)--(0.17858300600531926,-0.004973681536212152);
\draw[black, very thick] (0.19976901885816142,0.14065706287384386)--(0.15362496192050873,0.1841013922061019);
\draw[black, very thick] (0.19976901885816142,0.14065706287384386)--(0.153557661738525,0.18416475488949516);
\draw[black, very thick] (0.15362496192050873,0.1841013922061019)--(0.153557661738525,0.18416475488949516);
\draw[black, very thick] (0.15362496192050873,0.1841013922061019)--(0.14846174817708205,0.18896252434680957);
\draw[black, very thick] (0.153557661738525,0.18416475488949516)--(0.14846174817708205,0.18896252434680957);
\draw[black, very thick] (0.1872786858015116,0.01718136555514127)--(0.18417778144474048,0.00928081221598432);
\draw[black, very thick] (0.18417778144474048,0.00928081221598432)--(0.17858300600531926,-0.004973681536212152);
\draw[black, very thick] (0.694675762855622,0.11435450611503778)--(0.8118874384593696,0.1633564836689559);
\draw[black, very thick] (0.4744519247018341,0.0778293073440714)--(0.5263999923459467,0.10139917795079811);
\draw[black, very thick] (0.44625527759791844,0.20305910066999244)--(0.45162040934298203,0.19736030299255442);
\draw[black, very thick] (0.35675114120119855,0.3082552180370413)--(0.37932348476575223,0.34083217785400804);
\draw[black, very thick] (0.37932348476575223,0.34083217785400804)--(0.3849643126762702,0.3563026374685184);
\draw[black, very thick] (0.3781940958333602,0.1927108234591293)--(0.43025304052789864,0.18259516696622943);
\draw[black, very thick] (0.43025304052789864,0.18259516696622943)--(0.44500451138548675,0.17105528445839943);
\draw[black, very thick] (0.5646392290876473,0.35459424606590373)--(0.5570056583155196,0.3675261222104846);
\draw[black, very thick] (0.5646392290876473,0.35459424606590373)--(0.5464741713477421,0.2996958987076048);
\draw[black, very thick] (0.5570056583155196,0.3675261222104846)--(0.5326754741324017,0.40874338970209434);
\draw[black, very thick] (0.48139088714400735,0.22947900357015216)--(0.4952009553583087,0.24437837632781984);
\draw[black, very thick] (0.44500451138548675,0.17105528445839943)--(0.45877901890888895,0.13346796697448662);
\draw[black, very thick] (0.7221119059162853,0.4758477701896109)--(0.7711711320631728,0.41318386283041797);
\draw[black, very thick] (0.8272393746070328,0.4730369662060027)--(0.7711711320631728,0.41318386283041797);
\draw[black, very thick] (0.5800834392674282,0.554694787525681)--(0.6188652259059673,0.5462977123528099);
\draw[black, very thick] (0.6188652259059673,0.5462977123528099)--(0.6361928686502474,0.5410765590383083);
\draw[black, very thick] (0.7880511358690263,0.6120556809230696)--(0.7843733421016718,0.6092230402432683);
\draw[black, very thick] (0.7880511358690263,0.6120556809230696)--(0.7851149531728918,0.6097942298808102);
\draw[black, very thick] (0.7880511358690263,0.6120556809230696)--(0.7940377748297177,0.6166665966317921);
\draw[black, very thick] (0.7880511358690263,0.6120556809230696)--(0.8108312686319776,0.5941237704915503);
\draw[black, very thick] (0.7880511358690263,0.6120556809230696)--(0.7860546895085898,0.6214997109685254);
\draw[black, very thick] (0.7843733421016718,0.6092230402432683)--(0.7851149531728918,0.6097942298808102);
\draw[black, very thick] (0.7940377748297177,0.6166665966317921)--(0.7860546895085898,0.6214997109685254);
\draw[black, very thick] (0.7940377748297177,0.6166665966317921)--(0.8108312686319776,0.5941237704915503);
\draw[black, very thick] (0.8851632825205141,0.5493531197369946)--(0.9421973439307684,0.5584544111718442);
\draw[black, very thick] (1.0241664694608577,0.57153475091379)--(1.0378006206639054,0.5737104399756852);
\draw[black, very thick] (1.0823433206558173,0.580818404653238)--(1.0605948498489768,0.5980528201790608);
\draw[black, very thick] (1.0823433206558173,0.580818404653238)--(1.0588598410834567,0.5992695764466222);
\draw[black, very thick] (1.0605948498489768,0.5980528201790608)--(1.0588598410834567,0.5992695764466222);
\draw[black, very thick] (1.0605948498489768,0.5980528201790608)--(1.0418135837704467,0.6112240638453353);
\draw[black, very thick] (0.9625363197915181,0.6668209632607429)--(0.9338936187195976,0.6869080007107792);
\draw[black, very thick] (1.0418135837704467,0.6112240638453353)--(1.0588598410834567,0.5992695764466222);
\draw[black, very thick] (0.9168260577530807,0.6988774282956411)--(0.9338936187195976,0.6869080007107792);
\draw[black, very thick] (0.84050884537445,0.7523984527019172)--(0.8032960357300718,0.6954122436021932);
\draw[black, very thick] (0.642114287806097,0.6304560849618224)--(0.6273026388087081,0.6064132747783626);
\draw[black, very thick] (0.7630519194473339,0.6944316881617858)--(0.7580659591276837,0.6931913911438812);
\draw[black, very thick] (0.7630519194473339,0.6944316881617858)--(0.8032960357300718,0.6954122436021932);
\draw[black, very thick] (0.7580659591276837,0.6931913911438812)--(0.8032960357300718,0.6954122436021932);
\draw[black, very thick] (0.7097088513509334,0.6811621785700015)--(0.6984784752351365,0.6783685337927906);
\draw[black, very thick] (0.7097088513509334,0.6811621785700015)--(0.6595268144771232,0.6905536626697233);
\draw[black, very thick] (0.6984784752351365,0.6783685337927906)--(0.6595268144771232,0.6905536626697233);
\fill[blue] (0.1813061501361128,0.42139697266953957) circle (0.1pt);\draw(0.1813061501361128,0.42139697266953957) circle (0.1pt);
\fill[blue] (0.14872563829292249,0.5541217325950919) circle (0.1pt);\draw(0.14872563829292249,0.5541217325950919) circle (0.1pt);
\fill[blue] (0.05789114625497006,0.7591447714783306) circle (0.1pt);\draw(0.05789114625497006,0.7591447714783306) circle (0.1pt);
\fill[blue] (0.025416319820836576,0.7229033871050434) circle (0.1pt);\draw(0.025416319820836576,0.7229033871050434) circle (0.1pt);
\fill[blue] (0.07059719372583992,0.7820235001551941) circle (0.1pt);\draw(0.07059719372583992,0.7820235001551941) circle (0.1pt);
\fill[blue] (0.816871669547315,0.8207650844535032) circle (0.1pt);\draw(0.816871669547315,0.8207650844535032) circle (0.1pt);
\fill[blue] (0.797244196123064,0.8483019422161482) circle (0.1pt);\draw(0.797244196123064,0.8483019422161482) circle (0.1pt);
\fill[blue] (0.5150163133378375,0.8659430761214151) circle (0.1pt);\draw(0.5150163133378375,0.8659430761214151) circle (0.1pt);
\fill[blue] (0.47164836244225017,0.9109390813944592) circle (0.1pt);\draw(0.47164836244225017,0.9109390813944592) circle (0.1pt);
\fill[blue] (0.38653207447931925,0.9992506770062238) circle (0.1pt);\draw(0.38653207447931925,0.9992506770062238) circle (0.1pt);
\fill[blue] (0.7024556644193323,0.853730955492912) circle (0.1pt);\draw(0.7024556644193323,0.853730955492912) circle (0.1pt);
\fill[blue] (0.29224699803261267,0.34151062925235404) circle (0.1pt);\draw(0.29224699803261267,0.34151062925235404) circle (0.1pt);
\fill[blue] (0.25616296056720894,0.325003496798576) circle (0.1pt);\draw(0.25616296056720894,0.325003496798576) circle (0.1pt);
\fill[blue] (0.3693375759616444,0.43463383746958717) circle (0.1pt);\draw(0.3693375759616444,0.43463383746958717) circle (0.1pt);
\fill[blue] (0.37181151999371226,0.4343189103317002) circle (0.1pt);\draw(0.37181151999371226,0.4343189103317002) circle (0.1pt);
\fill[blue] (0.4100600707523283,0.4294499616155625) circle (0.1pt);\draw(0.4100600707523283,0.4294499616155625) circle (0.1pt);
\fill[blue] (0.25532816389877266,0.4461882000106892) circle (0.1pt);\draw(0.25532816389877266,0.4461882000106892) circle (0.1pt);
\fill[blue] (0.15937905866236854,0.9047088031490487) circle (0.1pt);\draw(0.15937905866236854,0.9047088031490487) circle (0.1pt);
\fill[blue] (0.1590067661594326,0.9455878746166037) circle (0.1pt);\draw(0.1590067661594326,0.9455878746166037) circle (0.1pt);
\fill[blue] (0.17361789056520271,0.9419355930212361) circle (0.1pt);\draw(0.17361789056520271,0.9419355930212361) circle (0.1pt);
\fill[blue] (0.06284152846679965,0.9189337987707326) circle (0.1pt);\draw(0.06284152846679965,0.9189337987707326) circle (0.1pt);
\fill[blue] (0.2921345429365123,0.9608102969999692) circle (0.1pt);\draw(0.2921345429365123,0.9608102969999692) circle (0.1pt);
\fill[blue] (0.2526856156454257,0.7769398874936125) circle (0.1pt);\draw(0.2526856156454257,0.7769398874936125) circle (0.1pt);
\fill[blue] (0.27892638396565794,0.7990815186159101) circle (0.1pt);\draw(0.27892638396565794,0.7990815186159101) circle (0.1pt);
\fill[blue] (0.19387030084463816,0.7007344954508157) circle (0.1pt);\draw(0.19387030084463816,0.7007344954508157) circle (0.1pt);
\fill[blue] (0.37293603650876134,0.6988993954844914) circle (0.1pt);\draw(0.37293603650876134,0.6988993954844914) circle (0.1pt);
\fill[blue] (0.3982187058250188,0.6582318996866904) circle (0.1pt);\draw(0.3982187058250188,0.6582318996866904) circle (0.1pt);
\fill[blue] (0.18134071416677522,0.6074653564194161) circle (0.1pt);\draw(0.18134071416677522,0.6074653564194161) circle (0.1pt);
\fill[blue] (0.1726754312359303,0.6109192601286578) circle (0.1pt);\draw(0.1726754312359303,0.6109192601286578) circle (0.1pt);
\fill[blue] (0.34509732830012835,0.6165202358288133) circle (0.1pt);\draw(0.34509732830012835,0.6165202358288133) circle (0.1pt);
\fill[blue] (0.232149332632162,0.5956374112256204) circle (0.1pt);\draw(0.232149332632162,0.5956374112256204) circle (0.1pt);
\fill[blue] (0.36683993239497326,0.6217975775006084) circle (0.1pt);\draw(0.36683993239497326,0.6217975775006084) circle (0.1pt);
\fill[blue] (0.3010446550260852,0.8408546046323935) circle (0.1pt);\draw(0.3010446550260852,0.8408546046323935) circle (0.1pt);
\fill[blue] (0.30700725570527704,0.9512527117166875) circle (0.1pt);\draw(0.30700725570527704,0.9512527117166875) circle (0.1pt);
\fill[blue] (0.304643617111923,0.9074897229510144) circle (0.1pt);\draw(0.304643617111923,0.9074897229510144) circle (0.1pt);
\fill[blue] (0.29981130217595947,0.8180189617063004) circle (0.1pt);\draw(0.29981130217595947,0.8180189617063004) circle (0.1pt);
\fill[blue] (0.30055791219730454,0.8318425158098488) circle (0.1pt);\draw(0.30055791219730454,0.8318425158098488) circle (0.1pt);
\fill[blue] (0.18563030063821997,0.7444433737589387) circle (0.1pt);\draw(0.18563030063821997,0.7444433737589387) circle (0.1pt);
\fill[blue] (0.5337666871494476,0.7563578295951302) circle (0.1pt);\draw(0.5337666871494476,0.7563578295951302) circle (0.1pt);
\fill[blue] (0.4828113417610102,0.6605018615470278) circle (0.1pt);\draw(0.4828113417610102,0.6605018615470278) circle (0.1pt);
\fill[blue] (0.2585381570265892,0.4541014013513852) circle (0.1pt);\draw(0.2585381570265892,0.4541014013513852) circle (0.1pt);
\fill[blue] (0.23156648132676858,0.5534121318594152) circle (0.1pt);\draw(0.23156648132676858,0.5534121318594152) circle (0.1pt);
\fill[blue] (0.2277040583383429,0.5644397202114) circle (0.1pt);\draw(0.2277040583383429,0.5644397202114) circle (0.1pt);
\fill[blue] (0.26648526291374264,0.4998123096218199) circle (0.1pt);\draw(0.26648526291374264,0.4998123096218199) circle (0.1pt);
\fill[blue] (0.293302842223799,0.5346545129197736) circle (0.1pt);\draw(0.293302842223799,0.5346545129197736) circle (0.1pt);
\fill[blue] (0.3498485323164778,0.5891158637590141) circle (0.1pt);\draw(0.3498485323164778,0.5891158637590141) circle (0.1pt);
\fill[blue] (0.2919655586376268,0.49690587581737233) circle (0.1pt);\draw(0.2919655586376268,0.49690587581737233) circle (0.1pt);
\fill[blue] (0.34045815137782165,0.574156595807007) circle (0.1pt);\draw(0.34045815137782165,0.574156595807007) circle (0.1pt);
\fill[blue] (0.5084688879484911,0.5290406992216823) circle (0.1pt);\draw(0.5084688879484911,0.5290406992216823) circle (0.1pt);
\fill[blue] (0.4169904242473743,0.42937728712224627) circle (0.1pt);\draw(0.4169904242473743,0.42937728712224627) circle (0.1pt);
\fill[blue] (0.49024128124118205,0.4678639258176544) circle (0.1pt);\draw(0.49024128124118205,0.4678639258176544) circle (0.1pt);
\fill[blue] (0.38355244884048856,0.06922894316707184) circle (0.1pt);\draw(0.38355244884048856,0.06922894316707184) circle (0.1pt);
\fill[blue] (0.44145691925859737,0.03871333161972445) circle (0.1pt);\draw(0.44145691925859737,0.03871333161972445) circle (0.1pt);
\fill[blue] (0.37508174221498614,0.07369299906895496) circle (0.1pt);\draw(0.37508174221498614,0.07369299906895496) circle (0.1pt);
\fill[blue] (0.35676393151346963,0.08334647127462609) circle (0.1pt);\draw(0.35676393151346963,0.08334647127462609) circle (0.1pt);
\fill[blue] (0.27052606386950184,0.04422519848677735) circle (0.1pt);\draw(0.27052606386950184,0.04422519848677735) circle (0.1pt);
\fill[blue] (0.33402090497369247,0.09314723041436684) circle (0.1pt);\draw(0.33402090497369247,0.09314723041436684) circle (0.1pt);
\fill[blue] (0.19012226347068706,-0.017724988674975257) circle (0.1pt);\draw(0.19012226347068706,-0.017724988674975257) circle (0.1pt);
\fill[blue] (0.25415869288763704,0.031614330708849434) circle (0.1pt);\draw(0.25415869288763704,0.031614330708849434) circle (0.1pt);
\fill[blue] (0.19976901885816142,0.14065706287384386) circle (0.1pt);\draw(0.19976901885816142,0.14065706287384386) circle (0.1pt);
\fill[blue] (0.15362496192050873,0.1841013922061019) circle (0.1pt);\draw(0.15362496192050873,0.1841013922061019) circle (0.1pt);
\fill[blue] (0.153557661738525,0.18416475488949516) circle (0.1pt);\draw(0.153557661738525,0.18416475488949516) circle (0.1pt);
\fill[blue] (0.14846174817708205,0.18896252434680957) circle (0.1pt);\draw(0.14846174817708205,0.18896252434680957) circle (0.1pt);
\fill[blue] (0.1872786858015116,0.01718136555514127) circle (0.1pt);\draw(0.1872786858015116,0.01718136555514127) circle (0.1pt);
\fill[blue] (0.18417778144474048,0.00928081221598432) circle (0.1pt);\draw(0.18417778144474048,0.00928081221598432) circle (0.1pt);
\fill[blue] (0.17858300600531926,-0.004973681536212152) circle (0.1pt);\draw(0.17858300600531926,-0.004973681536212152) circle (0.1pt);
\fill[blue] (0.694675762855622,0.11435450611503778) circle (0.1pt);\draw(0.694675762855622,0.11435450611503778) circle (0.1pt);
\fill[blue] (0.6796560855289289,0.2444362909933815) circle (0.1pt);\draw(0.6796560855289289,0.2444362909933815) circle (0.1pt);
\fill[blue] (0.4744519247018341,0.0778293073440714) circle (0.1pt);\draw(0.4744519247018341,0.0778293073440714) circle (0.1pt);
\fill[blue] (0.5263999923459467,0.10139917795079811) circle (0.1pt);\draw(0.5263999923459467,0.10139917795079811) circle (0.1pt);
\fill[blue] (0.44625527759791844,0.20305910066999244) circle (0.1pt);\draw(0.44625527759791844,0.20305910066999244) circle (0.1pt);
\fill[blue] (0.35675114120119855,0.3082552180370413) circle (0.1pt);\draw(0.35675114120119855,0.3082552180370413) circle (0.1pt);
\fill[blue] (0.37932348476575223,0.34083217785400804) circle (0.1pt);\draw(0.37932348476575223,0.34083217785400804) circle (0.1pt);
\fill[blue] (0.3781940958333602,0.1927108234591293) circle (0.1pt);\draw(0.3781940958333602,0.1927108234591293) circle (0.1pt);
\fill[blue] (0.43025304052789864,0.18259516696622943) circle (0.1pt);\draw(0.43025304052789864,0.18259516696622943) circle (0.1pt);
\fill[blue] (0.5646392290876473,0.35459424606590373) circle (0.1pt);\draw(0.5646392290876473,0.35459424606590373) circle (0.1pt);
\fill[blue] (0.5570056583155196,0.3675261222104846) circle (0.1pt);\draw(0.5570056583155196,0.3675261222104846) circle (0.1pt);
\fill[blue] (0.5326754741324017,0.40874338970209434) circle (0.1pt);\draw(0.5326754741324017,0.40874338970209434) circle (0.1pt);
\fill[blue] (0.3849643126762702,0.3563026374685184) circle (0.1pt);\draw(0.3849643126762702,0.3563026374685184) circle (0.1pt);
\fill[blue] (0.48139088714400735,0.22947900357015216) circle (0.1pt);\draw(0.48139088714400735,0.22947900357015216) circle (0.1pt);
\fill[blue] (0.4952009553583087,0.24437837632781984) circle (0.1pt);\draw(0.4952009553583087,0.24437837632781984) circle (0.1pt);
\fill[blue] (0.5464741713477421,0.2996958987076048) circle (0.1pt);\draw(0.5464741713477421,0.2996958987076048) circle (0.1pt);
\fill[blue] (0.45162040934298203,0.19736030299255442) circle (0.1pt);\draw(0.45162040934298203,0.19736030299255442) circle (0.1pt);
\fill[blue] (0.44500451138548675,0.17105528445839943) circle (0.1pt);\draw(0.44500451138548675,0.17105528445839943) circle (0.1pt);
\fill[blue] (0.45877901890888895,0.13346796697448662) circle (0.1pt);\draw(0.45877901890888895,0.13346796697448662) circle (0.1pt);
\fill[blue] (0.7221119059162853,0.4758477701896109) circle (0.1pt);\draw(0.7221119059162853,0.4758477701896109) circle (0.1pt);
\fill[blue] (0.8272393746070328,0.4730369662060027) circle (0.1pt);\draw(0.8272393746070328,0.4730369662060027) circle (0.1pt);
\fill[blue] (0.7711711320631728,0.41318386283041797) circle (0.1pt);\draw(0.7711711320631728,0.41318386283041797) circle (0.1pt);
\fill[blue] (0.5800834392674282,0.554694787525681) circle (0.1pt);\draw(0.5800834392674282,0.554694787525681) circle (0.1pt);
\fill[blue] (0.6188652259059673,0.5462977123528099) circle (0.1pt);\draw(0.6188652259059673,0.5462977123528099) circle (0.1pt);
\fill[blue] (0.6361928686502474,0.5410765590383083) circle (0.1pt);\draw(0.6361928686502474,0.5410765590383083) circle (0.1pt);
\fill[blue] (0.7880511358690263,0.6120556809230696) circle (0.1pt);\draw(0.7880511358690263,0.6120556809230696) circle (0.1pt);
\fill[blue] (0.7843733421016718,0.6092230402432683) circle (0.1pt);\draw(0.7843733421016718,0.6092230402432683) circle (0.1pt);
\fill[blue] (0.7940377748297177,0.6166665966317921) circle (0.1pt);\draw(0.7940377748297177,0.6166665966317921) circle (0.1pt);
\fill[blue] (0.7851149531728918,0.6097942298808102) circle (0.1pt);\draw(0.7851149531728918,0.6097942298808102) circle (0.1pt);
\fill[blue] (0.8108312686319776,0.5941237704915503) circle (0.1pt);\draw(0.8108312686319776,0.5941237704915503) circle (0.1pt);
\fill[blue] (0.8118874384593696,0.1633564836689559) circle (0.1pt);\draw(0.8118874384593696,0.1633564836689559) circle (0.1pt);
\fill[blue] (0.8851632825205141,0.5493531197369946) circle (0.1pt);\draw(0.8851632825205141,0.5493531197369946) circle (0.1pt);
\fill[blue] (1.0241664694608577,0.57153475091379) circle (0.1pt);\draw(1.0241664694608577,0.57153475091379) circle (0.1pt);
\fill[blue] (1.0378006206639054,0.5737104399756852) circle (0.1pt);\draw(1.0378006206639054,0.5737104399756852) circle (0.1pt);
\fill[blue] (1.0823433206558173,0.580818404653238) circle (0.1pt);\draw(1.0823433206558173,0.580818404653238) circle (0.1pt);
\fill[blue] (0.9421973439307684,0.5584544111718442) circle (0.1pt);\draw(0.9421973439307684,0.5584544111718442) circle (0.1pt);
\fill[blue] (1.0605948498489768,0.5980528201790608) circle (0.1pt);\draw(1.0605948498489768,0.5980528201790608) circle (0.1pt);
\fill[blue] (0.9625363197915181,0.6668209632607429) circle (0.1pt);\draw(0.9625363197915181,0.6668209632607429) circle (0.1pt);
\fill[blue] (1.0418135837704467,0.6112240638453353) circle (0.1pt);\draw(1.0418135837704467,0.6112240638453353) circle (0.1pt);
\fill[blue] (0.9168260577530807,0.6988774282956411) circle (0.1pt);\draw(0.9168260577530807,0.6988774282956411) circle (0.1pt);
\fill[blue] (1.0588598410834567,0.5992695764466222) circle (0.1pt);\draw(1.0588598410834567,0.5992695764466222) circle (0.1pt);
\fill[blue] (0.9338936187195976,0.6869080007107792) circle (0.1pt);\draw(0.9338936187195976,0.6869080007107792) circle (0.1pt);
\fill[blue] (0.84050884537445,0.7523984527019172) circle (0.1pt);\draw(0.84050884537445,0.7523984527019172) circle (0.1pt);
\fill[blue] (0.642114287806097,0.6304560849618224) circle (0.1pt);\draw(0.642114287806097,0.6304560849618224) circle (0.1pt);
\fill[blue] (0.7630519194473339,0.6944316881617858) circle (0.1pt);\draw(0.7630519194473339,0.6944316881617858) circle (0.1pt);
\fill[blue] (0.7580659591276837,0.6931913911438812) circle (0.1pt);\draw(0.7580659591276837,0.6931913911438812) circle (0.1pt);
\fill[blue] (0.7097088513509334,0.6811621785700015) circle (0.1pt);\draw(0.7097088513509334,0.6811621785700015) circle (0.1pt);
\fill[blue] (0.6984784752351365,0.6783685337927906) circle (0.1pt);\draw(0.6984784752351365,0.6783685337927906) circle (0.1pt);
\fill[blue] (0.714648332376192,0.6171165250868014) circle (0.1pt);\draw(0.714648332376192,0.6171165250868014) circle (0.1pt);
\fill[blue] (0.7860546895085898,0.6214997109685254) circle (0.1pt);\draw(0.7860546895085898,0.6214997109685254) circle (0.1pt);
\fill[blue] (0.8032960357300718,0.6954122436021932) circle (0.1pt);\draw(0.8032960357300718,0.6954122436021932) circle (0.1pt);
\fill[blue] (0.6595268144771232,0.6905536626697233) circle (0.1pt);\draw(0.6595268144771232,0.6905536626697233) circle (0.1pt);
\fill[blue] (0.6273026388087081,0.6064132747783626) circle (0.1pt);\draw(0.6273026388087081,0.6064132747783626) circle (0.1pt);
 \end{scope} 
 \draw (0.2,0.2) rectangle (0.8,0.8); 
 \end{tikzpicture}

%% file: VoronoiSINR5.tex
\begin{tikzpicture}[scale=8] 
 \begin{scope} 
\clip(0.2,0.2) rectangle (0.8,0.8);
\draw[red, thick] (0.8340383503796273,0.7758938836779715)--(0.822841928213028,0.8123889598178787);
\draw[red, thick] (0.16402452764015735,0.3738974712746203)--(0.17418575103719508,0.3127497753447539);
\draw[red, thick] (0.148397993788087,0.28793013543623097)--(0.17418575103719508,0.3127497753447539);
\draw[red, thick] (0.16402452764015735,0.3738974712746203)--(0.18588103912022907,0.40276003273976374);
\draw[red, thick] (0.18588103912022907,0.40276003273976374)--(0.13638721173493218,0.6043853598410402);
\draw[red, thick] (0.0072051335010195755,0.8663122452666777)--(0.07225031953352958,0.7988984371260657);
\draw[red, thick] (0.06964041109724696,0.7722567622110954)--(0.07225031953352958,0.7988984371260657);
\draw[red, thick] (0.822841928213028,0.8123889598178787)--(0.767274141410929,0.8903491846659611);
\draw[red, thick] (0.767274141410929,0.8903491846659611)--(0.6347894819842471,0.8155039577188894);
\draw[red, thick] (0.6347894819842471,0.8155039577188894)--(0.5646470338381138,0.8144491937677845);
\draw[red, thick] (0.148397993788087,0.28793013543623097)--(0.13656917844354743,0.2001593013620731);
\draw[red, thick] (0.19799452036188725,0.2983935547149186)--(0.17418575103719508,0.3127497753447539);
\draw[red, thick] (0.19799452036188725,0.2983935547149186)--(0.39386035741641695,0.3879950389869909);
\draw[red, thick] (0.2614932165060834,0.4483621656048845)--(0.25909724262654277,0.44854516125233673);
\draw[red, thick] (0.2614932165060834,0.4483621656048845)--(0.4157501376828479,0.4287256297513061);
\draw[red, thick] (0.18588103912022907,0.40276003273976374)--(0.25909724262654277,0.44854516125233673);
\draw[red, thick] (0.4157501376828479,0.4287256297513061)--(0.39386035741641695,0.3879950389869909);
\draw[red, thick] (0.2125410814435029,0.932206125956973)--(0.1503709212537891,0.9000494683595954);
\draw[red, thick] (0.2125410814435029,0.932206125956973)--(0.049259228651907655,0.9730210070877323);
\draw[red, thick] (0.1503709212537891,0.9000494683595954)--(0.07950484344558738,0.9060196908690696);
\draw[red, thick] (0.07950484344558738,0.9060196908690696)--(0.03418643383988759,0.9411415638172064);
\draw[red, thick] (0.049259228651907655,0.9730210070877323)--(0.03418643383988759,0.9411415638172064);
\draw[red, thick] (0.3078280725794059,0.9664502126130512)--(0.2125410814435029,0.932206125956973);
\draw[red, thick] (0.0072051335010195755,0.8663122452666777)--(0.07950484344558738,0.9060196908690696);
\draw[red, thick] (0.07225031953352958,0.7988984371260657)--(0.1503709212537891,0.9000494683595954);
\draw[red, thick] (0.20785657217564413,0.7391137029871157)--(0.299736887125702,0.8166411601178638);
\draw[red, thick] (0.20785657217564413,0.7391137029871157)--(0.1626026765847196,0.6149341693854876);
\draw[red, thick] (0.299736887125702,0.8166411601178638)--(0.410824171102973,0.6379558481429204);
\draw[red, thick] (0.1626026765847196,0.6149341693854876)--(0.21771196590924116,0.5929681031496008);
\draw[red, thick] (0.21771196590924116,0.5929681031496008)--(0.3581748076390561,0.6189381158954018);
\draw[red, thick] (0.3581748076390561,0.6189381158954018)--(0.37128451309713617,0.6232642739680077);
\draw[red, thick] (0.410824171102973,0.6379558481429204)--(0.37128451309713617,0.6232642739680077);
\draw[red, thick] (0.3078280725794059,0.9664502126130512)--(0.299736887125702,0.8166411601178638);
\draw[red, thick] (0.06964041109724696,0.7722567622110954)--(0.20785657217564413,0.7391137029871157);
\draw[red, thick] (0.13638721173493218,0.6043853598410402)--(0.1626026765847196,0.6149341693854876);
\draw[red, thick] (0.5646470338381138,0.8144491937677845)--(0.4573568415412147,0.6126174690235537);
\draw[red, thick] (0.4573568415412147,0.6126174690235537)--(0.410824171102973,0.6379558481429204);
\draw[red, thick] (0.25909724262654277,0.44854516125233673)--(0.2553892027923351,0.48539597554160624);
\draw[red, thick] (0.2553892027923351,0.48539597554160624)--(0.21771196590924116,0.5929681031496008);
\draw[red, thick] (0.2553892027923351,0.48539597554160624)--(0.3581748076390561,0.6189381158954018);
\draw[red, thick] (0.2614932165060834,0.4483621656048845)--(0.37128451309713617,0.6232642739680077);
\draw[red, thick] (0.5092587354725757,0.5327225311464738)--(0.4573568415412147,0.6126174690235537);
\draw[red, thick] (0.5092587354725757,0.5327225311464738)--(0.4959932207799768,0.47088604498338843);
\draw[red, thick] (0.4959932207799768,0.47088604498338843)--(0.4157501376828479,0.4287256297513061);
\draw[red, thick] (0.45995291048752523,0.028965958412618953)--(0.33570477350580863,0.0944446314103765);
\draw[red, thick] (0.13656917844354743,0.2001593013620731)--(0.22603564326488124,0.11592720692047373);
\draw[red, thick] (0.679280231333423,0.07685494809943949)--(0.7039360788752047,0.13691025596481718);
\draw[red, thick] (0.679280231333423,0.07685494809943949)--(0.5401162699392226,0.10592861733374764);
\draw[red, thick] (0.7039360788752047,0.13691025596481718)--(0.692984094105217,0.16657177838006781);
\draw[red, thick] (0.692984094105217,0.16657177838006781)--(0.5705030851078778,0.14956730017396283);
\draw[red, thick] (0.5705030851078778,0.14956730017396283)--(0.5401162699392226,0.10592861733374764);
\draw[red, thick] (0.6704261078411188,0.2983594041110991)--(0.692984094105217,0.16657177838006781);
\draw[red, thick] (0.6704261078411188,0.2983594041110991)--(0.609626416494816,0.3156521591180685);
\draw[red, thick] (0.5705030851078778,0.14956730017396283)--(0.609626416494816,0.3156521591180685);
\draw[red, thick] (0.45995291048752523,0.028965958412618953)--(0.4765623254864447,0.08494160121449679);
\draw[red, thick] (0.4765623254864447,0.08494160121449679)--(0.5401162699392226,0.10592861733374764);
\draw[red, thick] (0.38150587398438834,0.343981855497234)--(0.4497840091409947,0.1953790520155802);
\draw[red, thick] (0.38150587398438834,0.343981855497234)--(0.29113000848283804,0.21354918994579553);
\draw[red, thick] (0.297180677272192,0.20845266990432376)--(0.2922904993448189,0.21081962561616713);
\draw[red, thick] (0.297180677272192,0.20845266990432376)--(0.4415822666940058,0.18039376697821286);
\draw[red, thick] (0.4415822666940058,0.18039376697821286)--(0.4497840091409947,0.1953790520155802);
\draw[red, thick] (0.29113000848283804,0.21354918994579553)--(0.2922904993448189,0.21081962561616713);
\draw[red, thick] (0.4959932207799768,0.47088604498338843)--(0.5773694923567108,0.33302816738559193);
\draw[red, thick] (0.39386035741641695,0.3879950389869909)--(0.38150587398438834,0.343981855497234);
\draw[red, thick] (0.5773694923567108,0.33302816738559193)--(0.4497840091409947,0.1953790520155802);
\draw[red, thick] (0.19799452036188725,0.2983935547149186)--(0.29113000848283804,0.21354918994579553);
\draw[red, thick] (0.33570477350580863,0.0944446314103765)--(0.297180677272192,0.20845266990432376);
\draw[red, thick] (0.22603564326488124,0.11592720692047373)--(0.2922904993448189,0.21081962561616713);
\draw[red, thick] (0.4765623254864447,0.08494160121449679)--(0.4415822666940058,0.18039376697821286);
\draw[red, thick] (0.609626416494816,0.3156521591180685)--(0.5773694923567108,0.33302816738559193);
\draw[red, thick] (0.6966078171943033,0.517084675521012)--(0.8394170635303543,0.5386328082928883);
\draw[red, thick] (0.6966078171943033,0.517084675521012)--(0.7629987978835259,0.4097388057906942);
\draw[red, thick] (0.8394170635303543,0.5386328082928883)--(0.8199448304662863,0.43374447248530434);
\draw[red, thick] (0.7629987978835259,0.4097388057906942)--(0.8199448304662863,0.43374447248530434);
\draw[red, thick] (0.5092587354725757,0.5327225311464738)--(0.5854609328882098,0.5563630707966651);
\draw[red, thick] (0.6704261078411188,0.2983594041110991)--(0.706255398851686,0.31142037155747987);
\draw[red, thick] (0.706255398851686,0.31142037155747987)--(0.7629987978835259,0.4097388057906942);
\draw[red, thick] (0.6966078171943033,0.517084675521012)--(0.6791059840045244,0.528146008831672);
\draw[red, thick] (0.5854609328882098,0.5563630707966651)--(0.6791059840045244,0.528146008831672);
\draw[red, thick] (0.8394170635303543,0.5386328082928883)--(0.8421316549547286,0.5424862865901321);
\draw[red, thick] (0.6791059840045244,0.528146008831672)--(0.7961709115289796,0.6183095407912882);
\draw[red, thick] (0.8421316549547286,0.5424862865901321)--(0.7961709115289796,0.6183095407912882);
\draw[red, thick] (0.8180743974620243,0.1360181565046783)--(0.7847271846057482,0.2833695475015336);
\draw[red, thick] (0.7039360788752047,0.13691025596481718)--(0.8180743974620243,0.1360181565046783);
\draw[red, thick] (0.706255398851686,0.31142037155747987)--(0.7847271846057482,0.2833695475015336);
\draw[red, thick] (0.8340383503796273,0.7758938836779715)--(0.8378339831936373,0.7542743252554434);
\draw[red, thick] (0.663654233906528,0.6697057286556454)--(0.6320097831139766,0.6120438657752406);
\draw[red, thick] (0.663654233906528,0.6697057286556454)--(0.8042869452731249,0.7046892265869095);
\draw[red, thick] (0.6320097831139766,0.6120438657752406)--(0.7954628450760505,0.6220772182985786);
\draw[red, thick] (0.8042869452731249,0.7046892265869095)--(0.7954628450760505,0.6220772182985786);
\draw[red, thick] (0.6347894819842471,0.8155039577188894)--(0.663654233906528,0.6697057286556454);
\draw[red, thick] (0.5854609328882098,0.5563630707966651)--(0.6320097831139766,0.6120438657752406);
\draw[red, thick] (0.8378339831936373,0.7542743252554434)--(0.8042869452731249,0.7046892265869095);
\draw[red, thick] (0.7961709115289796,0.6183095407912882)--(0.7954628450760505,0.6220772182985786);
\draw[black, very thick] (0.1813061501361128,0.42139697266953957)--(0.25532816389877266,0.4461882000106892);
\draw[black, very thick] (0.1813061501361128,0.42139697266953957)--(0.2585381570265892,0.4541014013513852);
\draw[black, very thick] (0.14872563829292249,0.5541217325950919)--(0.18134071416677522,0.6074653564194161);
\draw[black, very thick] (0.14872563829292249,0.5541217325950919)--(0.1726754312359303,0.6109192601286578);
\draw[black, very thick] (0.05789114625497006,0.7591447714783306)--(0.07059719372583992,0.7820235001551941);
\draw[black, very thick] (0.05789114625497006,0.7591447714783306)--(0.025416319820836576,0.7229033871050434);
\draw[black, very thick] (0.025416319820836576,0.7229033871050434)--(0.07059719372583992,0.7820235001551941);
\draw[black, very thick] (0.816871669547315,0.8207650844535032)--(0.797244196123064,0.8483019422161482);
\draw[black, very thick] (0.816871669547315,0.8207650844535032)--(0.84050884537445,0.7523984527019172);
\draw[black, very thick] (0.797244196123064,0.8483019422161482)--(0.7024556644193323,0.853730955492912);
\draw[black, very thick] (0.5150163133378375,0.8659430761214151)--(0.47164836244225017,0.9109390813944592);
\draw[black, very thick] (0.38653207447931925,0.9992506770062238)--(0.2921345429365123,0.9608102969999692);
\draw[black, very thick] (0.38653207447931925,0.9992506770062238)--(0.30700725570527704,0.9512527117166875);
\draw[black, very thick] (0.29224699803261267,0.34151062925235404)--(0.25616296056720894,0.325003496798576);
\draw[black, very thick] (0.3693375759616444,0.43463383746958717)--(0.37181151999371226,0.4343189103317002);
\draw[black, very thick] (0.4100600707523283,0.4294499616155625)--(0.4169904242473743,0.42937728712224627);
\draw[black, very thick] (0.25532816389877266,0.4461882000106892)--(0.2585381570265892,0.4541014013513852);
\draw[black, very thick] (0.15937905866236854,0.9047088031490487)--(0.1590067661594326,0.9455878746166037);
\draw[black, very thick] (0.15937905866236854,0.9047088031490487)--(0.17361789056520271,0.9419355930212361);
\draw[black, very thick] (0.15937905866236854,0.9047088031490487)--(0.06284152846679965,0.9189337987707326);
\draw[black, very thick] (0.1590067661594326,0.9455878746166037)--(0.17361789056520271,0.9419355930212361);
\draw[black, very thick] (0.1590067661594326,0.9455878746166037)--(0.06284152846679965,0.9189337987707326);
\draw[black, very thick] (0.2921345429365123,0.9608102969999692)--(0.30700725570527704,0.9512527117166875);
\draw[black, very thick] (0.2921345429365123,0.9608102969999692)--(0.304643617111923,0.9074897229510144);
\draw[black, very thick] (0.2526856156454257,0.7769398874936125)--(0.27892638396565794,0.7990815186159101);
\draw[black, very thick] (0.27892638396565794,0.7990815186159101)--(0.29981130217595947,0.8180189617063004);
\draw[black, very thick] (0.19387030084463816,0.7007344954508157)--(0.18563030063821997,0.7444433737589387);
\draw[black, very thick] (0.37293603650876134,0.6988993954844914)--(0.3982187058250188,0.6582318996866904);
\draw[black, very thick] (0.3982187058250188,0.6582318996866904)--(0.36683993239497326,0.6217975775006084);
\draw[black, very thick] (0.3982187058250188,0.6582318996866904)--(0.4828113417610102,0.6605018615470278);
\draw[black, very thick] (0.18134071416677522,0.6074653564194161)--(0.1726754312359303,0.6109192601286578);
\draw[black, very thick] (0.34509732830012835,0.6165202358288133)--(0.36683993239497326,0.6217975775006084);
\draw[black, very thick] (0.34509732830012835,0.6165202358288133)--(0.3498485323164778,0.5891158637590141);
\draw[black, very thick] (0.232149332632162,0.5956374112256204)--(0.23156648132676858,0.5534121318594152);
\draw[black, very thick] (0.232149332632162,0.5956374112256204)--(0.2277040583383429,0.5644397202114);
\draw[black, very thick] (0.3010446550260852,0.8408546046323935)--(0.30055791219730454,0.8318425158098488);
\draw[black, very thick] (0.30700725570527704,0.9512527117166875)--(0.304643617111923,0.9074897229510144);
\draw[black, very thick] (0.29981130217595947,0.8180189617063004)--(0.30055791219730454,0.8318425158098488);
\draw[black, very thick] (0.23156648132676858,0.5534121318594152)--(0.2277040583383429,0.5644397202114);
\draw[black, very thick] (0.26648526291374264,0.4998123096218199)--(0.2919655586376268,0.49690587581737233);
\draw[black, very thick] (0.26648526291374264,0.4998123096218199)--(0.293302842223799,0.5346545129197736);
\draw[black, very thick] (0.293302842223799,0.5346545129197736)--(0.2919655586376268,0.49690587581737233);
\draw[black, very thick] (0.3498485323164778,0.5891158637590141)--(0.34045815137782165,0.574156595807007);
\draw[black, very thick] (0.5084688879484911,0.5290406992216823)--(0.49024128124118205,0.4678639258176544);
\draw[black, very thick] (0.38355244884048856,0.06922894316707184)--(0.37508174221498614,0.07369299906895496);
\draw[black, very thick] (0.38355244884048856,0.06922894316707184)--(0.44145691925859737,0.03871333161972445);
\draw[black, very thick] (0.44145691925859737,0.03871333161972445)--(0.4744519247018341,0.0778293073440714);
\draw[black, very thick] (0.37508174221498614,0.07369299906895496)--(0.35676393151346963,0.08334647127462609);
\draw[black, very thick] (0.35676393151346963,0.08334647127462609)--(0.33402090497369247,0.09314723041436684);
\draw[black, very thick] (0.27052606386950184,0.04422519848677735)--(0.25415869288763704,0.031614330708849434);
\draw[black, very thick] (0.19012226347068706,-0.017724988674975257)--(0.18417778144474048,0.00928081221598432);
\draw[black, very thick] (0.19012226347068706,-0.017724988674975257)--(0.17858300600531926,-0.004973681536212152);
\draw[black, very thick] (0.19976901885816142,0.14065706287384386)--(0.15362496192050873,0.1841013922061019);
\draw[black, very thick] (0.19976901885816142,0.14065706287384386)--(0.153557661738525,0.18416475488949516);
\draw[black, very thick] (0.19976901885816142,0.14065706287384386)--(0.14846174817708205,0.18896252434680957);
\draw[black, very thick] (0.15362496192050873,0.1841013922061019)--(0.153557661738525,0.18416475488949516);
\draw[black, very thick] (0.15362496192050873,0.1841013922061019)--(0.14846174817708205,0.18896252434680957);
\draw[black, very thick] (0.153557661738525,0.18416475488949516)--(0.14846174817708205,0.18896252434680957);
\draw[black, very thick] (0.1872786858015116,0.01718136555514127)--(0.18417778144474048,0.00928081221598432);
\draw[black, very thick] (0.18417778144474048,0.00928081221598432)--(0.17858300600531926,-0.004973681536212152);
\draw[black, very thick] (0.694675762855622,0.11435450611503778)--(0.8118874384593696,0.1633564836689559);
\draw[black, very thick] (0.4744519247018341,0.0778293073440714)--(0.5263999923459467,0.10139917795079811);
\draw[black, very thick] (0.4744519247018341,0.0778293073440714)--(0.45877901890888895,0.13346796697448662);
\draw[black, very thick] (0.44625527759791844,0.20305910066999244)--(0.45162040934298203,0.19736030299255442);
\draw[black, very thick] (0.44625527759791844,0.20305910066999244)--(0.43025304052789864,0.18259516696622943);
\draw[black, very thick] (0.35675114120119855,0.3082552180370413)--(0.37932348476575223,0.34083217785400804);
\draw[black, very thick] (0.35675114120119855,0.3082552180370413)--(0.3849643126762702,0.3563026374685184);
\draw[black, very thick] (0.37932348476575223,0.34083217785400804)--(0.3849643126762702,0.3563026374685184);
\draw[black, very thick] (0.3781940958333602,0.1927108234591293)--(0.43025304052789864,0.18259516696622943);
\draw[black, very thick] (0.43025304052789864,0.18259516696622943)--(0.45162040934298203,0.19736030299255442);
\draw[black, very thick] (0.43025304052789864,0.18259516696622943)--(0.44500451138548675,0.17105528445839943);
\draw[black, very thick] (0.5646392290876473,0.35459424606590373)--(0.5570056583155196,0.3675261222104846);
\draw[black, very thick] (0.5646392290876473,0.35459424606590373)--(0.5326754741324017,0.40874338970209434);
\draw[black, very thick] (0.5646392290876473,0.35459424606590373)--(0.5464741713477421,0.2996958987076048);
\draw[black, very thick] (0.5570056583155196,0.3675261222104846)--(0.5326754741324017,0.40874338970209434);
\draw[black, very thick] (0.5570056583155196,0.3675261222104846)--(0.5464741713477421,0.2996958987076048);
\draw[black, very thick] (0.48139088714400735,0.22947900357015216)--(0.4952009553583087,0.24437837632781984);
\draw[black, very thick] (0.45162040934298203,0.19736030299255442)--(0.44500451138548675,0.17105528445839943);
\draw[black, very thick] (0.44500451138548675,0.17105528445839943)--(0.45877901890888895,0.13346796697448662);
\draw[black, very thick] (0.7221119059162853,0.4758477701896109)--(0.7711711320631728,0.41318386283041797);
\draw[black, very thick] (0.8272393746070328,0.4730369662060027)--(0.7711711320631728,0.41318386283041797);
\draw[black, very thick] (0.5800834392674282,0.554694787525681)--(0.6188652259059673,0.5462977123528099);
\draw[black, very thick] (0.6188652259059673,0.5462977123528099)--(0.6361928686502474,0.5410765590383083);
\draw[black, very thick] (0.7880511358690263,0.6120556809230696)--(0.7843733421016718,0.6092230402432683);
\draw[black, very thick] (0.7880511358690263,0.6120556809230696)--(0.7851149531728918,0.6097942298808102);
\draw[black, very thick] (0.7880511358690263,0.6120556809230696)--(0.7940377748297177,0.6166665966317921);
\draw[black, very thick] (0.7880511358690263,0.6120556809230696)--(0.8108312686319776,0.5941237704915503);
\draw[black, very thick] (0.7880511358690263,0.6120556809230696)--(0.7860546895085898,0.6214997109685254);
\draw[black, very thick] (0.7843733421016718,0.6092230402432683)--(0.7851149531728918,0.6097942298808102);
\draw[black, very thick] (0.7843733421016718,0.6092230402432683)--(0.8108312686319776,0.5941237704915503);
\draw[black, very thick] (0.7940377748297177,0.6166665966317921)--(0.7860546895085898,0.6214997109685254);
\draw[black, very thick] (0.7940377748297177,0.6166665966317921)--(0.8108312686319776,0.5941237704915503);
\draw[black, very thick] (0.7851149531728918,0.6097942298808102)--(0.8108312686319776,0.5941237704915503);
\draw[black, very thick] (0.7851149531728918,0.6097942298808102)--(0.7860546895085898,0.6214997109685254);
\draw[black, very thick] (0.8851632825205141,0.5493531197369946)--(0.9421973439307684,0.5584544111718442);
\draw[black, very thick] (1.0241664694608577,0.57153475091379)--(1.0378006206639054,0.5737104399756852);
\draw[black, very thick] (1.0823433206558173,0.580818404653238)--(1.0605948498489768,0.5980528201790608);
\draw[black, very thick] (1.0823433206558173,0.580818404653238)--(1.0588598410834567,0.5992695764466222);
\draw[black, very thick] (1.0605948498489768,0.5980528201790608)--(1.0588598410834567,0.5992695764466222);
\draw[black, very thick] (1.0605948498489768,0.5980528201790608)--(1.0418135837704467,0.6112240638453353);
\draw[black, very thick] (0.9625363197915181,0.6668209632607429)--(0.9338936187195976,0.6869080007107792);
\draw[black, very thick] (1.0418135837704467,0.6112240638453353)--(1.0588598410834567,0.5992695764466222);
\draw[black, very thick] (0.9168260577530807,0.6988774282956411)--(0.9338936187195976,0.6869080007107792);
\draw[black, very thick] (0.84050884537445,0.7523984527019172)--(0.8032960357300718,0.6954122436021932);
\draw[black, very thick] (0.642114287806097,0.6304560849618224)--(0.6273026388087081,0.6064132747783626);
\draw[black, very thick] (0.7630519194473339,0.6944316881617858)--(0.7580659591276837,0.6931913911438812);
\draw[black, very thick] (0.7630519194473339,0.6944316881617858)--(0.8032960357300718,0.6954122436021932);
\draw[black, very thick] (0.7580659591276837,0.6931913911438812)--(0.8032960357300718,0.6954122436021932);
\draw[black, very thick] (0.7097088513509334,0.6811621785700015)--(0.6984784752351365,0.6783685337927906);
\draw[black, very thick] (0.7097088513509334,0.6811621785700015)--(0.6595268144771232,0.6905536626697233);
\draw[black, very thick] (0.6984784752351365,0.6783685337927906)--(0.6595268144771232,0.6905536626697233);
\fill[blue] (0.1813061501361128,0.42139697266953957) circle (0.1pt);\draw(0.1813061501361128,0.42139697266953957) circle (0.1pt);
\fill[blue] (0.14872563829292249,0.5541217325950919) circle (0.1pt);\draw(0.14872563829292249,0.5541217325950919) circle (0.1pt);
\fill[blue] (0.05789114625497006,0.7591447714783306) circle (0.1pt);\draw(0.05789114625497006,0.7591447714783306) circle (0.1pt);
\fill[blue] (0.025416319820836576,0.7229033871050434) circle (0.1pt);\draw(0.025416319820836576,0.7229033871050434) circle (0.1pt);
\fill[blue] (0.07059719372583992,0.7820235001551941) circle (0.1pt);\draw(0.07059719372583992,0.7820235001551941) circle (0.1pt);
\fill[blue] (0.816871669547315,0.8207650844535032) circle (0.1pt);\draw(0.816871669547315,0.8207650844535032) circle (0.1pt);
\fill[blue] (0.797244196123064,0.8483019422161482) circle (0.1pt);\draw(0.797244196123064,0.8483019422161482) circle (0.1pt);
\fill[blue] (0.5150163133378375,0.8659430761214151) circle (0.1pt);\draw(0.5150163133378375,0.8659430761214151) circle (0.1pt);
\fill[blue] (0.47164836244225017,0.9109390813944592) circle (0.1pt);\draw(0.47164836244225017,0.9109390813944592) circle (0.1pt);
\fill[blue] (0.38653207447931925,0.9992506770062238) circle (0.1pt);\draw(0.38653207447931925,0.9992506770062238) circle (0.1pt);
\fill[blue] (0.7024556644193323,0.853730955492912) circle (0.1pt);\draw(0.7024556644193323,0.853730955492912) circle (0.1pt);
\fill[blue] (0.29224699803261267,0.34151062925235404) circle (0.1pt);\draw(0.29224699803261267,0.34151062925235404) circle (0.1pt);
\fill[blue] (0.25616296056720894,0.325003496798576) circle (0.1pt);\draw(0.25616296056720894,0.325003496798576) circle (0.1pt);
\fill[blue] (0.3693375759616444,0.43463383746958717) circle (0.1pt);\draw(0.3693375759616444,0.43463383746958717) circle (0.1pt);
\fill[blue] (0.37181151999371226,0.4343189103317002) circle (0.1pt);\draw(0.37181151999371226,0.4343189103317002) circle (0.1pt);
\fill[blue] (0.4100600707523283,0.4294499616155625) circle (0.1pt);\draw(0.4100600707523283,0.4294499616155625) circle (0.1pt);
\fill[blue] (0.25532816389877266,0.4461882000106892) circle (0.1pt);\draw(0.25532816389877266,0.4461882000106892) circle (0.1pt);
\fill[blue] (0.15937905866236854,0.9047088031490487) circle (0.1pt);\draw(0.15937905866236854,0.9047088031490487) circle (0.1pt);
\fill[blue] (0.1590067661594326,0.9455878746166037) circle (0.1pt);\draw(0.1590067661594326,0.9455878746166037) circle (0.1pt);
\fill[blue] (0.17361789056520271,0.9419355930212361) circle (0.1pt);\draw(0.17361789056520271,0.9419355930212361) circle (0.1pt);
\fill[blue] (0.06284152846679965,0.9189337987707326) circle (0.1pt);\draw(0.06284152846679965,0.9189337987707326) circle (0.1pt);
\fill[blue] (0.2921345429365123,0.9608102969999692) circle (0.1pt);\draw(0.2921345429365123,0.9608102969999692) circle (0.1pt);
\fill[blue] (0.2526856156454257,0.7769398874936125) circle (0.1pt);\draw(0.2526856156454257,0.7769398874936125) circle (0.1pt);
\fill[blue] (0.27892638396565794,0.7990815186159101) circle (0.1pt);\draw(0.27892638396565794,0.7990815186159101) circle (0.1pt);
\fill[blue] (0.19387030084463816,0.7007344954508157) circle (0.1pt);\draw(0.19387030084463816,0.7007344954508157) circle (0.1pt);
\fill[blue] (0.37293603650876134,0.6988993954844914) circle (0.1pt);\draw(0.37293603650876134,0.6988993954844914) circle (0.1pt);
\fill[blue] (0.3982187058250188,0.6582318996866904) circle (0.1pt);\draw(0.3982187058250188,0.6582318996866904) circle (0.1pt);
\fill[blue] (0.18134071416677522,0.6074653564194161) circle (0.1pt);\draw(0.18134071416677522,0.6074653564194161) circle (0.1pt);
\fill[blue] (0.1726754312359303,0.6109192601286578) circle (0.1pt);\draw(0.1726754312359303,0.6109192601286578) circle (0.1pt);
\fill[blue] (0.34509732830012835,0.6165202358288133) circle (0.1pt);\draw(0.34509732830012835,0.6165202358288133) circle (0.1pt);
\fill[blue] (0.232149332632162,0.5956374112256204) circle (0.1pt);\draw(0.232149332632162,0.5956374112256204) circle (0.1pt);
\fill[blue] (0.36683993239497326,0.6217975775006084) circle (0.1pt);\draw(0.36683993239497326,0.6217975775006084) circle (0.1pt);
\fill[blue] (0.3010446550260852,0.8408546046323935) circle (0.1pt);\draw(0.3010446550260852,0.8408546046323935) circle (0.1pt);
\fill[blue] (0.30700725570527704,0.9512527117166875) circle (0.1pt);\draw(0.30700725570527704,0.9512527117166875) circle (0.1pt);
\fill[blue] (0.304643617111923,0.9074897229510144) circle (0.1pt);\draw(0.304643617111923,0.9074897229510144) circle (0.1pt);
\fill[blue] (0.29981130217595947,0.8180189617063004) circle (0.1pt);\draw(0.29981130217595947,0.8180189617063004) circle (0.1pt);
\fill[blue] (0.30055791219730454,0.8318425158098488) circle (0.1pt);\draw(0.30055791219730454,0.8318425158098488) circle (0.1pt);
\fill[blue] (0.18563030063821997,0.7444433737589387) circle (0.1pt);\draw(0.18563030063821997,0.7444433737589387) circle (0.1pt);
\fill[blue] (0.5337666871494476,0.7563578295951302) circle (0.1pt);\draw(0.5337666871494476,0.7563578295951302) circle (0.1pt);
\fill[blue] (0.4828113417610102,0.6605018615470278) circle (0.1pt);\draw(0.4828113417610102,0.6605018615470278) circle (0.1pt);
\fill[blue] (0.2585381570265892,0.4541014013513852) circle (0.1pt);\draw(0.2585381570265892,0.4541014013513852) circle (0.1pt);
\fill[blue] (0.23156648132676858,0.5534121318594152) circle (0.1pt);\draw(0.23156648132676858,0.5534121318594152) circle (0.1pt);
\fill[blue] (0.2277040583383429,0.5644397202114) circle (0.1pt);\draw(0.2277040583383429,0.5644397202114) circle (0.1pt);
\fill[blue] (0.26648526291374264,0.4998123096218199) circle (0.1pt);\draw(0.26648526291374264,0.4998123096218199) circle (0.1pt);
\fill[blue] (0.293302842223799,0.5346545129197736) circle (0.1pt);\draw(0.293302842223799,0.5346545129197736) circle (0.1pt);
\fill[blue] (0.3498485323164778,0.5891158637590141) circle (0.1pt);\draw(0.3498485323164778,0.5891158637590141) circle (0.1pt);
\fill[blue] (0.2919655586376268,0.49690587581737233) circle (0.1pt);\draw(0.2919655586376268,0.49690587581737233) circle (0.1pt);
\fill[blue] (0.34045815137782165,0.574156595807007) circle (0.1pt);\draw(0.34045815137782165,0.574156595807007) circle (0.1pt);
\fill[blue] (0.5084688879484911,0.5290406992216823) circle (0.1pt);\draw(0.5084688879484911,0.5290406992216823) circle (0.1pt);
\fill[blue] (0.4169904242473743,0.42937728712224627) circle (0.1pt);\draw(0.4169904242473743,0.42937728712224627) circle (0.1pt);
\fill[blue] (0.49024128124118205,0.4678639258176544) circle (0.1pt);\draw(0.49024128124118205,0.4678639258176544) circle (0.1pt);
\fill[blue] (0.38355244884048856,0.06922894316707184) circle (0.1pt);\draw(0.38355244884048856,0.06922894316707184) circle (0.1pt);
\fill[blue] (0.44145691925859737,0.03871333161972445) circle (0.1pt);\draw(0.44145691925859737,0.03871333161972445) circle (0.1pt);
\fill[blue] (0.37508174221498614,0.07369299906895496) circle (0.1pt);\draw(0.37508174221498614,0.07369299906895496) circle (0.1pt);
\fill[blue] (0.35676393151346963,0.08334647127462609) circle (0.1pt);\draw(0.35676393151346963,0.08334647127462609) circle (0.1pt);
\fill[blue] (0.27052606386950184,0.04422519848677735) circle (0.1pt);\draw(0.27052606386950184,0.04422519848677735) circle (0.1pt);
\fill[blue] (0.33402090497369247,0.09314723041436684) circle (0.1pt);\draw(0.33402090497369247,0.09314723041436684) circle (0.1pt);
\fill[blue] (0.19012226347068706,-0.017724988674975257) circle (0.1pt);\draw(0.19012226347068706,-0.017724988674975257) circle (0.1pt);
\fill[blue] (0.25415869288763704,0.031614330708849434) circle (0.1pt);\draw(0.25415869288763704,0.031614330708849434) circle (0.1pt);
\fill[blue] (0.19976901885816142,0.14065706287384386) circle (0.1pt);\draw(0.19976901885816142,0.14065706287384386) circle (0.1pt);
\fill[blue] (0.15362496192050873,0.1841013922061019) circle (0.1pt);\draw(0.15362496192050873,0.1841013922061019) circle (0.1pt);
\fill[blue] (0.153557661738525,0.18416475488949516) circle (0.1pt);\draw(0.153557661738525,0.18416475488949516) circle (0.1pt);
\fill[blue] (0.14846174817708205,0.18896252434680957) circle (0.1pt);\draw(0.14846174817708205,0.18896252434680957) circle (0.1pt);
\fill[blue] (0.1872786858015116,0.01718136555514127) circle (0.1pt);\draw(0.1872786858015116,0.01718136555514127) circle (0.1pt);
\fill[blue] (0.18417778144474048,0.00928081221598432) circle (0.1pt);\draw(0.18417778144474048,0.00928081221598432) circle (0.1pt);
\fill[blue] (0.17858300600531926,-0.004973681536212152) circle (0.1pt);\draw(0.17858300600531926,-0.004973681536212152) circle (0.1pt);
\fill[blue] (0.694675762855622,0.11435450611503778) circle (0.1pt);\draw(0.694675762855622,0.11435450611503778) circle (0.1pt);
\fill[blue] (0.6796560855289289,0.2444362909933815) circle (0.1pt);\draw(0.6796560855289289,0.2444362909933815) circle (0.1pt);
\fill[blue] (0.4744519247018341,0.0778293073440714) circle (0.1pt);\draw(0.4744519247018341,0.0778293073440714) circle (0.1pt);
\fill[blue] (0.5263999923459467,0.10139917795079811) circle (0.1pt);\draw(0.5263999923459467,0.10139917795079811) circle (0.1pt);
\fill[blue] (0.44625527759791844,0.20305910066999244) circle (0.1pt);\draw(0.44625527759791844,0.20305910066999244) circle (0.1pt);
\fill[blue] (0.35675114120119855,0.3082552180370413) circle (0.1pt);\draw(0.35675114120119855,0.3082552180370413) circle (0.1pt);
\fill[blue] (0.37932348476575223,0.34083217785400804) circle (0.1pt);\draw(0.37932348476575223,0.34083217785400804) circle (0.1pt);
\fill[blue] (0.3781940958333602,0.1927108234591293) circle (0.1pt);\draw(0.3781940958333602,0.1927108234591293) circle (0.1pt);
\fill[blue] (0.43025304052789864,0.18259516696622943) circle (0.1pt);\draw(0.43025304052789864,0.18259516696622943) circle (0.1pt);
\fill[blue] (0.5646392290876473,0.35459424606590373) circle (0.1pt);\draw(0.5646392290876473,0.35459424606590373) circle (0.1pt);
\fill[blue] (0.5570056583155196,0.3675261222104846) circle (0.1pt);\draw(0.5570056583155196,0.3675261222104846) circle (0.1pt);
\fill[blue] (0.5326754741324017,0.40874338970209434) circle (0.1pt);\draw(0.5326754741324017,0.40874338970209434) circle (0.1pt);
\fill[blue] (0.3849643126762702,0.3563026374685184) circle (0.1pt);\draw(0.3849643126762702,0.3563026374685184) circle (0.1pt);
\fill[blue] (0.48139088714400735,0.22947900357015216) circle (0.1pt);\draw(0.48139088714400735,0.22947900357015216) circle (0.1pt);
\fill[blue] (0.4952009553583087,0.24437837632781984) circle (0.1pt);\draw(0.4952009553583087,0.24437837632781984) circle (0.1pt);
\fill[blue] (0.5464741713477421,0.2996958987076048) circle (0.1pt);\draw(0.5464741713477421,0.2996958987076048) circle (0.1pt);
\fill[blue] (0.45162040934298203,0.19736030299255442) circle (0.1pt);\draw(0.45162040934298203,0.19736030299255442) circle (0.1pt);
\fill[blue] (0.44500451138548675,0.17105528445839943) circle (0.1pt);\draw(0.44500451138548675,0.17105528445839943) circle (0.1pt);
\fill[blue] (0.45877901890888895,0.13346796697448662) circle (0.1pt);\draw(0.45877901890888895,0.13346796697448662) circle (0.1pt);
\fill[blue] (0.7221119059162853,0.4758477701896109) circle (0.1pt);\draw(0.7221119059162853,0.4758477701896109) circle (0.1pt);
\fill[blue] (0.8272393746070328,0.4730369662060027) circle (0.1pt);\draw(0.8272393746070328,0.4730369662060027) circle (0.1pt);
\fill[blue] (0.7711711320631728,0.41318386283041797) circle (0.1pt);\draw(0.7711711320631728,0.41318386283041797) circle (0.1pt);
\fill[blue] (0.5800834392674282,0.554694787525681) circle (0.1pt);\draw(0.5800834392674282,0.554694787525681) circle (0.1pt);
\fill[blue] (0.6188652259059673,0.5462977123528099) circle (0.1pt);\draw(0.6188652259059673,0.5462977123528099) circle (0.1pt);
\fill[blue] (0.6361928686502474,0.5410765590383083) circle (0.1pt);\draw(0.6361928686502474,0.5410765590383083) circle (0.1pt);
\fill[blue] (0.7880511358690263,0.6120556809230696) circle (0.1pt);\draw(0.7880511358690263,0.6120556809230696) circle (0.1pt);
\fill[blue] (0.7843733421016718,0.6092230402432683) circle (0.1pt);\draw(0.7843733421016718,0.6092230402432683) circle (0.1pt);
\fill[blue] (0.7940377748297177,0.6166665966317921) circle (0.1pt);\draw(0.7940377748297177,0.6166665966317921) circle (0.1pt);
\fill[blue] (0.7851149531728918,0.6097942298808102) circle (0.1pt);\draw(0.7851149531728918,0.6097942298808102) circle (0.1pt);
\fill[blue] (0.8108312686319776,0.5941237704915503) circle (0.1pt);\draw(0.8108312686319776,0.5941237704915503) circle (0.1pt);
\fill[blue] (0.8118874384593696,0.1633564836689559) circle (0.1pt);\draw(0.8118874384593696,0.1633564836689559) circle (0.1pt);
\fill[blue] (0.8851632825205141,0.5493531197369946) circle (0.1pt);\draw(0.8851632825205141,0.5493531197369946) circle (0.1pt);
\fill[blue] (1.0241664694608577,0.57153475091379) circle (0.1pt);\draw(1.0241664694608577,0.57153475091379) circle (0.1pt);
\fill[blue] (1.0378006206639054,0.5737104399756852) circle (0.1pt);\draw(1.0378006206639054,0.5737104399756852) circle (0.1pt);
\fill[blue] (1.0823433206558173,0.580818404653238) circle (0.1pt);\draw(1.0823433206558173,0.580818404653238) circle (0.1pt);
\fill[blue] (0.9421973439307684,0.5584544111718442) circle (0.1pt);\draw(0.9421973439307684,0.5584544111718442) circle (0.1pt);
\fill[blue] (1.0605948498489768,0.5980528201790608) circle (0.1pt);\draw(1.0605948498489768,0.5980528201790608) circle (0.1pt);
\fill[blue] (0.9625363197915181,0.6668209632607429) circle (0.1pt);\draw(0.9625363197915181,0.6668209632607429) circle (0.1pt);
\fill[blue] (1.0418135837704467,0.6112240638453353) circle (0.1pt);\draw(1.0418135837704467,0.6112240638453353) circle (0.1pt);
\fill[blue] (0.9168260577530807,0.6988774282956411) circle (0.1pt);\draw(0.9168260577530807,0.6988774282956411) circle (0.1pt);
\fill[blue] (1.0588598410834567,0.5992695764466222) circle (0.1pt);\draw(1.0588598410834567,0.5992695764466222) circle (0.1pt);
\fill[blue] (0.9338936187195976,0.6869080007107792) circle (0.1pt);\draw(0.9338936187195976,0.6869080007107792) circle (0.1pt);
\fill[blue] (0.84050884537445,0.7523984527019172) circle (0.1pt);\draw(0.84050884537445,0.7523984527019172) circle (0.1pt);
\fill[blue] (0.642114287806097,0.6304560849618224) circle (0.1pt);\draw(0.642114287806097,0.6304560849618224) circle (0.1pt);
\fill[blue] (0.7630519194473339,0.6944316881617858) circle (0.1pt);\draw(0.7630519194473339,0.6944316881617858) circle (0.1pt);
\fill[blue] (0.7580659591276837,0.6931913911438812) circle (0.1pt);\draw(0.7580659591276837,0.6931913911438812) circle (0.1pt);
\fill[blue] (0.7097088513509334,0.6811621785700015) circle (0.1pt);\draw(0.7097088513509334,0.6811621785700015) circle (0.1pt);
\fill[blue] (0.6984784752351365,0.6783685337927906) circle (0.1pt);\draw(0.6984784752351365,0.6783685337927906) circle (0.1pt);
\fill[blue] (0.714648332376192,0.6171165250868014) circle (0.1pt);\draw(0.714648332376192,0.6171165250868014) circle (0.1pt);
\fill[blue] (0.7860546895085898,0.6214997109685254) circle (0.1pt);\draw(0.7860546895085898,0.6214997109685254) circle (0.1pt);
\fill[blue] (0.8032960357300718,0.6954122436021932) circle (0.1pt);\draw(0.8032960357300718,0.6954122436021932) circle (0.1pt);
\fill[blue] (0.6595268144771232,0.6905536626697233) circle (0.1pt);\draw(0.6595268144771232,0.6905536626697233) circle (0.1pt);
\fill[blue] (0.6273026388087081,0.6064132747783626) circle (0.1pt);\draw(0.6273026388087081,0.6064132747783626) circle (0.1pt);
 \end{scope} 
 \draw (0.2,0.2) rectangle (0.8,0.8); 
 \end{tikzpicture}